\documentclass[12pt,twoside]{amsart}
\usepackage{amssymb}
\usepackage{amscd}
\usepackage{xypic}

\title[Log Minimal Model 
Program]
{Fundamental theorems for the log minimal model program} 
\author{Osamu Fujino} 
\date{2010/8/13, Version 6.09}
\address{Department of Mathematics, Faculty of 
Science, Kyoto University, 
Kyoto 606-8502, Japan} 
\email{fujino@math.kyoto-u.ac.jp}
\newcommand{\Supp}[0]{{\operatorname{Supp}}}
\newcommand{\Bs}[0]{{\operatorname{Bs}}}
\newcommand{\Exc}[0]{{\operatorname{Exc}}}
\newcommand{\mult}[0]{{\operatorname{mult}}}
\newcommand{\Nlc}[0]{{\operatorname{Nlc}}}
\newcommand{\Nklt}[0]{{\operatorname{Nklt}}}
\newcommand{\Spec}[0]{{\operatorname{Spec}}}
\newcommand{\codim}[0]{{\operatorname{codim}}}
\newcommand{\Coker}[0]{{\operatorname{Coker}}}
\newcommand{\xIm}[0]{{\operatorname{Im}}}
\newcommand{\Pic}[0]{{\operatorname{Pic}}}

\newtheorem{thm}{Theorem}[section]
\newtheorem{lem}[thm]{Lemma}
\newtheorem{cor}[thm]{Corollary}
\newtheorem{prop}[thm]{Proposition}

\newtheorem{cla}{Claim}

\theoremstyle{definition}
\newtheorem{ex}[thm]{Example}
\newtheorem{defn}[thm]{Definition}
\newtheorem{rem}[thm]{Remark}
\newtheorem*{ack}{Acknowledgments}      
         
\newtheorem{say}[thm]{}
\newtheorem{step}{Step}
\begin{document}
\bibliographystyle{amsalpha+}

\maketitle

\begin{abstract}
In this paper, we prove the cone theorem and 
the contraction theorem for pairs $(X, B)$, 
where $X$ is a normal variety and $B$ is an effective 
$\mathbb R$-divisor on $X$ such that 
$K_X+B$ is $\mathbb R$-Cartier. 
\end{abstract} 

\tableofcontents

\section{Introduction}\label{sec1} 
The main purpose of this paper is 
to prove the following cone and contraction theorem. 
It is the culmination of the works of 
several authors:~Ambro, 
Benveniste, Birkar, Kawamata, Koll\'ar, 
Mori, Reid, Shokurov, and others. 
It is indispensable for the study of the log minimal model program 
for varieties with bad singularities (cf.~\cite{fujino16}). 

\begin{thm}[{cf.~Theorems \ref{thm143}, \ref{thm144}, \ref{bir-prop}, and 
\ref{thm-la}}]\label{thm1.1}
Let $X$ be a normal variety defined over $\mathbb C$ and 
let $B$ be an effective $\mathbb R$-divisor 
such that $K_X+B$ is $\mathbb R$-Cartier, and let 
$\pi:X\to S$ be a projective 
morphism onto a variety $S$. 
Then we have 
$$\overline {NE}(X/S)=\overline {NE}(X/S)_{K_X+B\geq 0} 
+\overline {NE}(X/S)_{\Nlc (X, B)}+\sum R_j$$ 
with the following properties.  
\begin{itemize}
\item[(1)] $\Nlc (X, B)$ is the non-lc locus 
of $(X, B)$ and 
$$
\overline {NE}(X/S)_{\Nlc(X, B)}=
\xIm (\overline {NE}(\Nlc(X, B)/S)
\to \overline {NE}(X/S)). 
$$
\item[(2)] 
$R_j$ is a $(K_X+B)$-negative 
extremal ray of $\overline {NE}(X/S)$ such that 
$R_j\cap \overline {NE}(X/S)_{\Nlc (X, B)}=\{0\}$ for 
every $j$. 
\item[(3)] Let $A$ be a $\pi$-ample $\mathbb R$-divisor 
on $X$. 
Then there are only finitely many $R_j$'s included in 
$(K_X+B+A)_{<0}$. In particular, 
the $R_j$'s are discrete in the half-space 
$(K_X+B)_{<0}$. 
\item[(4)] 
Let $F$ be a face of $\overline {NE}(X/S)$ such 
that $$F\cap (\overline {NE}(X/S)_{K_X+B\geq 0}+
\overline {NE}(X/S)_{\Nlc (X, B)})=\{0\}. $$ 
Then there exists a contraction morphism $\varphi_F:X\to Y$ 
over $S$. 
\begin{itemize}
\item[(i)] 
Let $C$ be an integral curve on $X$ such that 
$\pi(C)$ is a point. 
Then $\varphi_F(C)$ is a point if and 
only if 
$[C]\in F$. 
\item[(ii)] $\mathcal O_Y\simeq (\varphi_F)_*\mathcal O_X$.  
\item[(iii)] Let $L$ be a line bundle on $X$ such 
that $L\cdot C=0$ for 
every curve $C$ with $[C]\in F$. 
Then there is a line bundle $L_Y$ on $Y$ such 
that $L\simeq \varphi^*_FL_Y$. 
\end{itemize}
\item[(5)] 
Every $(K_X+B)$-negative extremal 
ray $R$ with $$R\cap \overline {NE}(X/S)_{\Nlc(X, B)}=\{0\}$$ 
is spanned by a rational 
curve $C$ with $0<-(K_X+B)\cdot C\leq 2\dim X$. 
\end{itemize}

From now on, we further assume that 
$(X, B)$ is log canonical, 
that is, $\Nlc (X, B)=\emptyset$. 
Then we have the following properties. 
\begin{itemize}
\item[(6)] 
Let $H$ be an effective $\mathbb R$-Cartier 
$\mathbb R$-divisor on $X$ such that 
$K_X+B+H$ is $\pi$-nef and $(X, B+H)$ is log canonical. 
Then, either $K_X+B$ is also $\pi$-nef or there 
is a $(K_X+B)$-negative 
extremal ray $R$ such that $(K_X+B+\lambda H)\cdot R=0$ where 
$$\lambda:=\inf\{ t\geq 0 \, |\,  K_X+B+tH \ {\text{is $\pi$-nef}} \,\}. 
$$ 
Of course, $K_X+B+\lambda H$ is $\pi$-nef.  
\end{itemize}
\end{thm}
The first half of Theorem \ref{thm1.1}, that is, 
(1), (2), (3), and (4) in Theorem \ref{thm1.1},  
is the main result of \cite{ambro}. 
His proof depends on the theory of {\em{quasi-log varieties}}. 
Unfortunately, the theory of quasi-log varieties is inaccessible 
even for experts because it requires very technical 
arguments on reducible varieties. 
In this paper, we give a proof of the above cone and contraction 
theorem 
without using the notion of quasi-log varieties. 
Our approach is much more direct than Ambro's 
theory of quasi-log varieties. 
We note that the reader does not have to refer to 
\cite{ambro} nor the book \cite{book} in order to 
read this paper. 
The latter half of Theorem \ref{thm1.1}, that is, (5) and (6) in 
Theorem \ref{thm1.1},  
is a generalization of the results obtained by 
Koll\'ar, Kawamata, Shokurov, and Birkar. 
We note that the formulation of (5) is new. 
It will play important roles in the log minimal model program 
with scaling. So, we include this 
part in our cone and contraction theorem. 

Here we would like to compare our results with the theory of quasi-log
varieties (\cite{ambro}, \cite{book}).
The first part of Theorem \ref{thm1.1}, the main theorem of our article, was first
proved
by the theory of quasi-log varieties;
the proof reduced a problem on an irreducible normal variety to
one on the union of certain proper closed subvarieties (called non-klt centers)
of various dimensions. 
The notion of quasi-log varieties was a framework to treat such reducible
varieties. 
Thus various strong vanishing theorems on quasi-log varieties were needed and
a significant part of \cite{book} was devoted to the proof of such theorems.
Delicate arguments were also needed to overcome several technical difficulties
including partial resolutions of reducible closed subvarieties.

The main idea of our paper first appeared in \cite{non-va}, which treated
a special kind of non-klt centers called minimal lc centers and proved their existence
and normality. The point of our approach was a fully general treatment of
minimal lc centers,
whose existence and normality were settled earlier under special
assumptions.

Vanishing theorems required for these proofs are stronger than
the Kawamata--Viehweg--Nadel vanishing theorem but not as difficult as
the one on quasi-log varieties.
The next step is to reduce a problem on a normal variety to one on its
minimal lc centers,
where a vanishing theorem
plays a central role. Thus it is enough to 
consider
only
normal varieties in our treatment.
Though we need to prepare vanishing theorems stronger than the
Kawamata--Viehweg--Nadel vanishing theorem, they are all proved
in our paper (without quoting \cite{ambro} or \cite{book}). 
It is our view that the
most important contribution of our paper and \cite{non-va} is the correct
formulation of various vanishing theorems and non-vanishing
theorem (Theorem \ref{thm111}), by which the cone and contraction
theorems can be proved without any difficulties. 
It is not needed
to treat reducible varieties or precise partial resolutions of
singularities of reducible varieties.
As already mentioned, the vanishing theorems needed and formulated in our
paper are stronger than the Kawamata--Viehweg--Nadel vanishing theorem.
It is our belief that this advancement of vanishing theorems distinguishes
our treatment from those in \cite{kmm}, \cite{km} and \cite{laza}.

Let us briefly recall the history of the cone and 
contraction theorem. 
In the epoch-making paper \cite{mori-th}, 
Mori invented the cone theorem for smooth 
projective varieties and 
the contraction theorem for smooth projective threefolds by 
his ingenious ideas. 
See, for example, \cite[Theorems 1.24 and 1.32]{km}. 
After Mori's pioneering works, the cone and 
contraction theorem was proved and generalized for singular varieties 
by the completely different method, which 
is now called X-method (cf.~\cite{kawa-cone}, \cite{kollar-cone}, 
\cite{reid}, and \cite{shokurov}). 
In \cite{ambro}, 
Ambro introduced the notion of quasi-log varieties and 
generalized the cone and contraction theorem. See, for example, 
\cite[Chapter 3]{book}. For the details of the history of the 
cone and contraction theorem up to \cite{kmm}, 
we recommend the reader to see the introductions of Chapters 2, 3, and 
4 of \cite{kmm}.  

We summarize the contents of this paper. 
Section \ref{sec2} is a warm-up. 
Here, we discuss the base point free theorem 
for projective log canonical surfaces to motivate 
the reader. 
This section clarifies the difference between our new approach 
and Ambro's theory of quasi-log varieties. 
In Section \ref{sec2.5}, we explain our philosophy on 
various vanishing theorems. This section helps the reader to 
understand the subsequent sections on our new vanishing theorems. 
Section \ref{sec3} collects the preliminary 
definitions and results. 
In Section \ref{sec4}, we explain 
the Hodge theoretic aspect of the injectivity theorem. 
It is an easy consequence of the 
theory of mixed Hodge 
structures on compact support cohomology groups 
of smooth quasi-projective varieties. 
Section \ref{sec5} treats generalizations 
of Koll\'ar's injectivity, torsion-free, and vanishing theorems. 
These results play crucial roles in the following sections. 
They replace the Kawamata--Viehweg vanishing theorem. 
In Section \ref{sec6}, we introduce the notion of 
non-lc ideal sheaves. 
It is an analogue of the well-known 
multiplier ideal sheaves. 
Section \ref{sec7} contains a very important 
vanishing theorem. 
It is a generalization of the Nadel vanishing theorem. 
It is very useful for the study of log canonical pairs. 
In Section \ref{sec8}, we recall the basic properties of lc centers. 
Section \ref{sec9} treats the dlt blow-up 
following Hacon and its slight refinement, 
which will be useful for future studies (cf.~\cite{gongyo}). 
Here, we need \cite{bchm}. 
In Section \ref{sec10}, we give a vanishing 
theorem for minimal lc centers. 
By the dlt blow-up obtained in Section \ref{sec9}, 
we can easily prove this very important vanishing theorem. 
Section \ref{sec11} is devoted to the proof of the non-vanishing 
theorem. 
In Section \ref{sec12}, we prove the base point free theorem. 
It is a direct consequence of the non-vanishing theorem. 
In Section \ref{sec13}, we quickly recall Shokurov's 
differents. 
Section \ref{sec14} is devoted to 
the rationality theorem. 
In Section \ref{sec15}, we obtain the cone theorem 
and contraction theorem by using 
the rationality theorem and base point free theorem. 
Section \ref{sec-n17} is a supplement to the 
base point free theorem. 
In Section \ref{sec16.5}, we discuss estimates 
of lengths of extremal rays. 
It is very important for the study of 
the log minimal model program with scaling. 
Our formulation for non-lc pairs is new. 
In Section \ref{sec16}, 
we quickly explain 
some results which were obtained by the theory of 
quasi-log varieties but can not be 
covered by our new approach. 
In the final section:~Section \ref{sec19}, 
we briefly explain some related topics obtained by the 
author for the reader's convenience. 

This paper grew out from the ideas in \cite{non-va}. 
The result in Section \ref{sec9} heavily depends on \cite{bchm}. 
We use it to prove the vanishing theorem for minimal 
lc centers in Section \ref{sec10}. We note that 
we can prove the result in Section \ref{sec10} 
without applying \cite{bchm} if we discuss the theory of mixed Hodge 
structures on compact support 
cohomology groups of {\em{reducible}} varieties. 
It was carried out in \cite[Chapter 2]{book}. 
We note that \cite[Chapter 2]{book} 
is independent of the log minimal model program 
for klt pairs. 
So, the non-vanishing theorem:~Theorem \ref{thm111}, 
the base point free theorem:~Theorem \ref{thm121}, 
the rationality theorem:~Theorem \ref{thm131}, 
and the cone theorem:~Theorem \ref{thm144} 
do not depend on the corresponding results for klt pairs. 
Therefore, our proofs are new even for klt pairs. 
In Section \ref{sec16.5}, we need Theorem \ref{thm91}, which is a consequence 
of \cite{bchm}, to prove Theorems \ref{prop146} and \ref{thm-la}. 
At present there are no proofs of Theorems \ref{prop146} and \ref{thm-la} 
without using 
\cite{bchm}. 
However, Theorem \ref{prop146} can be directly proved if 
we have an appropriate vanishing theorem 
for projective morphisms between analytic spaces. 
For the details, see \cite[Remark 3.22]{book}. 

\begin{ack}
The author 
was partially supported by The Inamori Foundation and by 
the Grant-in-Aid for Young Scientists (A) $\sharp$20684001 from 
JSPS. 
He thanks Natsuo Saito for drawing a beautiful picture 
of a Kleiman-Mori cone. 
He also thanks Takeshi Abe for useful discussions 
and Yoshinori Gongyo for some questions. 
Finally, he thanks Professor Shigefumi Mori for useful comments, information, 
and warm encouragement. 
\end{ack}

We will work over $\mathbb C$, the complex number field, throughout this 
paper. 

\section{Warm-ups}\label{sec2}

In this section, we explain the base point free theorem 
for projective log canonical surfaces to motivate 
the reader. This section clarifies the difference between our 
new approach and Ambro's theory of quasi-log 
varieties. We recommend the reader to see 
\cite[Section 4]{fuji-lec} for 
Ambro's approach. 
The following theorem is a very special case of Theorem \ref{thm121}. 

\begin{thm}[Base point free theorem for lc surfaces]\label{thm-n1}
Let $(X, B)$ be a projective log canonical surface. 
Let $L$ be a nef Cartier divisor on $X$ such that
$aL-(K_X+B)$ is ample for some $a>0$.
Then $|mL|$ is base point free for $m\gg 0$.
\end{thm}

It can not be proved by the traditional X-method. 
A key ingredient of this paper is 
the following generalization of Koll\'ar's vanishing theorem. 
We will describe it in Section \ref{sec7}. 

\begin{thm}[{cf.~Theorem \ref{thm71}}]\label{thm-n3}
Let $(X, B)$ be a projective log canonical pair. 
Let $D$ be a Cartier divisor on $X$ such that
$D-(K_X+B)$ is ample.
Let $C$ be an lc center of $(X, B)$ with
a reduced scheme structure.
Then
$$
H^i(X, \mathcal I_C\otimes \mathcal O_X(D))= 0
$$
for every $i>0$, where
$\mathcal I_C$ is the defining ideal sheaf of $C$.
In particular, the restriction map
$$
H^0(X, \mathcal O_X(D))\to H^0(C, \mathcal O_C(D))
$$
is surjective.
\end{thm}

In Theorem \ref{thm-n3}, we do not assume that 
$C$ is isolated in the non-klt locus of the pair $(X, B)$, 
neither do we assume that there exists another boundary $\mathbb R$-divisor 
$B'$ on $X$ such that $(X, B')$ is klt. 
Therefore, it can not be proved by the traditional arguments 
depending on the Kawamata--Viehweg--Nadel vanishing theorem. 

The next theorem is a special case of Theorem \ref{thm111}. 
This formulation was first introduced in \cite{non-va}. 
We will see that it is equivalent to Theorem \ref{thm-n1}. 

\begin{thm}[Non-vanishing theorem for lc surfaces]\label{thm-n2}
Let $X$ be a projective log canonical surface. 
Let $L$ be a nef Cartier divisor on $X$ such that
$aL-(K_X+B)$ is ample for some $a>0$.
Then the base locus $\Bs|mL|$ of $|mL|$ contains no lc centers of
$(X, B)$ for $m \gg 0$.
\end{thm}

\begin{proof}
It is sufficient to check that $\Bs|mL|$ contains no minimal lc 
centers
of $(X, B)$ for $m \gg 0$.
Let $C$ be a minimal lc center of $(X, B)$.
If $C$ is a point $P$, then $\Bs|mL|$ does not contain $C$ for every 
$m \geq a$.
It is because
the evaluation map
$$
H^0(X, \mathcal O_X(mL))\to \mathbb C(P)\simeq H^0(P, \mathcal O_P(mL))
$$
is surjective for every $m\geq a$ by Theorem \ref{thm-n3}.
If $C$ is a curve, then $C\subset \llcorner B\lrcorner$ and $(X, B)$ is
plt around $C$.
Therefore,
$$
K_C+B_C=(K_X+B)|_C
$$
is klt by adjunction.
Since $aL|_C-(K_C+B_C)$ is ample,
there exists $m_1$ such that
$|mL|_C|$ is base point free for every $m\geq m_1$.
By Theorem \ref{thm-n3},
the restriction map
$$
H^0(X, \mathcal O_X(mL))\to H^0(C, \mathcal O_C(mL))
$$
is surjective for every $m \geq a$.
Thus, $\Bs|mL|$ does not contain $C$ for $m \gg 0$.
So, we finish the proof since there are only finitely many minimal 
lc centers.
\end{proof}

In the above proof, $C$ is a point or a divisor on $X$. 
So, there are no difficulties in investigating 
minimal lc centers. 
When $\dim X\geq 3$, we need a more powerful 
vanishing theorem (cf.~Theorem \ref{thm101}) to study 
linear systems on minimal lc centers. 

Let us explain the proof of Theorem \ref{thm-n1}. 

\begin{proof}[Proof of {\em{Theorem \ref{thm-n1}}}]
If $(X, B)$ is klt, then the statement is well-known as 
the Kawamata--Shokurov base point free
theorem (cf.~\cite[Theorem 3.3]{km}).
So, we assume that $(X, B)$ is lc but not klt for simplicity.
By Theorem \ref{thm-n2}, we can take
general members $D_1, D_2, D_3\in |m_1L|$ for some $m_1>0$.
If $\Bs|m_1L|=\emptyset$, then $L$ is semi-ample. 
Therefore, we assume that $\Bs|m_1L|\ne \emptyset$. 
We note that
$(X, B+D)$, where $D=D_1+D_2+D_3$, is log canonical
outside $\Bs|m_1L|$ and
that $(X, B+D)$ is not log canonical at the generic point of every 
irreducible
component of $\Bs |m_1L|$. Let $c$ be the log canonical
threshold of $(X, B)$ with respect to $D$.
Then $0<c$ by Theorem \ref{thm-n2}
and $c<1$ because $(X, B+D)$ is not log canonical. 
By the above construction, $(X, B+cD)$ is log canonical and
there is an lc center $C$ of $(X, B+cD)$ such that
$C$ is contained in $\Bs|m_1L|$. By applying Theorem
\ref{thm-n2} to $$(3cm_1+a)L-(K_X+B+cD)\sim _{\mathbb Q}
aL-(K_X+B)$$ on $(X, B+cD)$,
we see that $\Bs|m_2m_1L|$ does not contain $C$ for
$m_2\gg 0$.
Therefore, $\Bs|m_2m_1L|\subsetneq \Bs|m_1L|$ holds.
By the noetherian induction, we obtain
that $L$ is semi-ample. 
With a little care, we can check that $|mL|$ is base point free 
for $m \gg 0$. We omit some details here. 
For details, see the proof of Theorem \ref{thm121}. 
\end{proof}

In Ambro's framework of quasi-log varieties (cf.~\cite{ambro}, \cite{book}, 
and \cite{fuji-lec}), 
we have to discuss the base point free theorem 
for certain reducible curves
(cf.~\cite{fuji0}) to prove Theorem \ref{thm-n1}. 
We note that the ultimate generalization of 
Theorem \ref{thm-n1} for surfaces is proved in \cite{fujino16}.

One of the main purposes of this paper 
is to generalize 
Theorem \ref{thm-n2} 
for pairs $(X, B)$, where 
$X$ is an $n$-dimensional 
normal variety and $B$ is an effective 
$\mathbb R$-divisor on $X$ such that 
$K_X+B$ is $\mathbb R$-Cartier 
(see Theorems \ref{thm111} and 
\ref{thm121}). 

\section{Kawamata--Viehweg, Nadel, Koll\'ar, $\cdots$}\label{sec2.5}

In this section, we explain our philosophy on 
vanishing theorems. 
There exists a big conceptual difference between our 
new approach described in this paper and the traditional 
arguments based on the Kawamata--Viehweg--Nadel 
vanishing theorem (cf.~\cite{kmm}, \cite{km}, 
and \cite{laza}). 

In the traditional X-method, the 
following type of the Kawamata--Viehweg vanishing theorem 
plays crucial roles (cf.~\cite[Theorem 3.1]{km}, 
\cite[Theorem 9.1.18]{laza}). 

\begin{say}[The Kawamata--Viehweg vanishing theorem]\label{v11}  
Let $X$ be a smooth projective variety and let $B$ be an 
effective $\mathbb Q$-divisor such that 
$\Supp B$ is simple normal crossing and 
$\llcorner B\lrcorner =0$. 
Let $L$ be a Cartier divisor on $X$ such that 
$L-(K_X+B)$ is nef and big. 
Then 
$$
H^i(X, \mathcal O_X(L))=0
$$ 
for every $i>0$. 
\end{say}

Recently, the (algebraic version of) Nadel vanishing theorem, which is 
a generalization of the above Kawamata--Viehweg vanishing 
theorem, is very often used for the study of linear systems 
(cf.~\cite[Theorem 9.4.17]{laza}). 
 
\begin{say}[The Nadel vanishing theorem]\label{v12} 
Let $X$ be a normal projective variety 
and let $B$ be an effective $\mathbb Q$-divisor on $X$ such that 
$K_X+B$ is $\mathbb Q$-Cartier. 
Let $L$ be a Cartier divisor on $X$ such that 
$L-(K_X+B)$ is nef and big. 
Then 
$$
H^i(X, \mathcal O_X(L)\otimes \mathcal J(X, B))=0
$$ 
for every $i>0$, 
where $\mathcal J(X, B)$ is the multiplier ideal sheaf of 
the pair $(X, B)$ (see Remark \ref{rem63} below). 
\end{say}

The following relative version of 
the Kawamata--Viehweg vanishing theorem 
sometimes plays very important roles implicitly 
(cf.~\cite[Theorem 9.4.17]{laza}, 
\cite[Corollary 2.68]{km}). 

\begin{say}[The relative Kawamata--Viehweg vanishing theorem]\label{v13}
Let $X$ be a normal projective variety and 
let $B$ be an effective $\mathbb Q$-divisor on $X$ such that 
$K_X+B$ is $\mathbb Q$-Cartier. 
Let $f:Y\to X$ be a projective resolution 
such that $K_Y+B_Y=f^*(K_X+B)$ 
and that $\Supp B_Y$ is simple normal crossing. 
Then 
$$
R^if_*\mathcal O_Y(-\llcorner B_Y\lrcorner)=0
$$ 
for every $i>0$. 
\end{say}

It is obvious that \ref{v11} is a special case of 
\ref{v12}. 
It is a routine exercise to prove \ref{v13} by \ref{v11}. 
We note that \ref{v12} can be obtained 
as a consequence of \ref{v11} and 
\ref{v13} by Hironaka's resolution theorem and 
Leray's spectral sequence. 
In this paper, we see the Nadel vanishing theorem \ref{v12}
(resp.~the relative Kawamata--Viehweg vanishing theorem \ref{v13}) 
as a special case of Koll\'ar's vanishing theorem \ref{v14} (ii) 
(resp.~Koll\'ar's torsion-free theorem \ref{v14} (i)). 

Let us recall Koll\'ar's theorems (cf.~\cite[10.15 Corollary]{kollar}). 

\begin{say}[Koll\'ar's torsion-free and vanishing theorems]\label{v14}
Let $Y$ be a smooth 
projective variety and let $\Delta$ be an effective $\mathbb Q$-divisor 
on $Y$ such that 
$\Supp \Delta$ is simple normal crossing 
and $\llcorner \Delta\lrcorner=0$. Let 
$f:Y\to X$ be a surjective 
morphism onto a projective variety $X$ and 
let $D$ be a Cartier divisor on $Y$. 
\begin{itemize}
\item[(i)] If $D-(K_Y+\Delta)\sim _{\mathbb Q}f^*M$ for some 
$\mathbb Q$-Cartier $\mathbb Q$-divisor $M$ on $X$, 
then $R^if_*\mathcal O_Y(D)$ is torsion-free for every $i\geq 0$. 
In particular, $R^if_*\mathcal O_Y(D)=0$ for every $i>0$ if $f$ is 
birational. 
\item[(ii)] If $D-(K_Y+\Delta)\sim_{\mathbb Q}f^*M$, 
where $M$ is an ample $\mathbb Q$-divisor on $X$, 
then 
$$
H^i(X, R^jf_*\mathcal O_Y(D))=0
$$ 
for every $i>0$ and $j\geq 0$. 
\end{itemize}
\end{say}

We will completely generalize it in Theorem \ref{thm53}. 
As we stated above, 
in this paper, \ref{v12} is not seen as a combination of \ref{v11} 
and \ref{v13}. It should be recognized as 
a special case of Koll\'ar's vanishing 
theorem \ref{v14} (ii). 
We do not see the vanishing theorem \ref{v13} 
as a relative vanishing theorem 
but as a special case of Koll\'ar's torsion-free theorem 
\ref{v14} (i). This change of viewpoint opens the 
door to the study of 
log canonical pairs. 

\begin{say}[Philosophy]
We note that \ref{v14} follows from the theory of pure Hodge structures. 
In our philosophy, 
we have the following correspondences. 

{\small
\begin{equation*}
 \fbox{
\begin{tabular}{c}
\begin{minipage}{4.8cm}
Kawamata log terminal pairs
\end{minipage}
\end{tabular}
}
\Longleftrightarrow 
\fbox{
 \begin{tabular}{c}
\begin{minipage}{3.8cm}
Pure Hodge structures
\end{minipage}
\end{tabular}
}
\end{equation*}
}

and 

{\small
\begin{equation*}
 \fbox{
\begin{tabular}{c}
\begin{minipage}{3.2cm}
Log canonical pairs
\end{minipage}
\end{tabular}
}
\Longleftrightarrow 
\fbox{
 \begin{tabular}{c}
\begin{minipage}{4cm}
Mixed Hodge structures
\end{minipage}
\end{tabular}
}
\end{equation*}
}

Therefore, it is very natural to prove a \lq\lq mixed\rq\rq\ 
version of \ref{v14} for the 
study of log canonical pairs. 
We will carry it out in Sections \ref{sec4} and \ref{sec5}.  
There is a big difference between our framework discussed in this 
paper (cf.~Sections \ref{sec11}, \ref{sec12}, and \ref{sec14}) 
and the traditional 
X-method from the Hodge theoretic viewpoint. 
We believe that 
all the results for klt pairs 
can be proved without using the theory of mixed Hodge 
structures (cf.~\cite{fuji-kawamata}). 
\end{say}

\begin{say}[Further discussions]
When we consider various extension theorems, 
which play crucial roles in the proof of 
the existence of pl flips (cf.~\cite{hm}), 
we think that the following correspondence is natural. 

{\small
\begin{equation*}
 \fbox{
\begin{tabular}{c}
\begin{minipage}{4.8cm}
Kawamata log terminal pairs
\end{minipage}
\end{tabular}
}
\Longleftrightarrow 
\fbox{
 \begin{tabular}{c}
\begin{minipage}{1.9cm}
$L^2$-method
\end{minipage}
\end{tabular}
}
\end{equation*}
}

The extension theorem in \cite{hm} can be proved as a consequence of 
the usual vanishing theorems. 
However, we note that the origin of the extension theorem 
is the Ohsawa--Takegoshi $L^2$ extension theorem. The Nadel 
vanishing theorem also has its origin in the $L^2$-method. 
It is very natural to try to generalize the above correspondence for log 
canonical pairs. 
However, we do not know what should be in the 
right box in the correspondence below. 

{\small
\begin{equation*}
 \fbox{
\begin{tabular}{c}
\begin{minipage}{3.2cm}
Log canonical pairs
\end{minipage}
\end{tabular}
}
\Longleftrightarrow 
\fbox{
 \begin{tabular}{c}
\begin{minipage}{1.5cm}
\begin{center}
?
\end{center}
\end{minipage}
\end{tabular}
}
\end{equation*}
}

It is very desirable to fill the right box 
correctly. 
Here, we do not discuss this topic any more. 
\end{say}

\section{Preliminaries}\label{sec3}

We work over the complex number field $\mathbb C$ throughout 
this paper. 
But we note that by using the Lefschetz principle, we can 
extend almost everything to the case where 
the base field is 
an algebraically closed field of characteristic 
zero. 
In this paper, an {\em{algebraic scheme}} denotes 
a scheme which is separated and of finite type over $\mathbb C$.  
We collect the basic notation and definitions. 

\begin{say}[$m\gg 0$]
The expression \lq ... for $m\gg 0$\rq \ means that 
\lq there exists a positive number $m_0$ such that 
... for every $m\geq m_0$.\rq
\end{say}

\begin{say}[Operations on $\mathbb R$-divisors] 
For an $\mathbb R$-Weil divisor 
$D=\sum _{j=1}^r d_j D_j$ such that 
$D_j$ is a prime divisor for 
every $j$ and $D_i\ne D_j$ for $i\ne j$, we define 
the {\em{round-up}} $\ulcorner D\urcorner =\sum _{j=1}^{r} 
\ulcorner d_j\urcorner D_j$ 
(resp.~the {\em{round-down}} $\llcorner D\lrcorner 
=\sum _{j=1}^{r} \llcorner d_j \lrcorner D_j$), 
where for every real number $x$, 
$\ulcorner x\urcorner$ (resp.~$\llcorner x\lrcorner$) is 
the integer defined by $x\leq \ulcorner x\urcorner <x+1$ 
(resp.~$x-1<\llcorner x\lrcorner \leq x$). 
The {\em{fractional part}} $\{D\}$ 
of $D$ denotes $D-\llcorner D\lrcorner$. 
We define 
\begin{align*}&
D^{=1}=\sum _{d_j=1}D_j, \ \ D^{\leq 1}=\sum_{d_j\leq 1}d_j D_j, \\ &
D^{<1}=\sum_{d_j< 1}d_j D_j, \ \ \text{and}\ \ \ 
D^{>1}=\sum_{d_j>1}d_j D_j. 
\end{align*}
We call $D$ a {\em{boundary}} 
$\mathbb R$-divisor if 
$0\leq d_j\leq 1$ 
for every $j$. 
We note that $\sim _{\mathbb Q}$ (resp.~$\sim _{\mathbb R}$) 
denotes the $\mathbb Q$-linear (resp.~$\mathbb R$-linear) equivalence 
of $\mathbb Q$-divisors (resp.~$\mathbb R$-divisors). 
Let $D_1$ and $D_2$ be $\mathbb R$-Cartier 
$\mathbb R$-divisors on 
$X$ and 
let $f:X\to Y$ be a morphism. 
We say that $D_1$ and $D_2$ are {\em{$\mathbb R$-linearly 
$f$-equivalent}}, denoted by $D_1\sim _{\mathbb R, f}D_2$, 
if and only if there is an 
$\mathbb R$-Cartier $\mathbb R$-divisor $B$ 
on $Y$ such that $D_1\sim _{\mathbb R}D_2+f^*B$.  
We can define $D_1\sim _{\mathbb Q, f}D_2$ for $\mathbb 
Q$-Cartier $\mathbb Q$-divisors $D_1$ and $D_2$ similarly. 
\end{say} 

\begin{defn}[Exceptional locus]
For a proper birational morphism $f:X\to Y$, 
the {\em{exceptional locus}} $\Exc (f)\subset X$ is the locus where 
$f$ is not an isomorphism. 
\end{defn}

\begin{say}[Discrepancy, singularities of pairs, etc.]
Let $X$ be a normal variety and let $B$ be an effective $\mathbb R$-divisor 
on $X$ such that $K_X+B$ is $\mathbb R$-Cartier. 
Let $f:Y\to X$ be a resolution such that 
$\Exc (f)\cup f^{-1}_*B$ has a simple normal crossing 
support, where $f^{-1}_*B$ is the strict transform of $B$ on $Y$. 
We write $$K_Y=f^*(K_X+B)+\sum _i a_i E_i$$ and 
$a(E_i, X, B)=a_i$. 
We say that $(X, B)$ is {\em{lc}} (resp.~{\em{klt}}) if and only if 
$a_i\geq -1$ (resp.~$a_i>-1$) for every $i$. 
Note that the {\em{discrepancy}} $a(E, X, B)\in \mathbb R$ can be 
defined for every prime divisor $E$ {\em{over}} $X$. 
If $a(E, X, B)>-1$ for every 
exceptional divisor $E$ over $X$, 
then the pair $(X, B)$ is called {\em{plt}}. 
Here, lc (resp.~klt, plt) is an abbreviation of 
{\em{log canonical}} (resp.~{\em{Kawamata log terminal}}, 
{\em{purely log terminal}}). 
By the definition, there exists the largest Zariski open set 
$U$ (resp.~$U'$) of $X$ such that $(X, B)$ 
is lc (resp.~klt) on $U$ (resp.~$U'$). 
We put $\Nlc(X, B)=X\setminus U$ (resp.~$\Nklt (X, B)=X\setminus U'$) 
and call it 
the {\em{non-lc locus}} (resp.~{\em{non-klt locus}}) of the pair 
$(X, B)$. 
We sometimes simply denote $\Nlc (X, B)$ by $X_{NLC}$. 

Let $(X, B)$ be a log canonical pair and let $M$ be an effective 
$\mathbb R$-Cartier $\mathbb R$-divisor 
on $X$. 
The {\em{log canonical threshold}} of $(X, B)$ with respect to 
$M$ is defined by 
$$
c=\sup \{ t \in \mathbb R\, | \, (X, B+tM) \, \text{is log canonical}\}. 
$$
\end{say}

\begin{defn}[Center] 
Let $E$ be a prime divisor over $X$. The closure 
of the image of $E$ on $X$ is denoted by 
$c_X(E)$ and 
called the {\em{center}} of $E$ on $X$. 
\end{defn}

\begin{defn}[Lc center]
Let $X$ be a normal variety and let $B$ be an effective 
$\mathbb R$-divisor on $X$ such that 
$K_X+B$ is $\mathbb R$-Cartier. 
If $a(E, X, B)=-1$ and 
$c_X(E)$ is not contained in $\Nlc(X, B)$, then 
$c_X(E)$ is called 
an {\em{lc center}} of $(X, B)$. 
It is obvious that there are at most finitely many lc centers. 
\end{defn} 

We note that 
our definition 
of lc centers is slightly different from the usual one. 
For details, see \cite[Section 3]{ft}. 

\begin{defn}[Stratum]Let $(X, B)$ be 
a log canonical pair. 
A {\em{stratum}} of $(X, B)$ denotes 
$X$ itself or an lc center of $(X, B)$. 

Let $T$ be a simple normal crossing divisor on a smooth 
variety $Y$. A {\em{stratum}} of $T$ denotes 
a stratum of the pair $(Y, T)$ contained in $T$. 
\end{defn}

\begin{say}[Kleiman--Mori cone] 
Let $X$ be an algebraic scheme over $\mathbb C$ and 
let $\pi:X\to S$ be a proper morphism 
to an algebraic scheme $S$. 
Let $\Pic (X)$ be the group of line bundles on $X$. 
Take a complete curve on $X$ which 
is mapped to a point by $\pi$. For 
$\mathcal L\in \Pic (X)$, we 
define the intersection number 
$\mathcal L\cdot C=\deg _{\overline C}f^*\mathcal L$, 
where $f:\overline C\to C$ is the normalization of $C$. 
Via this intersection pairing, we introduce 
a bilinear form 
$$
\cdot: \Pic (X)\times Z_1(X/S)\to \mathbb Z, 
$$
where $Z_1(X/S)$ is the free abelian group 
generated by integral curves which are mapped to 
points on $S$ by $\pi$. 

Now we have the notion of numerical equivalence both 
in $Z_1(X/S)$ and in $\Pic(X)$, which 
is denoted by $\equiv$, 
and 
we obtain a perfect pairing 
$$
N^1(X/S)\times N_1(X/S)\to \mathbb R,  
$$
where 
$$N^1(X/S)=\{\Pic (X)/\equiv\}\otimes \mathbb R \ \ \  \text{and} 
\ \ \ N_1(X/S)=\{Z_1(X/S)/\equiv\}\otimes \mathbb R, 
$$ 
namely $N^1(X/S)$ and $N_1(X/S)$ are dual to 
each other through this intersection pairing. 
It is well known that $$\dim _{\mathbb R}N^1(X/S)=
\dim _{\mathbb R}N_1(X/S)<\infty.$$  
We write $\rho (X/S)=\dim _{\mathbb R}N^1(X/S)=
\dim _{\mathbb R}N_1(X/S)$. 
We define the Kleiman--Mori cone $\overline {NE}(X/S)$ 
as the closed convex cone 
in $N_1(X/S)$ generated by integral curves 
on $X$ which are mapped to points on $S$ by $\pi$. 
When $S=\Spec \mathbb C$, we drop $/\Spec \mathbb C$ from 
the notation, e.g., we simply write $N_1(X)$ instead 
of $N_1(X/\Spec \mathbb C)$. 
\end{say}

\begin{defn}
An element $D\in N^1(X/S)$ is called {\em{$\pi$-nef}} 
(or {\em{relatively nef for $\pi$}}), if $D\geq 0$ 
on $\overline {NE}(X/S)$. When $S=\Spec \mathbb C$, 
we simply say that $D$ is {\em{nef}}. 
\end{defn}

\begin{thm}[Kleiman's criterion for ampleness]\label{klei-th}  
Let $\pi:X\to S$ be a {\em{projective}} morphism 
between algebraic schemes. 
Then $\mathcal L\in \Pic (X)$ is $\pi$-ample 
if and only if the numerical class of $\mathcal L$ in $N^1(X/S)$ 
gives a positive function on $\overline {NE}(X/S)\setminus 
\{0\}$. 
\end{thm}

In Theorem \ref{klei-th}, 
we note that the projectivity of $\pi$ is indispensable (cf.~\cite{fuji-klei}). 

\begin{defn}[Semi-ample $\mathbb R$-divisors]\label{defn4949} 
An $\mathbb R$-Cartier $\mathbb R$-divisor 
$D$ on $X$ is {\em{$\pi$-semi-ample}} 
if $D\sim _{\mathbb R}\sum _ia_i D_i$, 
where $D_i$ is a $\pi$-semi-ample Cartier 
divisor 
on $X$ and $a_i$ is a positive 
real number for every $i$.  
\end{defn}

\begin{rem}
In Definition \ref{defn4949}, we can replace $D\sim _{\mathbb R}
\sum _i a_i D_i$ with $D=\sum _i a_i D_i$ since every principal 
Cartier divisor on $X$ is $\pi$-semi-ample. 
\end{rem}

The following two lemmas seem to be missing in the literature. 

\begin{lem}\label{49-1} 
Let $D$ be an $\mathbb R$-Cartier 
$\mathbb R$-divisor on $X$. 
Then the following conditions are equivalent. 
\begin{itemize}
\item[(1)] $D$ is $\pi$-semi-ample. 
\item[(2)] There exists a morphism $f:X\to Y$ over $S$ such that 
$D\sim _{\mathbb R}f^*A$, where 
$A$ is an $\mathbb R$-Cartier $\mathbb R$-divisor on $Y$ which is 
ample over $S$. 
\end{itemize}
\end{lem}
\begin{proof}
It is obvious that (1) follows from (2). 
If $D$ is $\pi$-semi-ample, then we can 
write $D\sim _{\mathbb R}\sum _i a_i D_i$ as in Definition \ref{defn4949}. 
By replacing $D_i$ with its multiple, we can assume that 
$\pi^*\pi_*\mathcal O_X(D_i)\to \mathcal O_X(D_i)$ is surjective 
for every $i$. Let $f:X\to Y$ be a morphism over $S$ 
obtained by the surjection $\pi^*\pi_*\mathcal O_X(\sum _i D_i)\to 
\mathcal O_X(\sum _i D_i)$. 
Then it is easy to see that $f:Y\to X$ has the desired property. 
\end{proof}

\begin{lem}\label{49-2}
Let $D$ be a Cartier divisor on $X$. 
If $D$ is $\pi$-semi-ample in the sense of {\em{Definition \ref{defn4949}}}, 
then $D$ is $\pi$-semi-ample 
in the usual sense, that is, $\pi^*\pi_*\mathcal O_X(mD)\to 
\mathcal O_X(mD)$ is surjective for some positive integer $m$. 
In particular, {\em{Definition \ref{defn4949}}} is well-defined. 
\end{lem}
\begin{proof}
We write $D\sim _{\mathbb R}\sum _i a_iD_i$ as in Definition \ref{defn4949}. 
Let $f:X\to Y$ be a morphism in Lemma \ref{49-1} (2). 
By taking the Stein factorization, we can assume that 
$f$ has connected fibers. By the construction, $D_i\sim _{\mathbb Q, f}0$ for every $i$. 
By replacing $D_i$ with its multiple, we can assume that 
$D_i\sim f^*D'_i$ for some Cartier divisor 
$D'_i$ on $Y$ for every $i$. 
Let $U$ be any Zariski open set of $Y$ on which $D'_i\sim 0$ for every 
$i$. On $f^{-1}(U)$, we have $D\sim _{\mathbb R}0$. 
This implies $D\sim _{\mathbb Q}0$ on $f^{-1}(U)$ since $D$ is 
Cartier. 
Therefore, there exists a positive integer $m$ such that 
$f^*f_*\mathcal O_X(mD)\to \mathcal O_X(mD)$ is surjective. 
By this surjection, we have $mD\sim f^*A$ for a Cartier divisor 
$A$ on $Y$ which is ample over $S$. 
This means that $D$ is $\pi$-semi-ample in the usual sense. 
\end{proof}

We will repeatedly use the following easy lemma. 
We give a detailed proof for the reader's convenience. 
 
\begin{lem}\label{lem21}
Let $X$ be a normal variety and 
let $B$ be an effective $\mathbb R$-Cartier $\mathbb R$-divisor on $X$ such 
that $\llcorner B\lrcorner =0$. Let $A$ be a Cartier 
divisor on $X$. 
Assume that $A\sim _{\mathbb R}B$. 
Then there exists a $\mathbb Q$-Cartier $\mathbb Q$-divisor 
$C$ on $X$ such that 
$A\sim _{\mathbb Q}C$, $\llcorner C\lrcorner=0$, 
and $\Supp C=\Supp B$. 
\end{lem}
\begin{proof}
We can write $B=A+\sum _{i=1}^{k} r_i (f_i)$, where 
$r_i\in \mathbb R$ and $f_i$ is a rational function on $X$ for every $i$. 
We put $$E=\Supp A\cup \Supp B\cup \underset{i=1}{\overset{k}{\bigcup}}\Supp (f_i). $$
Let $E=\sum _{j=1}^{n}E_j$ be the irreducible 
decomposition of $E$. 
We can write 
$$A=\sum _j a_j E_j,\ \   B=\sum _j b_j E_j, $$
and 
$$(f_i)=\sum _j m_{ij}E_j \ \ \text{for every} \ \ i. $$
We can assume that 
$b_j \in \mathbb Q$ for $1\leq j\leq l$ and 
$b_j \not \in \mathbb Q$ for $l+1\leq j\leq n$. 
We note that $a_j \in \mathbb Z$ for every $j$ and 
that $m_{ij}\in \mathbb Z$ for every $i, j$. 
We define 
\begin{align*}
\mathcal S=\left\{(v_1, \cdots, v_k) \in \mathbb R^k\ ; \ 
b_j =a_j +\sum _{i=1}^kv_i m_{ij} \ \text{for}\  1\leq j\leq l \right\}. 
\end{align*}
Then $\mathcal S$ is an affine subspace of $\mathbb R^k$ defined 
over $\mathbb Q$. We note that 
$\mathcal S$ is not empty since $(r_1, \cdots, r_k)\in \mathcal S$. 
If we take $(r'_1, \cdots, r'_k)\in \mathcal S\cap \mathbb Q^k$ which 
is very close to $(r_1, \cdots, r_k)$ and 
put $C=A+\sum _i r'_i (f_i)$, 
then it is obvious that $C$ satisfies the desired properties. 
\end{proof}

The next lemma is well known as the negativity lemma. 

\begin{lem}[Negativity lemma]\label{lem211}
Let $h:Z\to Y$ be a proper birational morphism between 
normal varieties. Let $-B$ be an $h$-nef 
$\mathbb R$-Cartier $\mathbb R$-divisor 
on $Z$. 
Then we have the following statements. 
\begin{itemize}
\item[(1)] $B$ is effective if and only if 
$h_*B$ is. 
\item[(2)] Assume that $B$ is 
effective. 
Then for every $y\in Y$, either 
$h^{-1}(y)\subset \Supp B$ or 
$h^{-1}(y)\cap \Supp B=\emptyset$. 
\end{itemize}
\end{lem}

Lemma \ref{lem211} is essentially an application 
of the Hodge index theorem for smooth 
projective surfaces. 
For the proof, 
see, for example, \cite[Lemma 3.39]{km}. 

We close this section with the following useful lemma. 
It is a consequence of Szab\'o's resolution lemma. 

\begin{lem}\label{lem25} 
Let $Z$ be a smooth variety and let 
$B$ be an $\mathbb R$-divisor 
on $Z$ such that 
$\Supp B$ is simple normal crossing. 
Let $f:Z\to X$ be a projective morphism and let $\overline X$ be a projective 
variety such that 
$\overline X$ contains $X$ as a Zariski open set. 
Then there exist a smooth projective 
variety $\overline Z$ and an $\mathbb R$-divisor $\overline B$ on 
$\overline Z$ such that 
\begin{itemize}
\item[(i)] $f:Z\to X$ is extended to $\overline f: \overline Z\to 
\overline X$. 
\item[(ii)]$\Supp \overline B$ is simple normal crossing. 
\item[(iii)] $\Supp \overline B\cup \Supp (\overline Z\setminus Z)$ is 
simple normal crossing. 
\item[(iv)] $\overline B|_Z=B$. 
\end{itemize}
\end{lem}
\begin{proof}
Let $Z'$ be an arbitrary compactification of 
$Z$. 
By taking the graph of $f:Z'\dashrightarrow \overline X$ and 
using Hironaka's resolution, we can assume that 
$Z'$ is smooth projective, 
$\Supp (Z'\setminus Z)$ is simple normal crossing, and $f:Z\to X$ is 
extended to $f':Z'\to \overline X$. 
Let $B'$ be the closure of $B$ on $Z'$. 
We apply Szab\'o's resolution lemma (see, 
for example, \cite{what}) to 
$\Supp B'\cup \Supp (Z'\setminus Z)$. 
Then we obtain the desired variety $\overline Z$ and 
$\overline B$. By the above construction, $f$ can be extended 
to $\overline f:\overline Z\to \overline X$. 
\end{proof}

\section{Hodge theoretic injectivity theorem}\label{sec4}

In this section, we will 
prove the following injectivity theorem, 
which is a generalization of \cite[5.1.~b)]{ev} 
for $\mathbb R$-divisors. We use the classical topology throughout 
this section. 

\begin{prop}[Fundamental injectivity theorem]\label{prop41} 
Let $X$ be a smooth projective variety and let $S+B$ be 
a boundary $\mathbb R$-divisor on $X$ such that 
the support of $S+B$ is simple normal crossing 
and $\llcorner S+B\lrcorner=S$. 
Let $L$ be a Cartier divisor on $X$ and let $D$ be an 
effective Cartier divisor whose support is contained 
in $\Supp B$. 
Assume that 
$L\sim_{\mathbb R}K_X+S+B$. Then 
the natural homomorphisms 
$$
H^q(X, \mathcal O_X(L))\to H^q(X, \mathcal O_X(L+D)) 
$$ 
which are induced by the natural inclusion 
$\mathcal O_X\to \mathcal O_X(D)$ are injective for all $q$. 
\end{prop}

Let us recall 
some results on the theory of mixed Hodge structures. 

\begin{say}[Mixed Hodge structures]\label{say52}
Let $V$ be a smooth 
projective variety and 
$\Sigma$ a simple normal crossing 
divisor on $V$. 
Let $\iota:V\setminus \Sigma\to V$ be the natural open 
immersion. 
Then $\iota_{!}\mathbb C_{V\setminus \Sigma}$ is quasi-isomorphic 
to the complex $\Omega^{\bullet}_V(\log \Sigma)\otimes \mathcal 
O_V(-\Sigma)$. 
By this quasi-isomorphism, 
we can construct the following 
spectral sequence 
\begin{align*}
E^{pq}_1=H^q(V, \Omega^{p}_V(\log \Sigma)
\otimes \mathcal O_V(-\Sigma))\Rightarrow H^{p+q}_c(V\setminus \Sigma, 
\mathbb C). 
\end{align*}
By the Serre duality, 
the right hand side 
$$H^q(V, \Omega^{p}_V(\log \Sigma)\otimes \mathcal O_V(-\Sigma))$$ is dual to 
$$H^{n-q}(V, \Omega^{n-p}_V(\log \Sigma)),$$ where 
$n=\dim V$. 
By the Poincar\'e duality, 
$H^{p+q}_c(V\setminus \Sigma, 
\mathbb C)$ is dual to 
$H^{2n-(p+q)}(V\setminus \Sigma, 
\mathbb C)$. 
Therefore, 
\begin{align*}
\dim H^{k}_c(V\setminus \Sigma, 
\mathbb C)=\sum _{p+q=k}\dim 
H^q(V, \Omega^{p}_V(\log \Sigma)\otimes 
\mathcal O_V(-\Sigma))
\end{align*} 
by Deligne (cf.~\cite[Corollaire (3.2.13) (ii)]{deligne}). 
Thus, the above spectral sequence 
degenerates 
at $E_1$. 
We will use this $E_1$-degeneration 
in the proof of Proposition \ref{prop41}. 
By the above $E_1$-degeneration, we obtain 
$$
H^k_c(V\setminus \Sigma, \mathbb C)\simeq 
\bigoplus_{p+q=k} H^q(V, \Omega^{p}_V(\log \Sigma)\otimes \mathcal O_V(-\Sigma)). 
$$ 
In particular, the natural 
inclusion $\iota_{!}\mathbb C_{V\setminus \Sigma}
\subset \mathcal O_V(-\Sigma)$ induces 
surjections 
$$
H^p_c(V\setminus \Sigma, \mathbb C)\simeq H^p(V, \iota_!
\mathbb C_{V\setminus \Sigma})\to H^p(V, \mathcal O_V(-\Sigma)) 
$$ 
for all $p$. 
\end{say}

\begin{proof}[Proof of {\em{Proposition \ref{prop41}}}] 
By Lemma \ref{lem21}, we can assume that 
$B$ is a $\mathbb Q$-divisor and 
that $L\sim _{\mathbb Q}K_X+S+B$. 
We put $\mathcal L=\mathcal O_X(L-K_X-S)$. 
Let $\nu$ be the smallest positive integer 
such that $\nu L\sim \nu (K_X+ S+ B)$. 
In particular, $\nu B$ is an integral Weil divisor. 
We take the $\nu$-fold cyclic cover 
$\pi': Y'=\Spec_X\!\bigoplus _{i=0}^{\nu-1} \mathcal L^{-i}\to 
X$ associated to 
the section $\nu B\in |\mathcal L^{\nu}|$. 
More precisely, let $s\in H^0(X, \mathcal L^{\nu})$ be a section 
whose zero divisor is $\nu B$. 
Then the dual of $s:\mathcal O_X\to \mathcal L^{\nu}$ 
defines an $\mathcal O_X$-algebra structure on 
$\bigoplus ^{\nu-1}_{i=0} \mathcal L^{-i}$. 
Let $Y\to Y'$ be the normalization and 
let $\pi:Y\to X$ be the composition morphism. 
For the details, see \cite[3.5.~Cyclic covers]{ev}. 
We can take a finite cover 
$\varphi:V\to Y$ such that 
$V$ is smooth 
and that $T$ is a simple normal crossing 
divisor on $V$, where 
$\psi=\pi\circ \varphi$ and 
$T=\psi^*S$,  
by Kawamata's covering trick (cf.~\cite[3.17.~Lemma]{ev}). 
Let $\iota':Y\setminus \pi^*S\to Y$ be the natural open immersion 
and let $U$ be the smooth locus of $Y$. 
We denote the natural open immersion $U\to Y$ by $j$. 
We put $\widetilde {\Omega}^p_Y(\log (\pi^*S))=j_*\Omega^p_U(\log 
(\pi^*S))$ for every $p$. 
Then it can be checked easily that 
$$
\iota'_{!}\mathbb C_{Y\setminus \pi^*S}\overset{qis}\longrightarrow 
\widetilde {\Omega}^{\bullet}_Y(\log (\pi^*S))\otimes 
\mathcal O_Y(-\pi^*S) 
$$
is a direct summand of 
$$
\varphi_*(\iota_{!}\mathbb C_{V\setminus T})\overset{qis}\longrightarrow 
\varphi_*(\Omega^{\bullet}_V(\log T)\otimes \mathcal O_V(-T)),  
$$ 
where $qis$ means a quasi-isomorphism. 
On the other hand, 
we can decompose 
$\pi_*(\widetilde{\Omega}^{\bullet}_Y(\log (\pi^*S))\otimes 
\mathcal O_Y(-\pi^*S))$ and $\pi_*
(\iota'_{!}\mathbb C_{Y\setminus \pi^*S})$ 
into eigen components of the Galois action 
of $\pi:Y\to X$. 
We write these decompositions as follows, 
\begin{align*}
\pi_*(\iota'_!\mathbb C_{Y\setminus \pi^*S})=\bigoplus _{i=0}^{\nu-1}
\mathcal C_i\subset \bigoplus _{i=0}^{\nu-1}\mathcal L^{-i}
(\llcorner iB\lrcorner-S)=\pi_*\mathcal O_Y(-\pi^*S), 
\end{align*}
where $\mathcal C_i\subset \mathcal L^{-i}(\llcorner iB\lrcorner 
-S)$ for every $i$. We put $\mathcal C=\mathcal C_1$. 
We have that  
$$
\mathcal C\overset{qis}\longrightarrow 
\Omega^{\bullet}_X(\log (S+B))\otimes 
\mathcal L^{-1}(-S) 
$$ 
is a direct summand of 
$$
\psi_*(\iota_{!}\mathbb C_{V\setminus T})\overset{qis}
\longrightarrow \psi_*(\Omega^{\bullet}_V(\log T)\otimes 
\mathcal O_V(-T)). 
$$ 
The $E_1$-degeneration of the spectral sequence 
\begin{eqnarray*}
E^{pq}_1&=&H^q(V, \Omega^p_V(\log T)\otimes \mathcal O_V(-T))
\\ &\Rightarrow& \mathbb H^{p+q}(V, 
\Omega^{\bullet}_V(\log T)\otimes 
\mathcal O_V(-T))\simeq H^{p+q}(V, 
\iota_!\mathbb C_{V\setminus T})
\end{eqnarray*}
(cf.~\ref{say52}) implies the $E_1$-degeneration of 
\begin{eqnarray*}
E^{pq}_1&=&H^q(X, \Omega^p_X(\log (S+B))
\otimes \mathcal L^{-1}(-S))\\
& \Rightarrow &\mathbb H^{p+q}(X, 
\Omega^{\bullet}_X(\log (S+B))\otimes 
\mathcal L^{-1}(-S))
\simeq H^{p+q}(X, \mathcal C)
\end{eqnarray*}
Therefore, the inclusion $\mathcal C\subset \mathcal L^{-1}(-S)$ 
induces surjections 
$$
H^p(X, \mathcal C)\to H^p(X, \mathcal L^{-1}(-S))  
$$ 
for all $p$. 
We can check 
the following simple property by seeing the monodromy 
action of the Galois group of 
$\pi:Y\to X$ 
on $\mathcal C$ around $\Supp B$. 

\begin{cor}[{cf.~\cite[Corollary 2.54]{km}}]\label{cor43}
Let $U\subset X$ be 
a connected open set such that 
$U\cap \Supp B\ne \emptyset$. 
Then $H^0(U, \mathcal C|_{U})=0$. 
\end{cor}
This property is utilized via the following fact. 
The proof is obvious. 

\begin{lem}[{cf.~\cite[Lemma 2.55]{km}}]\label{lem44}
Let $F$ be a sheaf of Abelian groups on a topological 
space 
$X$ and let $F_1, F_2\subset F$ be subsheaves. 
Let $Z\subset X$ be a closed subset. Assume 
that 
\begin{itemize}
\item[(1)] $F_2|_{X\setminus Z}=F|_{X\setminus Z}$, and 
\item[(2)] if $U$ is connected, open 
and $U\cap Z\ne \emptyset$, then 
$H^0(U, F_1|U)=0$. 
\end{itemize}
Then $F_1$ is a subsheaf of $F_2$. 
\end{lem}

As a corollary, we obtain: 

\begin{cor}[{cf.~\cite[Corollary 2.56]{km}}]
Let $M\subset \mathcal L^{-1}(-S)$ be a subsheaf such that 
$M|_{X\setminus \Supp B}
=\mathcal L^{-1}(-S)|_{X\setminus \Supp B}$. 
Then the injection 
$$
\mathcal C\to \mathcal L^{-1}(-S) 
$$ 
factors as 
$$
\mathcal C \to M\to \mathcal L^{-1}(-S). 
$$ Therefore, 
$$H^i(X, M)\to H^i(X, \mathcal L^{-1}(-S))
$$ 
is surjective for every $i$. 
\end{cor}

\begin{proof}
The first part is clear from Corollary \ref{cor43} 
and Lemma \ref{lem44}. 
This implies that we have maps 
$$
H^i(X, \mathcal C)\to H^i(X, M)\to H^i(X, \mathcal L^{-1}(-S)). 
$$ 
As we saw above, the composition is surjective. 
Hence so is the map on the right. 
\end{proof}

Therefore, we obtain that 
$$H^q(X, \mathcal L^{-1}(-S-D))\to H^q(X, 
\mathcal L^{-1}(-S))$$ is 
surjective for every $q$. By the Serre duality, we obtain 
\begin{align*}
H^q(X, \mathcal O_X(K_X)\otimes \mathcal L(S))\to H^q(X, \mathcal 
O_X(K_X)\otimes \mathcal L(S+D))
\end{align*} is 
injective for every $q$. 
This means that $$H^q(X, \mathcal O_X(L))\to 
H^q(X, \mathcal O_X(L+D))$$ is injective for every $q$.
\end{proof}

\section{Injectivity, torsion-free, and vanishing theorems}\label{sec5}

In this section, we prove generalizations of 
Koll\'ar's torsion-freeness and vanishing theorem (cf.~Theorem \ref{thm53}). 
First, we prove a generalization of 
Koll\'ar's injectivity 
theorem (cf.~\cite[Theorem 3.1]{ambro}). 
It is a straightforward consequence of Proposition \ref{prop41} 
and will produce the desired torsion-free and 
vanishing theorems. 

\begin{thm}[Injectivity theorem]\label{thm51} Let $X$ 
be a smooth projective variety and let 
$\Delta$ be a boundary 
$\mathbb R$-divisor such that 
$\Supp \Delta$ is simple normal crossing. 
Let $L$ be a Cartier 
divisor on $X$ and let $D$ be an effective 
Cartier divisor that contains no lc centers of 
$(X, \Delta)$. 
Assume the following conditions. 
\begin{itemize}
\item[(i)] $L\sim _{\mathbb R}K_X+\Delta+H$, 
\item[(ii)] $H$ is a semi-ample $\mathbb R$-Cartier 
$\mathbb R$-divisor, and 
\item[(iii)] $tH\sim _{\mathbb R} D+D'$ for some 
positive real number $t$, where 
$D'$ is an effective $\mathbb R$-Cartier 
$\mathbb R$-divisor 
whose support contains no lc centers of $(X, \Delta)$. 
\end{itemize}
Then the homomorphisms 
$$
H^q(X, \mathcal O_X(L))\to H^q(X, \mathcal O_X(L+D))
$$ 
which are induced by the natural inclusion 
$\mathcal O_X\to \mathcal O_X(D)$ are injective for all $q$. 
\end{thm}

\begin{proof}
We put $S=\llcorner \Delta\lrcorner$ and $B=\{\Delta\}$. 
We can take a resolution $f:Y\to X$ 
such that $f$ is an isomorphism 
outside $\Supp (D+D'+B)$, and that the 
union of the support of $f^*(S+B+D+D')$ and the 
exceptional locus of $f$ has a 
simple normal crossing support on $Y$. 
Let $B'$ be the strict transform of $B$ on $Y$. 
We write $$K_Y+S'+B'=f^*(K_X+S+B)+E, $$ 
where $S'$ is the strict transform of $S$ and 
$E$ is $f$-exceptional. 
It is easy to see that $E_+=\ulcorner E\urcorner \geq 0$. 
We put $L'=f^*L+E_+$ and $E_-=E_+-E\geq 0$. 
We note that $E_+$ is Cartier and $E_-$ is 
an effective $\mathbb R$-Cartier $\mathbb R$-divisor 
with $\llcorner E_{-}\lrcorner=0$. 
Since $f^*H$ is semi-ample, 
we can write $f^*H\sim _{\mathbb R}\sum _i a_i H'_i$, where 
$0<a_i<1$ and $H'_i$ is a general 
Cartier divisor on $Y$ for every $i$. 
We put $$B''=B'+E_-+\frac{\varepsilon}{t} 
f^*(D+D')+(1-\varepsilon) \sum_i a_i  H'_i$$ for 
some $0<\varepsilon \ll 1$. 
Then $L'\sim _{\mathbb R}K_Y+S'+B''$. 
By the construction, $\llcorner B''\lrcorner =0$, the 
support of $S'+B''$ is simple normal crossing 
on $Y$, and $\Supp B''\supset \Supp f^*D$. 
So, Proposition \ref{prop41} implies that 
the homomorphisms 
$$H^q(Y, \mathcal O_Y(L'))\to H^q(Y, \mathcal O_Y(L'+f^*D))$$ are 
injective for all $q$. 
It is easy to see that $f_*\mathcal O_Y(L')\simeq 
\mathcal O_X(L)$. 
By Lemma \ref{lem21}, 
we can write $L'\sim _{\mathbb Q}K_Y+S'+B'''$, where 
$B'''$ is a $\mathbb Q$-divisor on $Y$ such that 
$\llcorner B'''\lrcorner =0$ and $\Supp B'''=\Supp B''$. 
Thus, by Lemma \ref{lem52} below, $R^qf_*\mathcal O_Y(L')=0$ for 
all $q>0$. 
By the 
Leray spectral sequence, 
the homomorphisms $$H^q(X, \mathcal O_X(L))\to 
H^q(X, \mathcal O_X(L+D))$$ are injective for all $q$. 
\end{proof}

Let us recall the following well-known easy lemma. 

\begin{lem}[Reid--Fukuda type]\label{lem52}
Let $V$ be a smooth projective variety and 
let $B$ be a boundary $\mathbb Q$-divisor 
on $V$ such that 
$\Supp B$ is simple normal crossing. 
Let $f:V\to W$ be a projective birational morphism 
onto a variety $W$. 
Assume that 
$f$ is an isomorphism 
at the generic point of 
every lc center of $(V, B)$ and 
that $D$ is a Cartier divisor on $V$ such that 
$D-(K_V+B)$ is nef. Then 
$R^if_*\mathcal O_V(D)=0$ for every $i>0$. 
\end{lem}
\begin{proof}
We use the induction on the number of 
irreducible components of $\llcorner B\lrcorner$ and 
on the dimension of $V$. 
If $\llcorner B\lrcorner =0$, then 
the lemma follows from the Kawamata--Viehweg 
vanishing theorem (cf.~\cite[Corollary 2.68]{km}). Therefore, 
we can assume that there is an irreducible 
divisor $S\subset \llcorner B\lrcorner$. 
We 
consider the following short exact sequence 
$$
0\to \mathcal O_V(D-S)\to \mathcal O_V(D)\to 
\mathcal O_S(D)\to 0. 
$$ 
By induction, we see that 
$R^if_*\mathcal O_V(D-S)=0$ 
and $R^if_*\mathcal O_S(D)=0$ for 
every $i>0$. Thus, we have 
$R^if_*\mathcal O_V(D)=0$ for $i>0$. 
\end{proof}

The next theorem is the main theorem of this section 
(cf.~\cite{ambro}). 
See also \cite{fujino-higher}. 

\begin{thm}[Torsion-freeness and 
vanishing theorem]\label{thm53}Let 
$Y$ be a smooth variety and 
let $B$ be a boundary $\mathbb R$-divisor such that 
$\Supp B$ is simple normal crossing. 
Let $f:Y\to X$ be a projective morphism and let $L$ be a Cartier 
divisor on $Y$ such that 
$L-(K_Y+B)$ is $f$-semi-ample. 
\begin{itemize}
\item[(i)] Let $q$ be an 
arbitrary non-negative integer. 
Every non-zero local section of $R^qf_*\mathcal O_Y(L)$ contains 
in its support the $f$-image of 
some stratum of $(Y, B)$. 
\item[(ii)] 
Let $\pi:X\to S$ be a projective morphism. 
Assume that $L-(K_X+B)\sim _{\mathbb R}f^*H$ for 
some $\pi$-ample $\mathbb R$-Cartier 
$\mathbb R$-divisor $H$ on $X$. 
Then $R^p\pi_*R^qf_*\mathcal O_Y(L)=0$ for every $p>0$ 
and $q\geq 0$. 
\end{itemize}
\end{thm}

\begin{rem}
It is obvious that the statement of Theorem \ref{thm53} (i) is equivalent to 
the following one. 
\begin{itemize}
\item[(i$^\prime$)] Let $q$ be an arbitrary non-negative integer. 
Every associated prime of $R^qf_*\mathcal O_Y(L)$ is the generic point of 
the $f$-image of some stratum of $(Y, B)$. 
\end{itemize}
\end{rem}

Let us start the 
proof of Theorem \ref{thm53}. 

\begin{proof}[Proof of {\em{Theorem \ref{thm53}}}]
We take an $f$-semi-ample $\mathbb R$-Cartier 
$\mathbb R$-divisor $M$ on $Y$ such that 
$M\sim _{\mathbb R}L-(K_Y+B)$. 

(i) We divide the proof into two steps. 
\setcounter{step}{0}
\begin{step}\label{p1}
First, we assume that $X$ is projective. 
We can assume that $M$ is semi-ample by replacing 
$L$ (resp.~$M$) with 
$L+f^*A'$ (resp.~$M+f^*A'$), where 
$A'$ is a very ample Cartier divisor on $X$. 
Assume that $R^qf_*\mathcal O_Y(L)$ 
has a local section whose support does not contain 
the images of any 
$(Y, B)$-strata. More precisely, 
let $U$ be a non-empty Zariski open set and let $s\in 
\Gamma (U, R^qf_*\mathcal O_Y(L))$ be a non-zero 
section of $R^qf_*\mathcal O_Y(L)$ 
on $U$ whose support $V\subset U$ does not contain 
the $f$-images of any strata of $(Y, B)$. 
Let $\overline V$ be the closure of $V$ in $X$. We note that 
$\overline V\setminus V$ may contain the $f$-image of some stratum of $(Y, B)$. 
By replacing $Y$ with its blow-up along an lc center which is mapped into 
$\overline V\setminus V$, we can assume that 
an irreducible component $B_0$ of $\llcorner B\lrcorner$ is mapped into 
$\overline V\setminus V$ by $f$. 
We note that $M\sim _{\mathbb R}L-B_0-(K_X+B-B_0)$. 
We replace $L$ (resp.~$B$) with 
$L-B_0$ (resp.~$B-B_0$). By repeating this process finitely 
many times, we can assume that $\overline V$ does 
not contain the $f$-images of any strata of $(Y, B)$. 
Then 
we can find a very ample Cartier divisor 
$A$ with the following properties. 
\begin{itemize}
\item[(a)] $f^*A$ contains no lc centers of $(Y, B)$, and 
\item[(b)] $R^qf_*\mathcal O_Y(L)\to R^qf_*\mathcal O_Y(L)\otimes 
\mathcal O_X(A)$ is not injective. 
\end{itemize}
We can assume that $M-f^*A$ is semi-ample by replacing 
$L$ (resp.~$M$) with 
$L+f^*A$ (resp.~$M+f^*A$). 
If necessary, we replace $L$ (resp.~$M$) with 
$L+f^*A''$ (resp.~$M+f^*A''$), where $A''$ is a very ample 
Cartier divisor on $X$. 
Then, we have 
$$H^0(X, R^qf_*\mathcal O_Y(L))\simeq 
H^q(Y, \mathcal O_Y(L))$$ and 
$$H^0(X, R^qf_*\mathcal O_Y(L)\otimes \mathcal O_X(A))\simeq 
H^q(Y, \mathcal O_Y(L+f^*A)). 
$$
We see that 
$$H^0(X, R^qf_*\mathcal O_Y(L))
\to H^0(X, R^qf_*\mathcal O_Y(L)\otimes \mathcal O_X(A))
$$ 
is not injective by (b) if $A''$ is sufficiently ample. 
So, 
$$H^q(Y, \mathcal O_Y(L))
\to 
H^q(Y, \mathcal O_Y(L+f^*A))$$ is not injective. 
It contradicts Theorem \ref{thm51}. 
We finish the proof when $X$ is projective. 
\end{step}
\begin{step}\label{p2}
Next, we assume that $X$ is not projective. 
Note that the problem is local. So, we can shrink $X$ and assume that 
$X$ is affine. 
By the argument similar to the one in Step \ref{p1} in the proof of 
(ii) below, we can assume that $M$ is a semi-ample $\mathbb Q$-Cartier 
$\mathbb Q$-divisor. 
We compactify $X$ and apply Lemma \ref{lem25}. 
Then we obtain a compactification $\overline f: \overline Y\to \overline X$ of 
$f:Y\to X$. 
Let $\overline M$ be the closure of $M$ on 
$\overline Y$. If $\overline M$ is not a semi-ample $\mathbb Q$-Cartier 
$\mathbb Q$-divisor, then we take blowing-ups 
of $\overline Y$ inside $\overline Y\setminus Y$ 
and obtain a semi-ample $\mathbb Q$-Cartier $\mathbb Q$-divisor 
$\overline M$ on 
$\overline Y$ such that $\overline M|_{Y}=M$. 
Let $\overline L$ (resp.~$\overline B$) be the 
closure of $L$ (resp.~$B$) on $\overline Y$. 
We note that $\overline M\sim _{\mathbb R}\overline L-(K_{\overline Y}
+\overline B)$ does not necessarily hold. 
We can write $M+\sum _i a_i (f_i)=L-(K_Y+B)$, 
where $a_i$ is a real number and $f_i$ is a rational 
function on $Y$ for 
every $i$. 
We put $$E=\overline M+\sum _i a_i (f_i)-
(\overline L-(K_{\overline Y}+
\overline B)). $$  
We replace $\overline L$ (resp.~$\overline B)$ with 
$\overline L+\ulcorner E\urcorner$ (resp.~$\overline B+\{-E\}$). 
Then we obtain the desired property of 
$R^q\overline f_*\mathcal O_{\overline Y}
(\overline L)$ since $\overline X$ is projective. We 
note that $\Supp E$ is in $\overline Y\setminus Y$. 
So, this completes the whole proof. 
\end{step}

(ii) We divide the proof into three steps. 
\setcounter{step}{0}
\begin{step}\label{o1} 
We assume that $\dim S=0$. 
The following arguments are well known and 
standard. 
We describe them for the reader's convenience. 
In this case, we can write $H\sim _{\mathbb R}H_1+H_2$, 
where $H_1$ (resp.~$H_2$) is a $\pi$-ample $\mathbb Q$-Cartier 
$\mathbb Q$-divisor 
(resp.~$\pi$-ample $\mathbb R$-Cartier $\mathbb R$-divisor) on $X$. 
So, we can write $H_2\sim _{\mathbb R}\sum _i a_i H'_i$, 
where $0<a_i<1$ and $H'_i$ is a general very ample Cartier 
divisor on $X$ for every $i$. 
Replacing $B$ (resp.~$H$) with $B+\sum _i a_i f^*H'_i$ 
(resp.~$H_1$), we can assume that 
$H$ is a $\pi$-ample $\mathbb Q$-Cartier $\mathbb Q$-divisor. 
We take a general member $A\in |mH|$, 
where $m$ is a sufficiently divisible positive integer, such 
that $A'=f^*A$ and 
$R^qf_*\mathcal O_Y(L+A')$ is $\pi_*$-acyclic for all $q$. 
By (i), we have the following short exact sequences, 
$$
0\to R^qf_*\mathcal O_Y(L)\to 
R^qf_*\mathcal O_Y(L+A')
\to R^qf_*\mathcal O_{A'}(L+A')\to 0. 
$$ 
for all $q$. 
Note that $R^qf_*\mathcal O_{A'}(L+A')$ is $\pi_*$-acyclic by 
induction on $\dim X$ and $R^qf_*\mathcal O_Y(L+A')$ is also 
$\pi_*$-acyclic by the above assumption. 
Thus, $E^{pq}_2=0$ for $p\geq 2$ in the 
following commutative diagram of spectral sequences. 
$$
\xymatrix{
&E^{pq}_2=R^p\pi_*R^qf_*\mathcal O_Y(L) \ar[d]_{\varphi^{pq}}
\ar@{=>}[r]&R^{p+q}(\pi\circ f)_*
\mathcal O_Y(L)\ar[d]_{\varphi^{p+q}}\\ 
&\overline{E}^{pq}_2=R^p\pi_*R^qf_*\mathcal O_Y(L+A') 
\ar@{=>}[r]& R^{p+q}
(\pi\circ f)_*
\mathcal O_Y(L+A')
}
$$
We note that $\varphi^{1+q}$ is injective by Theorem \ref{thm51}. 
We have $E^{1q}_2\to R^{1+q}(\pi\circ f)_*\mathcal O_Y(L)$ 
is injective by the fact that 
$E^{pq}_2=0$ for $p\geq 2$. 
We also have that $\overline 
{E}^{1q}_2=0$ by the above assumption. 
Therefore, we 
obtain $E^{1q}_2=0$ since the 
injection $E_2^{1q}\to R^{1+q}(\pi\circ f)_*\mathcal O_Y(L+A')$ factors 
through $\overline {E}_2^{1q}=0$. 
This implies that $R^p\pi_*R^qf_*\mathcal O_Y(L)=0$ for every 
$p>0$ and $q\geq 0$. 
\end{step}

\begin{step}\label{o2}  
We assume that $S$ is projective. 
By replacing $H$ (resp.~$L$) with 
$H+\pi^*G$ (resp.~$L+(\pi\circ f)^*G$), where 
$G$ is a very ample Cartier divisor on $S$, 
we can assume that $H$ is an ample $\mathbb R$-Cartier 
$\mathbb R$-divisor. 
By the same argument as in Step 1, we can assume that $H$ is an ample 
$\mathbb Q$-Cartier $\mathbb Q$-divisor and $M\sim _{\mathbb Q}f^*H$. 
If $G$ is a 
sufficiently ample Cartier 
divisor on $S$, $H^k(S, R^p\pi_*R^qf_*\mathcal O_Y(L)\otimes \mathcal O_S(G))
=0$ for 
every $k\geq 1$, 
\begin{align*}
&H^0(S, R^p\pi_*R^qf_*\mathcal O_Y(L)\otimes \mathcal O_S(G))
\\ &\simeq H^p(X, R^qf_*\mathcal O_Y(L)\otimes \mathcal O_X(\pi^*G))
\\ &\simeq H^p(X, R^qf_*\mathcal O_Y(L+f^*\pi^*G)), 
\end{align*} 
and $R^p\pi_*R^qf_*\mathcal O_Y(L)\otimes \mathcal O_S(G)$ is generated by its 
global sections. 
Since $$M+f^*\pi^*G\sim _{\mathbb R}
L+f^*\pi^*G-(K_Y+B), $$ 
$$M+f^*\pi^*G\sim _{\mathbb Q}f^*(H+\pi^*G), $$
and $H+\pi^*G$ is ample, 
we can apply Step 1 and obtain $H^p(X, R^qf_*\mathcal O_Y(L+
f^*\pi^*G))=0$ for every $p>0$ and $q\geq 0$. 
Thus, $R^p\pi_*R^qf_*\mathcal 
O_Y(L)=0$ for every $p>0$ and $q\geq 0$ by the above arguments.
\end{step}
\begin{step}\label{o3}  
When $S$ is not projective, we shrink $S$ and assume 
that $S$ is affine. 
By the same argument as in Step \ref{o1} above, we can assume that 
$H$ is $\mathbb Q$-Cartier. 
We compactify $S$ and $X$, and can assume that 
$S$ and $X$ are projective. 
By Lemma \ref{lem25}, we can reduce it to the case when 
$S$ is projective. This step is essentially the same as Step \ref{p2} 
in the proof of (i). So, we omit the details here.  
\end{step}
We obtained the statement (ii). 
\end{proof}

\section{Non-lc ideal sheaves}\label{sec6}

We introduce the notion of {\em{non-lc ideal sheaves}}. 
It is an analogue of the usual {\em{multiplier ideal 
sheaves}} (see, for example, \cite[Chapter 9]{laza}). 
For details, see \cite{non-lc} and 
\cite{ft}. 

\begin{defn}[Non-lc ideal sheaf]\label{def61}
Let $X$ be a normal variety and let $B$ be 
an $\mathbb R$-divisor on $X$ such that 
$K_X+B$ is $\mathbb R$-Cartier. 
Let $f:Y\to X$ be a resolution with $K_Y+B_Y=f^*(K_X+B)$ 
such that 
$\Supp B_Y$ is simple normal crossing. 
Then we put 
\begin{align*}
\mathcal J_{NLC}(X, B)
&=f_*\mathcal O_Y(\ulcorner -(B_Y^{<1})\urcorner-
\llcorner B_Y^{>1}\lrcorner)\\& =f_*\mathcal O_Y
(-\llcorner B_Y\lrcorner 
+B^{=1}_Y)
\end{align*} and 
call it the {\em{non-lc ideal sheaf associated to $(X, B)$}}. 
If $B$ is effective, then 
$\mathcal J_{NLC}(X, B)\subset \mathcal O_X$. 
\end{defn}

The ideal sheaf $\mathcal J_{NLC}(X, B)$ is well-defined 
by the following 
easy lemma. 

\begin{lem}\label{lem62}
Let $g:Z\to Y$ be a proper birational morphism between 
smooth varieties and let $B_Y$ be an 
$\mathbb  R$-divisor on $Y$ such 
that $\Supp B_Y$ is simple normal crossing. 
Assume that $K_Z+B_Z=g^*(K_Y+B_Y)$ and 
that $\Supp B_Z$ is simple normal crossing. 
Then we have $$g_*\mathcal O_Z(\ulcorner -(B^{<1}_Z)\urcorner 
-\llcorner B^{>1}_Z\lrcorner)\simeq 
\mathcal O_Y(\ulcorner -(B^{<1}_Y)\urcorner 
-\llcorner B^{>1}_Y\lrcorner). $$ 
\end{lem}
\begin{proof}
By $K_Z+B_Z=g^*(K_Y+B_Y)$, we 
obtain 
\begin{align*}
K_Z=&g^*(K_Y+B^{=1}_Y+\{B_Y\})\\&+g^*
(\llcorner B^{<1}_Y\lrcorner+\llcorner B^{>1}_Y\lrcorner)
-(\llcorner B^{<1}_Z\lrcorner+\llcorner B^{>1}_Z\lrcorner)
-B^{=1}_Z-\{B_Z\}.
\end{align*} 
If $a(\nu, Y, B^{=1}_Y+\{B_Y\})=-1$ for a prime divisor 
$\nu$ over $Y$, then 
we can check that $a(\nu, Y, B_Y)=-1$ by using 
\cite[Lemma 2.45]{km}. 
Since $g^*
(\llcorner B^{<1}_Y\lrcorner+\llcorner B^{>1}_Y\lrcorner)
-(\llcorner B^{<1}_Z\lrcorner+\llcorner B^{>1}_Z\lrcorner)$ is 
Cartier, we can easily see that 
$$g^*(\llcorner B^{<1}_Y\lrcorner+\llcorner B^{>1}_Y\lrcorner)
=\llcorner B^{<1}_Z\lrcorner+\llcorner B^{>1}_Z\lrcorner+E, $$ 
where $E$ is an effective $f$-exceptional Cartier divisor. 
Thus, we obtain 
$$g_*\mathcal O_Z(\ulcorner -(B^{<1}_Z)\urcorner 
-\llcorner B^{>1}_Z\lrcorner)\simeq 
\mathcal O_Y(\ulcorner -(B^{<1}_Y)\urcorner 
-\llcorner B^{>1}_Y\lrcorner).$$  
This completes the proof. 
\end{proof}

\begin{rem}\label{rem63} 
We use the same notation as in Definition \ref{def61}. 
We put 
$$
\mathcal J(X, B)=f_*\mathcal O_Y(-\llcorner B_Y\lrcorner). 
$$ 
This sheaf $\mathcal J(X, B)$ is well known 
as the (algebraic version of) {\em{multiplier ideal 
sheaf}} of the pair $(X, B)$. 
See, for example, \cite[Chapter 9]{laza}. 
\end{rem}

By the definition, the following proposition is obvious. 

\begin{prop}
Let $X$ be a normal variety and let $B$ be an effective 
$\mathbb R$-divisor on $X$ such that 
$K_X+B$ is $\mathbb R$-Cartier. 
Then $(X, B)$ is log canonical if and only if 
$\mathcal J_{NLC}(X, B)=\mathcal O_X$. 
\end{prop}

The next proposition is a kind of Bertini's theorem. 

\begin{prop}\label{prop63} 
Let $X$ be a smooth 
variety and let $B$ be an effective $\mathbb R$-divisor 
on $X$ such that 
$K_X+B$ is $\mathbb R$-Cartier. 
Let $\Lambda$ be a linear system on $X$ and 
let $D\in \Lambda$ be a general member of $\Lambda$. 
Then 
$$
\mathcal J_{NLC}(X, B+tD)=\mathcal J_{NLC}(X, B)
$$ 
outside the 
base locus $\Bs \Lambda$ of 
$\Lambda$ for all $0\leq t\leq 1$. 
\end{prop}
\begin{proof}
By replacing $X$ with $X\setminus \Bs \Lambda$, we 
can assume that $\Bs \Lambda=\emptyset$. 
Let $f:Y\to X$ be a resolution as in Definition \ref{def61}. 
Since $D$ is a general member of $\Lambda$, 
$f^*D=f^{-1}_*D$ is a smooth divisor on $Y$ 
such that $\Supp f^*D\cup \Supp B_Y$ is 
simple normal crossing. Therefore, 
we can check that $\mathcal J_{NLC}(X, B+tD)=\mathcal J_{NLC}(X, B)$ for 
all $0\leq t\leq 1$. 
\end{proof}

We close this section with an important remark. 

\begin{rem}\label{rem64} 
In the subsequent sections (Sections \ref{sec7}, \ref{sec11}, 
\ref{sec12}, and \ref{sec14}), 
we consider the scheme structure of $\Nlc(X, B)$ 
defined by $\mathcal J_{NLC}(X, B)$. 
However, we can use $\mathcal J'(X, B)$ or 
$\mathcal J'_l(X, B)$ for any negative integer $l$ 
in place of $\mathcal J_{NLC}(X, B)$. 
For the definitions and basic properties of 
$\mathcal J'(X, B)$ and $\mathcal J'_l(X, B)$, see 
\cite{ft}. 
We adopt $\mathcal J_{NLC}(X, B)$ since 
we think $\mathcal J_{NLC}(X, B)$ is the most natural defining 
ideal sheaf of $\Nlc (X, B)$.  
\end{rem}

\section{Vanishing theorem}\label{sec7} 
The following vanishing theorem, which is a special case of 
\cite[Theorem 4.4]{ambro}, 
is one of the key results in this paper. 
We note that the importance of Theorem \ref{thm71} is 
in its formulation best suited for 
new frameworks explained in subsequent sections. 
For the details of Ambro's original 
statement, see \cite[Theorem 4.4]{ambro} or \cite[Theorem 3.39]{book}.
 
\begin{thm}\label{thm71} 
Let $X$ be a normal 
variety and $B$ an effective $\mathbb R$-divisor 
on $X$ such that 
$K_X+B$ is $\mathbb R$-Cartier. 
Let $D$ be a Cartier divisor on $X$. 
Assume that $D-(K_X+B)$ is $\pi$-ample, 
where $\pi:X\to S$ is a projective morphism onto a variety $S$. 
Let $\{C_i\}$ be {\em{any}} set of 
lc centers of the pair $(X, B)$. 
We put $W=\bigcup C_i$ with the reduced scheme structure. 
Assume that $W$ is disjoint from $\Nlc (X, B)$. 
Then we have 
$$
R^i\pi_*(\mathcal J\otimes \mathcal O_X(D))=0  
$$
for every $i>0$, where $\mathcal J=\mathcal I_W\cdot \mathcal 
J_{NLC}(X, B)\subset \mathcal O_X$ and 
$\mathcal I_W$ is the defining 
ideal sheaf of $W$ on $X$. 
Therefore, the restriction 
map 
$$
\pi_*\mathcal O_X(D)\to \pi_*\mathcal O_W(D)\oplus 
\pi_*\mathcal O_{\Nlc(X, B)}(D)
$$ 
is surjective and 
$$
R^i\pi_*\mathcal O_W(D)=0 
$$ 
for every $i>0$. 
In particular, the restriction maps 
$$\pi_*\mathcal O_X(D)\to \pi_*\mathcal O_W(D)$$
and $$
\pi_*\mathcal O_X(D)\to \pi_*\mathcal O_{\Nlc(X, B)}(D)
$$ 
are surjective. 
\end{thm}

\begin{proof}
Let $f:Y\to X$ be a resolution 
such that $\Supp f^{-1}_*B\cup \Exc (f)$ 
is a simple 
normal crossing 
divisor. 
We can further assume that 
$f^{-1}(W)$ is a simple normal crossing 
divisor on $Y$. We can write 
$$
K_Y+B_Y=f^*(K_X+B). 
$$ 
Let $T$ be the union of the irreducible 
components of $B^{=1}_Y$ that 
are mapped into $W$ by $f$. 
We consider 
the following 
short exact sequence 
$$
0\to \mathcal O_Y(A-N-T)\to \mathcal O_Y(A-N)\to \mathcal 
O_T(A-N)\to 0, 
$$ 
where $A=\ulcorner -(B^{<1}_Y)\urcorner$ 
and $N=\llcorner B^{>1}_Y\lrcorner$. 
Note that $A$ is an effective 
$f$-exceptional 
divisor. 
We obtain the following long exact sequence 
\begin{align*}
0&\to f_*\mathcal O_Y(A-N-T)\to f_*\mathcal O_Y(A-N)
\to f_*\mathcal O_T(A-N)\\ 
&\overset{\delta}\to R^1f_*\mathcal O_Y(A-N-T)\to \cdots.  
\end{align*}
Since 
\begin{align*}
&A-N-T-(K_Y+\{B_Y\}+B^{=1}_Y-T)=-(K_Y+B_Y)\\ &\sim_{\mathbb R}
-f^*(K_X+B), 
\end{align*}
every non-zero local section of $R^1f_*\mathcal O_Y(A-N-T)$ contains 
in its support the 
$f$-image of some stratum of $(Y, \{B_Y\}+B^{=1}_Y-T)$ by 
Theorem \ref{thm53} (i). 
On the other hand, $W=f(T)$. Therefore, 
the connecting homomorphism $\delta$ is a zero map. 
Thus, we have a short 
exact sequence 
\begin{align*}\label{siki-heart}
0\to f_*\mathcal O_Y(A-N-T)\to f_*\mathcal O_Y(A-N)\to 
f_*\mathcal O_T(A-N)\to 0. \tag{$\diamondsuit$} 
\end{align*}
We put $\mathcal J=f_*\mathcal O_Y(A-N-T)\subset \mathcal O_X$. 
Since $W$ is disjoint from $\Nlc(X, B)$, 
the ideal 
sheaf 
$\mathcal J$ coincides with $\mathcal I_W$ (resp.~$\mathcal J_{NLC}(X, B)$) 
in a neighborhood of $W$ (resp.~$\Nlc (X, B)$). Therefore, 
$\mathcal J=\mathcal I_W\cdot \mathcal J_{NLC}(X, B)$. 
We put $U=X\setminus \Nlc (X, B)$ and $V=f^{-1}(U)$. 
By restricting (\ref{siki-heart}) to $U$, 
we obtain 
$$
0\to f_*\mathcal O_V(A-T)\to f_*\mathcal O_V(A)\to f_*\mathcal O_T(A)\to 0. 
$$
Since $f_*\mathcal O_V(A)\simeq \mathcal O_U$, we 
have $f_*\mathcal O_T(A)\simeq \mathcal O_W$. 
The isomorphism $f_*\mathcal O_T(A)\simeq \mathcal O_W$ 
plays crucial roles in the next section. 
Thus we write it as a proposition. 
\begin{prop}\label{lem72} 
We have $f_*\mathcal O_T(A)\simeq \mathcal O_W$. 
It obviously implies that $f_*\mathcal O_T\simeq \mathcal O_W$ since 
$A$ is effective. 
\end{prop}
\begin{rem}
We did not use $D$ nor $\pi:X\to S$ to obtain Proposition \ref{lem72}. 
\end{rem}
Since 
$$
f^*D+A-N-T-(K_Y+\{B_Y\}+B^{=1}_Y-T)\sim_{\mathbb R}
f^*(D-(K_X+B)), 
$$ 
we have 
$$
R^i\pi_*(\mathcal J\otimes \mathcal O_X(D))\simeq 
R^i\pi_*(f_*\mathcal O_Y(A-N-T)\otimes \mathcal O_X(D))=0 
$$ 
for every $i>0$ by Theorem \ref{thm53} (ii). 
By considering the short exact sequence 
$$
0\to \mathcal J\to \mathcal J_{NLC}(X, B)\to \mathcal O_W\to 0, 
$$ 
we obtain 
\begin{align*}
\cdots &\to R^i\pi_*(\mathcal J_{NLC}(X, B)\otimes 
\mathcal O_X(D))\\ &\to R^i\pi_*\mathcal O_W(D)\to 
R^{i+1}\pi_*\mathcal (\mathcal J\otimes \mathcal O_X(D))\to \cdots.  
\end{align*} 
Since we have already 
checked 
$$
R^i\pi_*(\mathcal J_{NLC}(X, B)\otimes \mathcal O_X(D))
=R^i\pi_*(\mathcal J\otimes \mathcal O_X(D))=0 
$$ 
for every $i>0$, we have 
$R^i\pi_*\mathcal O_W(D)=0$ for 
all $i>0$. 
Finally, we consider the following 
short exact sequence 
$$
0\to \mathcal J\to \mathcal O_X\to \mathcal O_W\oplus \mathcal O_{\Nlc(X, B)}
\to 0. 
$$
By taking $\otimes \mathcal O_X(D)$ and $R^i\pi_*$, 
we obtain 
$$
0\to \pi_*(\mathcal J\otimes \mathcal O_X(D))\to
\pi_*\mathcal O_X(D)\to \pi_*\mathcal O_W(D)
\oplus \pi_*\mathcal O_{\Nlc(X, B)}(D)\to 0.  
$$ 
This completes the proof. 
\end{proof}

\section{Lc centers}\label{sec8}

We prove the basic properties of 
lc centers as an application of the result 
in the preceding section (cf.~Proposition \ref{lem72}). 
Theorem \ref{thm81} is very useful 
in the study of linear systems on log canonical pairs. 
It can not be proved by the traditional method based on 
the Kawamata--Viehweg--Nadel vanishing theorem in the sense
that the coefficients of $B$ cannot be perturbed in general.

\begin{thm}[{cf.~\cite[Propositions 4.7 and 4.8]{ambro}}]\label{thm81}
Let $X$ be a normal variety and let $B$ be an effective 
$\mathbb R$-divisor 
such that 
$(X, B)$ is log canonical. 
Then we have the following properties. 
\begin{itemize}
\item[(1)] $(X, B)$ has at most finitely many lc centers. 
\item[(2)] An intersection of two lc centers is 
a union of lc centers. 
\item[(3)] Any union of lc centers of $(X, B)$ is semi-normal. 
\item[(4)] Let $x\in X$ be a closed point such that 
$(X, B)$ is lc but not klt at $x$. Then 
there is a unique minimal lc center $W_x$ passing through 
$x$, and $W_x$ is normal at $x$. 
\end{itemize}
\end{thm}
\begin{proof} 
We use the notation in the proof of Theorem \ref{thm71}. 
(1) is obvious. 
(3) is also obvious by Proposition \ref{lem72} 
since $T$ is a simple normal crossing 
divisor. 
Let $C_1$ and $C_2$ be two lc centers of $(X, B)$. 
We fix a closed point $P\in C_1\cap C_2$. 
For the 
proof of (2), it is enough to find an lc center $C$ 
such that $P\in C\subset C_1\cap C_2$. 
We put 
$W=C_1\cup C_2$. 
By Proposition \ref{lem72}, 
we obtain $f_*\mathcal O_T\simeq \mathcal O_W$. This means 
that $f:T\to W$ has connected fibers. 
We note that $T$ is a simple normal crossing divisor on $Y$. 
Thus, there exist irreducible 
components $T_1$ and $T_2$ of $T$ such that  
$T_1\cap T_2\cap f^{-1}(P)\ne \emptyset$ and that 
$f(T_i)\subset C_i$ for $i=1, 2$. 
Therefore, we can find an lc center $C$ with 
$P\in C\subset C_1\cap C_2$. 
We finish the proof of (2). 
Finally, we will prove (4). 
The existence and the uniqueness of 
the minimal lc center follow from (2). 
We take the unique minimal lc center $W=W_x$ passing through 
$x$. By Proposition \ref{lem72}, 
we have $f_*\mathcal O_T\simeq \mathcal O_W$. 
By shrinking $W$ around $x$, we can assume that 
every stratum of $T$ dominates $W$. 
Thus, $f:T\to W$ factors through 
the normalization $W^{\nu}$ of $W$. Since 
$f_*\mathcal O_T\simeq \mathcal O_W$, we obtain 
that $W^{\nu}\to W$ is an isomorphism. So, we 
obtain (4). 
\end{proof}

\section{Dlt blow-ups}\label{sec9} 
In this section, we discuss {\em{dlt blow-ups}} 
by Hacon (cf.~Theorem \ref{thm91}). 
In the subsequent sections, we will only use Lemma \ref{lem92}
(well known to experts)
and Theorem \ref{thm91}. 
For details, see Sections \ref{sec10} and \ref{sec16.5}. 
We also discuss a slight refinement of 
dlt blow-ups (cf.~Theorem \ref{gong}),
which is useful for future studies of log canonical pairs and 
has already played crucial 
roles in the study of log canonical weak Fano pairs (cf.~\cite{gongyo}). 

Let us recall the definition of {\em{dlt pairs}}. For 
another definition and the basic properties 
of dlt pairs, see \cite[Section 2.3]{km} and 
\cite{what}. 

\begin{defn}[Dlt pair]\label{def91} 
Let $X$ be a normal variety and 
let $B$ be an effective $\mathbb R$-divisor 
on $X$ such that $K_X+B$ is $\mathbb R$-Cartier. 
If there exists a resolution $f:Y\to X$ such that 
\begin{itemize}
\item[(i)] both $\Exc (f)$ and $\Exc (f)\cup \Supp f^{-1}_*B$ are 
simple normal crossing divisors on $Y$, and 
\item[(ii)] $a(E, X, B)>-1$ for every exceptional divisor 
$E\subset Y$, 
\end{itemize}
then $(X, B)$ is called {\em{divisorial log terminal}} 
({\em{dlt}}, for short). 
\end{defn}

We will use the following lemma in Section \ref{sec10}. 
For the details, see \cite[3.9 Adjunction for 
dlt pairs]{what}. 

\begin{lem}\label{lem92}
Let $(X, B)$ be a dlt pair and 
let $V$ be an lc center of $(X, B)$. 
Then $K_V+B_V=(K_X+B)|_V$ is dlt by 
adjunction.
\end{lem}

We borrow the next  theorem from \cite{bchm}. 

\begin{thm}[{cf.~\cite[Theorem 1.2]{bchm}}]\label{thm92}
Let $(X, B)$ be a klt pair, 
where $K_X+B$ is $\mathbb R$-Cartier. 
Let $\pi:X\to S$ be a projective 
birational morphism of quasi-projective 
varieties. 
Then $(X, B)$ has a log terminal model 
over $S$. This means that 
there exists a projective birational morphism 
$f:X'\to S$ such that 
\begin{itemize}
\item[(i)] $X'$ is $\mathbb Q$-factorial, 
\item[(ii)] $\phi^{-1}$ has no exceptional divisors, 
where $\phi=f^{-1}\circ \pi:X\dashrightarrow X'$, 
\item[(iii)] $K_{X'}+B'$ is $f$-nef, 
where $B'=\phi_*B$, and  
\item[(iv)] $a(E, X, B)<a(E, X', B')$ for 
every $\phi$-exceptional divisor $E\subset X$. 
\end{itemize}
\end{thm}

The following theorem is very useful. It is 
a consequence of Theorem \ref{thm92}. 

\begin{thm}[Hacon]\label{thm91}
Let $X$ be a normal quasi-projective variety and 
let $B$ be a boundary $\mathbb R$-divisor on $X$ such 
that $K_X+B$ is $\mathbb R$-Cartier. 
In this case, we can construct 
a projective birational 
morphism $f:Y\to X$ from a normal quasi-projective 
variety $Y$ with the following properties. 
\begin{itemize}
\item[(i)] $Y$ is $\mathbb Q$-factorial. 
\item[(ii)] $a(E, X, B)\leq -1$ for every  
$f$-exceptional divisor $E$ on $Y$. 
\item[(iii)] We put 
$$
B_Y=f^{-1}_*B+\sum _{E: {\text{$f$-exceptional}}}E. 
$$ 
Then $(Y, B_Y)$ is dlt and  
$$
K_Y+B_Y=f^*(K_X+B)+\sum _{a(E, X, B)<-1}(a(E, X, B)+1)E. 
$$
In particular, if $(X, B)$ is lc, then 
$K_Y+B_Y=f^*(K_X+B)$. 
Moreover, if $(X, B)$ is dlt, then we can make $f$ small, 
that is, $f$ is an isomorphism in codimension one. 
\end{itemize}
\end{thm}
\begin{proof}
Let $\pi:V\to X$ be a resolution 
such that 
$\pi^{-1}_*B\cup \Exc (\pi)$ has a simple normal crossing support. 
We can assume that $\pi$ is a composite of blow-ups of centers 
of codimension at least two. 
Then there exists an effective 
$\pi$-exceptional Cartier divisor $C$ on $V$ such 
that $-C$ is $\pi$-ample. 
We put 
$$
F=\sum _{\begin{matrix}\scriptstyle{a(E, X, B)>-1, }\\ 
\scriptstyle {E: {\text{$\pi$-exceptional}}}\end{matrix}}E
$$ 
and 
$$
E^+=-\sum _{a(E, X, B)\leq -1}a(E, X, B)E. 
$$ 
We note that $E^+$ is not necessarily $\pi$-exceptional. 
We put $E=\Supp E^+$. 
We note that $E^+-E$ is $\pi$-exceptional. 

Let $H$ be a sufficiently ample Cartier divisor on $X$. 
We choose $0<\epsilon, \nu, \mu\ll 1$ and 
note that 
\begin{align*}\label{spade}
&E+(1-\nu)F+\mu(-C+\pi^*H)\tag{$\spadesuit$}
\\ &=(1-\epsilon \mu)E+(1-\nu)F+
\mu(\epsilon E-C+\pi^*H). 
\end{align*} 
Since $-C+\pi^*H$ and $\epsilon E-C+\pi^*H$ are ample, 
we can take effective $\mathbb Q$-divisors 
$H_1$ and $H_2$ with small coefficients 
such that 
$E+F+\pi^*B+H_1+H_2$ has a simple normal crossing 
support and 
that $-C+\pi^*H\sim _{\mathbb Q}H_1$, $\epsilon E-C+\pi^*H\sim _{\mathbb Q}H_2$. 
Then $(V, (1-\epsilon\mu)E+(1-\nu)F+\pi^{-1}_*B^{<1}+\mu H_2)$ 
is klt. 
By Theorem \ref{thm92}, 
it has a log terminal model $f:Y\to X$. 
By the above equation (\ref{spade}), this is also a relative minimal model 
of the pair $(V, E+(1-\nu)F+\pi^{-1}_*B^{<1}+\mu H_1)$, 
which is therefore dlt. 

For any divisor $G$ on $V$ appearing above, 
let $G'$ denote its transform on $Y$. 
By the above construction, 
$$
N=K_Y+(1-\epsilon \mu)E'+(1-\nu)F'+f^{-1}_*B^{<1}+\mu H'_2
$$ 
is $f$-nef and $K_Y+\overline B=f^*(K_X+B)$ is 
$\mathbb R$-linearly $f$-trivial. 
We put 
$$
D=\overline B-E'-(1-\nu)F'-f^{-1}_*B^{<1}+\mu C'. 
$$ 
Then 
\begin{align*}
-D&\sim _{\mathbb R, f}N-(K_Y+\overline B)\\ &=
-\overline B+(1-\epsilon \mu)E'+(1-\nu)F'
+f^{-1}_*B^{<1}+\mu H'_2, 
\end{align*}
hence it is $f$-nef. 
Since $f_*D=0$, 
we see that $D$ is effective by the negativity lemma (cf.~Lemma \ref{lem211}). 

Every divisor in $F$ has a negative coefficient in 
$$
\widetilde B-E-(1-\nu)F-\pi^{-1}_*B^{<1}+\mu C,
$$ 
where $K_V+\widetilde B=\pi^*(K_X+B)$. 
Therefore, $F$ is contracted on $Y$. 
So, every $f$-exceptional divisor has discrepancy $\leq -1$. By the 
above construction, $(Y, E'+f^{-1}_*B^{<1}+\mu H'_1)$ is 
dlt since $F'=0$. 
Therefore, $(Y, E'+f^{-1}_*B^{<1})$ is also dlt. 
This means that $(Y, B_Y)$ is 
dlt because $B_Y=E'+\sum f^{-1}_*B^{<1}$. 

When $(X, B)$ is dlt, we can assume that $E^+=\pi^{-1}_*B^{=1}$ by the 
definition of dlt pairs. 
Therefore, we can make $f$ small. 
\end{proof}

The following technical statement seems to be very useful for future 
studies (cf.~\cite{gongyo}), though 
we do not use it in this 
paper. 

\begin{thm}\label{gong}
Let $X$ be a normal quasi-projective variety and 
let $B$ be an effective $\mathbb R$-divisor on $X$ such 
that $(X, B)$ is lc. 
In this case, we can construct 
a projective birational 
morphism $f:Y\to X$ from a normal quasi-projective 
variety $Y$ with the following properties. 
\begin{itemize}
\item[(i)] $Y$ is $\mathbb Q$-factorial. 
\item[(ii)] $a(E, X, B)= -1$ for every  
$f$-exceptional divisor $E$ on $Y$. 
\item[(iii)] We put 
$$
B_Y=f^{-1}_*B+\sum _{E: {\text{$f$-exceptional}}}E. 
$$ 
Then $(Y, B_Y)$ is dlt and  
$
K_Y+B_Y=f^*(K_X+B)
$.
\item[(iv)] Let $\{C_i\}$ be {\em{any}} set of lc centers of $(X, B)$. We put 
$W=\bigcup C_i$ with the reduced scheme structure. 
Let $S$ be the union of the irreducible components 
of $B^{=1}_Y$ which are mapped into $W$ by $f$. Then 
$f_*\mathcal O_S\simeq \mathcal O_W$. 
\end{itemize}
\end{thm}
\begin{proof}
Let $\pi:V\to X$ be a resolution such that 
\begin{itemize}
\item[(1)] $\pi^{-1}(C)$ is a simple normal crossing divisor on $V$ for every 
lc center $C$ of $(X, B)$, and 
\item[(2)] $\pi^{-1}_*B\cup \Exc(\pi)\cup \pi^{-1}(\Nklt (X, B))$ has a simple normal crossing 
support. 
\end{itemize}
We apply the arguments in the proof of Theorem \ref{thm91}. From now on, 
we use the same notation as in the proof of Theorem \ref{thm91}. 
In this case, we have 
$$
E=\Supp E^+=E^+. 
$$ 
When we construct $f:Y\to X$, we can run the log minimal model program 
with scaling with respect to 
\begin{align*}
&K_V+E+(1-\nu)F+\pi^{-1}_*B^{<1}+\mu H_1\\ &\sim _{\mathbb R}K_V+(1-\varepsilon \mu)E+(1-\nu)F
+\pi^{-1}_*B^{<1}+\mu H_2 
\end{align*}
(cf.~\cite{bchm}). 
So, we can assume that $\varphi:V\dashrightarrow Y$ is a 
composition of $(K_V+E+(1-\nu)F+\pi^{-1}_*B^{<1}+\mu H_1)$-negative 
divisorial contractions and log flips. 
Let $\Sigma$ be an lc center of $(Y, B_Y)$. Then it is 
also an lc center of $(Y, B_Y+\mu H'_1)$. By the negativity lemma (cf.~Lemma \ref{lem211}), 
$\varphi:V\dashrightarrow Y$ is an isomorphism around the generic point of $\Sigma$. 
Therefore, if $f(\Sigma)\subset W$, 
then $\Sigma\subset S$ by the conditions (1) and (2) for 
$\pi:V\to X$. 
This means that no lc centers of $(Y, B_Y-S)$ are mapped into $W$ by $f$. 
Let $g:Z\to Y$ be a resolution such that 
\begin{itemize}
\item[(a)] $K_Z+B_Z=g^*(K_Y+B_Y)$, 
\item[(b)] $\Supp B_Y$ is a simple normal crossing divisor, and 
\item[(c)] $g$ is an isomorphism over the generic point of any lc center 
of $(Y, B_Y)$. 
\end{itemize}
Let $S_Z$ be the strict transform of $S$ on $Z$. 
We consider the following short exact sequence 
\begin{align*}\label{hhh}
0\to \mathcal O_Z(\ulcorner -(B^{<1}_Z)\urcorner -S_Z)&\to \mathcal O_Z(\ulcorner 
-(B^{<1}_Z)\urcorner )\tag{$\heartsuit$}\\ 
&\to \mathcal O_{S_Z}(\ulcorner -(B^{<1}_Z)\urcorner)\to 0. 
\end{align*}
We note that 
\begin{align*}
\ulcorner -(B^{<1}_Z)\urcorner -S_Z-(K_Z+\{B_Z\}+B^{=1}_Z-S_Z)\sim _{\mathbb R}-h^*(K_X+B), 
\end{align*}
where $h=f\circ g$. 
Then we obtain 
\begin{align*}
0&\to h_*\mathcal O_Z(\ulcorner -(B^{<1}_Z)\urcorner -S_Z)\to h_*\mathcal O_Z(\ulcorner 
-(B^{<1}_Z)\urcorner )\to 
h_*\mathcal O_{S_Z}(\ulcorner -(B^{<1}_Z)\urcorner)\\ &\overset {\delta}\to 
R^1h_*\mathcal O_Z(\ulcorner -(B^{<1}_Z)\urcorner-S_Z)\to \cdots. 
\end{align*}
Every associated prime of $R^1h_*\mathcal O_Z(\ulcorner -(B^{<1}_Z)\urcorner-S_Z)$ 
is the generic point of the $h$-image of some stratum of 
$(Z, \{B_Z\}+B^{=1}_Z-S_Z)$ by Theorem \ref{thm53} (i) and 
no lc centers of $(Z, \{B_Z\}+B^{=1}_Z-S_Z)$ are mapped into $W$ by $h$. 
Therefore, $\delta$ is a zero map. 
Thus, we obtain 
$$
0\to \mathcal I_W\to \mathcal O_X\to h_*\mathcal O_{S_Z}(\ulcorner -(B^{<1}_Z)\urcorner)\to 0 
$$ 
and $\mathcal O_W\simeq h_*\mathcal O_{S_Z}\simeq h_*\mathcal O_{S_Z}(\ulcorner 
-(B^{<1}_Z)\urcorner)$ (cf.~Proposition \ref{lem72}), 
where $\mathcal I_W$ is the defining ideal sheaf of $W$. Here, 
we used the fact that $\ulcorner -(B^{<1}_Z)\urcorner$ is effective 
and $h$-exceptional. 
By applying $g_*$ to (\ref{hhh}), we obtain 
$$
0\to \mathcal I_S\to \mathcal O_Y\to g_*\mathcal O_{S_Z}(\ulcorner -(B^{<1}_Z)\urcorner)\to 0
$$ 
and 
$\mathcal O_S\simeq g_*\mathcal O_{S_Z}\simeq g_*\mathcal O_{S_Z}(\ulcorner 
-(B^{<1}_Z)\urcorner)$ (cf.~Proposition \ref{lem72}), 
where $\mathcal I_S\simeq \mathcal O_Y(-S)$ is the defining ideal sheaf of $S$.
We note that 
$$R^1g_*\mathcal O_Z(\ulcorner -(B^{<1}_Z)\urcorner-S_Z)=0$$ by 
Theorem \ref{thm53} (i) since $g$ is an isomorphism 
at the generic point of any stratum of $(Z, \{B_Z\}+B^{=1}_Z-S_Z)$ and 
that $\ulcorner -(B^{<1}_Z)\urcorner$ is effective and $g$-exceptional. 
Therefore, $\mathcal O_W\simeq h_*\mathcal O_{S_Z}\simeq f_*g_*\mathcal O_{S_Z}
\simeq f_*\mathcal O_S$. 
\end{proof}

\section{Vanishing theorem for minimal lc centers}\label{sec10} 

In this section, we prove a vanishing theorem 
on minimal lc centers. 
It is very powerful and will play crucial roles in 
the proof of Theorem \ref{thm111}. 
We note that a key point of Theorem \ref{thm101} 
is in its formulation which
is best suited for our subsequent applications.

\begin{thm}[Vanishing theorem for minimal lc centers]\label{thm101} 
Let $X$ be a normal variety and 
let $B$ be an effective 
$\mathbb R$-divisor on $X$ such that 
$K_X+B$ is $\mathbb R$-Cartier. 
Let $W$ be a minimal lc center 
of $(X, B)$ such that 
$W$ is disjoint from $\Nlc(X, B)$. Let 
$\pi:X\to S$ be a projective morphism 
onto a variety $S$. 
Let $D$ be a Cartier divisor 
on $W$ such that 
$D-(K_X+B)|_W$ is $\pi$-ample. 
Then 
$$
R^i\pi_*\mathcal O_W(D)=0 
$$ 
for every $i>0$. 
\end{thm}
\begin{proof} 
Without loss of generality, we can assume that 
$S$ is quasi-projective. 
We shrink $X$ around $W$ and assume that 
$(X, B)$ is log canonical. 
By Theorem \ref{thm91}, 
we can make a projective 
birational morphism $f:Y\to X$ such that 
$K_Y+B_Y=f^*(K_X+B)$ and $(Y, B_Y)$ is dlt. 
We take an lc center $V$ of $(Y, B_Y)$ such that 
$f(V)=W$ and put $K_V+B_V=(K_Y+B_Y)|_V$. 
Then $(V, B_V)$ is dlt by Lemma \ref{lem92} and 
$K_V+B_V\sim _{\mathbb R}f^*((K_X+B)|_W)$. 
Let $g:Z\to V$ be a resolution such that 
$K_Z+B_Z=g^*(K_V+B_V)$ and $\Supp B_Z$ is simple normal 
crossing. Then we have 
$K_Z+B_Z\sim _{\mathbb R}h^*((K_X+B)|_W)$, where 
$h=f\circ g$. 
Since 
$$h^*(D-(K_X+B)|_W)\sim _{\mathbb R}h^*D+\ulcorner 
-(B^{<1}_Z)\urcorner-(K_Z+B^{=1}_Z+\{B_Z\}), $$
we obtain 
$$R^i\pi_*h_*\mathcal O_Z(h^*D+\ulcorner -(B^{<1}_Z)\urcorner)=0$$ for 
every $i>0$ by Theorem \ref{thm53} (ii). 
We note that 
$$h_*\mathcal O_Z(h^*D+\ulcorner -(B^{<1}_Z)\urcorner)\simeq 
f_*\mathcal O_V(f^*D)$$ by the projection formula since 
$\ulcorner -(B^{<1}_Z)\urcorner$ is effective and $g$-exceptional. 
We note that 
$\mathcal O_W(D)$ is a direct summand of $f_*\mathcal O_V(f^*D)\simeq 
\mathcal O_W(D)\otimes f_*\mathcal O_V$ since $W$ is normal 
(cf.~Theorem \ref{thm81} (4)). 
Therefore, we have $R^i\pi_* \mathcal O_W(D)=0$ for every $i>0$. 
\end{proof}

We close this section with a very important remark. 
 
\begin{rem}
The short proof of 
Theorem \ref{thm101} 
given in this section depends on
Theorem \ref{thm53} (ii), 
Theorem \ref{thm81} (4), and
Theorem \ref{thm91} 
which is a corollary to \cite{bchm}. 
However Theorem \ref{thm101}, a special case of \cite[Theorem 4.4]{ambro},
is independent of \cite{bchm} since it
can be proved without using Theorem \ref{thm91}. 
We refer the reader to \cite[Theorem 3.39]{book} as
for the independent proof of Theorem \ref{thm101}, which 
heavily depends on the theory of 
mixed Hodge structures on compact support cohomology groups 
of {\em{reducible}} varieties (cf.~\cite[Chapter 2]{book}).  
\end{rem}

\section{Non-vanishing theorem}\label{sec11}

In this section, we prove the non-vanishing theorem, 
which is a generalization of 
the main theorem of \cite{non-va}. 
In \cite{ambro}, Ambro does not discuss any generalization of Shokurov's 
non-vanishing theorem. Therefore, the result in this section 
is one of the main differences between 
the theory of quasi-log varieties and 
our new framework.
 
\begin{thm}[Non-vanishing theorem]\label{thm111}  
Let $X$ be a normal variety and 
let $B$ be an effective $\mathbb R$-divisor 
on $X$ such that 
$K_X+B$ is $\mathbb R$-Cartier. 
Let $\pi:X\to S$ be a projective 
morphism onto a variety $S$ and let $L$ be a $\pi$-nef Cartier 
divisor on $X$. 
Assume that 
\begin{itemize}
\item[(i)] $aL-(K_X+B)$ is $\pi$-ample for some 
real number $a>0$, and 
\item[(ii)] $\mathcal O_{\Nlc(X, B)}(mL)$ is $\pi|_{\Nlc(X, B)}$-generated 
for $m\gg 0$. 
\end{itemize}
Then the relative base locus $\Bs_{\pi} |mL|$ contains 
no lc centers of $(X, B)$ and 
is disjoint from $\Nlc(X, B)$ for 
$m\gg 0$. 
\end{thm}
\begin{proof}
Without loss of generality, we can assume that 
$S$ is affine. 
\setcounter{step}{0}
\begin{step}\label{step1}
In this step, we will prove that 
$\mathcal O_X(mL)$ is $\pi$-generated 
on an open neighborhood of $\Nlc (X, B)$ for $m \gg 0$. 

By the assumption, $\pi^*\pi_*\mathcal O_{\Nlc (X, B)}
(mL)\to \mathcal O_{\Nlc (X, B)}(mL)$ is 
surjective for $m \gg 0$. 
On the other hand, $\pi_*\mathcal O_X(mL)\to \pi_*\mathcal O_{\Nlc (X, B)}(mL)$ 
is surjective 
for $m \geq a$ since 
$$R^1\pi_*(\mathcal J_{NLC}(X, B)\otimes \mathcal O_X(mL))=0$$ for 
$m\geq a$ by Theorem \ref{thm71}. 
Therefore, for every large integer $m$, 
$\pi^*\pi_*\mathcal O_X(mL)\to \mathcal O_X(mL)$ is 
surjective on an open neighborhood of $\Nlc (X, B)$. 
See the following commutative diagram. 
$$
\xymatrix{
\pi^*\pi_*\mathcal O_X(mL)
\ar[d]\ar[r]& 
\pi^*\pi_*\mathcal O_{\Nlc (X, B)}(mL)
\ar[d]\ar[r]&0\\
\mathcal O_X(mL)
\ar[r]&
\mathcal O_{\Nlc(X, B)}(mL)
\ar[r]&0, 
}
$$
\end{step}
Let $W$ be a minimal lc center of $(X, B)$. 
Then it is sufficient to see that $W$ is not contained 
in $\Bs|mL|$ for $m\gg 0$. 
\begin{step}\label{step2} 
If $W\cap \Nlc (X, B)\ne \emptyset$, then 
$\Bs |mL|$ does not contain $W$ by Step \ref{step1}. 
So, from now on, we can assume that 
$W\cap \Nlc (X, B)=\emptyset$. 
\end{step}
\begin{step}\label{step3}
We assume that $L|_{W_\eta}$ is numerically trivial, where 
$W_\eta$ is the generic fiber of $W\to \pi(W)$. 
In this case, 
\begin{align*}
&h^0(W_\eta, \mathcal O_{W_\eta}(L))=\chi 
(W_\eta, \mathcal O_{W_\eta}(L))
\\&=\chi (W_\eta, \mathcal O_{W_\eta})=h^0(W_\eta, \mathcal O_{W_\eta})>0
\end{align*} 
by 
\cite[Chapter II \S2 Theorem 1]{kleiman} 
and the vanishing theorem:~Theorem \ref{thm101}. 
On the other hand, 
$$
\pi_*\mathcal O_X(mL)\to \pi_*\mathcal O_W(mL)\oplus 
\pi_*\mathcal O_{\Nlc (X, B)}(mL)
$$ 
is surjective for every $m \geq a$ by Theorem \ref{thm71}. 
In particular, the restriction map 
$\pi_*\mathcal O_X(mL)\to \pi_*\mathcal O_W(mL)$ is surjective 
for every $m \geq a$. 
Thus, $\Bs|mL|$ does not contain $W$ for every $m\geq a$. 
\end{step}
\begin{step}\label{step4} 
We assume that $L|_{W_{\eta}}$ is not numerically trivial. 
We take a general subvariety $V$ of $W$ such that 
$V\to \pi(W)$ is generically finite. 
If $l$ is a positive large integer, 
then we can write 
$$
lL-(K_X+B)=N_1+a_2N_2+\cdots +a_k N_k 
$$
with the following properties. 
\begin{itemize}
\item[(a)] $N_1$ is a $\pi$-ample 
$\mathbb Q$-Cartier $\mathbb Q$-divisor 
on $X$ such that 
$$((N_1|_W)|_F)^{\dim F}> d(\codim _WV)^{\dim F}, 
$$ 
where $d$ is the 
mapping degree of 
$V\to \pi(W)$ and 
$F$ is a general fiber of 
$W\to \pi(W)$. 
\item[(b)] $a_i$ is a positive 
real number and $N_i$ is a $\pi$-very ample 
Cartier divisor on $X$ for every $i\geq 2$. 
\end{itemize}
By Lemma \ref{lem112}, we can 
find an effective $\mathbb Q$-divisor 
$D_1$ on $W$ such that $D_1\sim _{\mathbb Q}N_1|_W$ with 
$\mult _VD_1>\codim _W V$. If 
$b$ is sufficiently large and divisible, 
then $bD_1\sim bN_1|_W$, 
$\mathcal I_W\otimes \mathcal O_X(bN_1)$ is $\pi$-generated, 
and $R^1\pi_*(\mathcal I_W\otimes \mathcal O_X(bN_1))=0$ since 
$N_1$ is $\pi$-ample, where 
$\mathcal I_W$ is the defining ideal sheaf of 
$W$. 
By using the following 
short exact sequence 
$$
0\to \pi_*(\mathcal I_W\otimes \mathcal O_X(bN_1))
\to \pi_*\mathcal O_X(bN_1)\to 
\pi_*\mathcal O_W(bN_1)\to 0, 
$$
we can find an effective $\mathbb Q$-divisor 
$M_1$ on $X$ with the following 
properties. 
\begin{itemize}
\item[(1)] $M_1|_W=D_1$. 
\item[(2)] $M_1\sim_{\mathbb Q}N_1$. 
\item[(3)] $(X, B+M_1)$ is 
lc outside $W\cup \Nlc (X, B)$. 
\item[(4)] $\mathcal J_{NLC}(X, B+M_1)=\mathcal J_{NLC}(X, B)$ outside 
$W$. 
\end{itemize}
Let $M_i$ be a general member of $|N_i|$ for every $i\geq 2$. 
We put $M=M_1+a_2M_2+\cdots +a_k M_k$. 
Then we have 
\begin{itemize}
\item[(i)] $M|_W\geq D_1$. 
\item[(ii)] $M\sim _{\mathbb R}lL-(K_X+B)$. 
\item[(iii)] $(X, B+M)$ is lc outside 
$W\cup \Nlc (X, B)$. 
\item[(iv)] $\mathcal J_{NLC}(X, B+M)=\mathcal J_{NLC}(X, B)$ outside 
$W$. 
\end{itemize}
We take the log canonical threshold 
$c$ of $(X, B)$ with respect 
to $M$ outside $\Nlc (X, B)$. 
By the above construction, we have $0<c<1$. 
More precisely, we see $0<c$ since $M$ contains no lc centers of $(X, B)$. The 
inequality $c<1$ follows from the fact that $M|_W\geq D_1$ and $\mult _V D_1>
\codim _WV$. 
We note that 
$$
(a-ac+cl)L-(K_X+B+cM)\sim _{\mathbb R}
(1-c)(aL-(K_X+B))
$$
is $\pi$-ample. Moreover, 
we can find a smaller lc center $W'$ of $(X, B+cM)$ contained 
in $W$ (cf.~Theorem \ref{thm81} (2)). 
Therefore, we replace 
$(X, B)$ with 
$(X, B+cM)$, $a$ with 
$a-ac+cl$, and consider the new lc center $W'$. 
By repeating this process, we reach the situation 
where $L|_{W_\eta}$ is numerically trivial. 
\end{step}
Anyway, we proved that $\Bs|mL|$ contains no lc centers 
of $(X, B)$ for $m \gg 0$. 
\end{proof}

The following lemma is a relative 
version of Shokurov's concentration 
method. We used it in the proof of Theorem \ref{thm111}. 

\begin{lem}\label{lem112}
Let $f:Y\to Z$ be a projective 
morphism from a normal variety $Y$ onto an affine 
variety $Z$. 
Let $V$ be a general 
closed subvariety of $Y$ such that 
$f:V\to Z$ 
is generically finite. 
Let $M$ be an $f$-ample $\mathbb R$-divisor 
on $Y$. Assume that 
$$
(M|_F)^d>km^d, 
$$ 
where $F$ is a general fiber of $f:Y\to Z$, $d=\dim F$, and $k$ is the mapping 
degree of $f:V\to Z$. Then 
we can find an effective 
$\mathbb R$-divisor $D$ on $Y$ such that 
$$
D\sim _{\mathbb R}M 
$$ 
and that 
$\mult _V D>m$. 
If $M$ is a $\mathbb Q$-divisor, then we can make 
$D$ a $\mathbb Q$-divisor with $D\sim _{\mathbb Q}M$. 
\end{lem}
\begin{proof}
We can write 
$$
M=M_1+a_2M_2+\cdots +a_lM_l, 
$$ 
where $M_1$ is an $f$-ample $\mathbb Q$-Cartier 
$\mathbb Q$-divisor such that 
$(M_1|_F)^d>km^d$, 
$a_i$ is a positive real number, and 
$M_i$ is an $f$-ample Cartier 
divisor for every $i$. 
If $M$ is a $\mathbb Q$-divisor, then we 
can assume that $l=2$ and $a_2$ is rational.  
Let $\mathcal I_V$ be the defining ideal sheaf of 
$V$ on $Y$. We consider the following 
exact sequence
\begin{align*}
0&\to f_*(\mathcal O_Y(pM_1)\otimes \mathcal I^{pm}_V)\to 
f_*\mathcal O_Y(pM_1)\\ &\to f_*(\mathcal O_Y(pM_1)\otimes 
\mathcal O_Y/\mathcal I^{pm}_V)\to \cdots   
\end{align*}
for a sufficiently large and divisible 
integer $p$. 
By restricting the above sequence to 
a general fiber $F$ of $f$, we can check 
that the rank of $f_*\mathcal O_Y(pM_1)$ is greater than 
that of $f_*(\mathcal O_Y(pM_1)\otimes 
\mathcal O_Y/\mathcal I^{pm}_V)$ by the usual 
estimates (see Lemma \ref{n-lem123} below). 
Therefore, $f_*(\mathcal O_Y(pM_1)\otimes \mathcal 
I^{pm}_V)\ne 0$. 
Let $D_1$ be a member of 
$$
H^0(Z, f_*(\mathcal O_Y(pM_1)\otimes \mathcal I^{pm}_V))=
H^0(Y, \mathcal O_Y(pM_1)\otimes \mathcal I^{pm}_V)
$$ 
and let $D_i$ be an effective $\mathbb Q$-Cartier 
$\mathbb Q$-divisor such that 
$D_i\sim _{\mathbb Q}M_i$ for $i\geq 2$. 
We can take $D_2$ with $\mult _VD_2>0$. 
Then $D=(1/p)D_1+a_2D_2+\cdots 
+a_l D_l$ satisfies 
the desired properties. 
\end{proof}

We close this section with the following well-known 
lemma. The proof is obvious. 

\begin{lem}\label{n-lem123}
Let $X$ be a normal projective variety with 
$\dim X=d$ and 
let $A$ be an ample 
$\mathbb Q$-divisor on $X$ such that 
$rA$ is Cartier for some positive integer $r$. 
Then 
\begin{align*}
h^0(X, \mathcal O_X(trA))&=\chi (X, \mathcal O_X(trA))
\\ &= \frac{(trA)^d}{d!}+{\text{$($lower terms in $t$$)$}} 
\end{align*} 
by the Riemann-Roch formula and the 
Serre vanishing theorem for $t \gg 0$. 

Let $P\in X$ be a smooth point. 
Then 
\begin{align*}
\dim _{\mathbb C}\mathcal O_X/m^{\alpha}_P&=\begin{pmatrix}
\alpha -1 +d \\ d
\end{pmatrix} 
\\&=\frac{\alpha^d}{d!}+{\text{$($lower terms in $\alpha$$)$}}  
\end{align*}
for all $\alpha\geq 1$, where $m_P$ is the maximal ideal associated to $P$. 
\end{lem}

\section{Base point free theorem}\label{sec12}
The base point free theorem is one of the most important 
theorems in the log minimal model program. 
Since we have already established 
the non-vanishing theorem (Theorem \ref{thm111}) in our framework, 
there are no difficulties in obtaining the base point free theorem (Theorem \ref{thm121}). 
Our approach is simpler than \cite{ambro}, 
though Theorem \ref{thm121} is a special 
case of the base point free theorem for quasi-log varieties obtained by 
Ambro (cf.~\cite[Theorem 5.1]{ambro} and \cite[Theorem 3.66]{book}). 
Indeed in the approach of \cite{ambro} it is necessary
to treat reducible non-equidimensional 
quasi-log varieties even for the proof of the base point free theorem for log canonical 
pairs.

\begin{thm}[Base point free theorem]\label{thm121} 
Let $X$ be a normal variety and let 
$B$ be an effective $\mathbb R$-divisor 
on $X$ such that 
$K_X+B$ is $\mathbb R$-Cartier. 
Let $\pi:X\to S$ be a projective 
morphism onto a variety $S$ and let $L$ be a $\pi$-nef Cartier 
divisor on $X$. 
Assume that 
\begin{itemize}
\item[(i)] $aL-(K_X+B)$ is $\pi$-ample for some 
real number $a>0$, and 
\item[(ii)] $\mathcal O_{\Nlc(X, B)}(mL)$ is $\pi|_{\Nlc(X, B)}$-generated 
for $m\gg 0$. 
\end{itemize}
Then $\mathcal O_X(mL)$ is $\pi$-generated 
for $m \gg 0$. 
\end{thm}

We will prove the base point free theorem for $\mathbb R$-divisors 
in Section \ref{sec-n17} as an application of the 
cone theorem:~Theorem \ref{thm144}. 

\begin{proof}
We can assume that 
$S$ is affine. 
\setcounter{step}{0}
\begin{step}\label{st1}
We assume that 
$(X, B)$ is klt and 
that 
$L_{\eta}$ is numerically trivial, 
where $L_{\eta}=L|_{X_{\eta}}$ and $X_\eta$ is the generic fiber of 
$\pi:X\to S$. Then 
we have 
\begin{align*}
&h^0(X_\eta, \mathcal O_{X_\eta}(L_\eta))
=\chi (X_\eta, \mathcal O_{X_\eta}(L_\eta))
\\ &=\chi (X_\eta, \mathcal O_{X_\eta})=h^0(X_\eta, \mathcal O_{X_\eta})>0
\end{align*}  
by 
\cite[Chapter II \S2 Theorem 1]{kleiman} and 
the vanishing theorem. Here, 
the Kawamata--Viehweg vanishing theorem is 
sufficient. Therefore, 
$|L|\ne \emptyset$. Let $D$ be a member of $|L|$. If $D=0$, then 
it is obvious 
that $|mL|$ is free for every $m$. 
Thus, we can assume that 
$D\ne 0$. 
Let $c$ be the log canonical threshold 
of $(X, B)$ with respect to 
$D$. We replace 
$(X, B)$ with $(X, B+cD)$, $a$ with 
$a+c$. 
Then we can assume that $(X, B)$ is lc but not klt. 
This case will be treated in Step \ref{st3}. 
\end{step}
\begin{step}\label{st2}
We assume that 
$(X, B)$ is klt and 
that $L_\eta$ is not numerically trivial. 
We take a general subvariety $V$ on $X$ such that 
$\pi:V\to S$ is generically finite. 
By Lemma \ref{lem112}, we can find 
an effective $\mathbb R$-divisor $D$ on $X$ such 
that 
$$
D\sim _{\mathbb R}lL-(K_X+B) 
$$ 
for some large $l$ and that 
$\mult _V D>\codim _XV$. Let $c$ be the 
log canonical threshold 
of $(X, B)$ with respect to 
$D$. 
By the above construction, we obtain $0<c<1$. We replace 
$(X, B)$ with 
$(X, B+cD)$, 
$a$ with $a-ac+cl$ and 
can assume that 
$(X, B)$ is lc but not klt. 
We note that 
$$
(a-ac+cl)L-(K_X+B+cD)\sim_{\mathbb R}(1-c)(aL-(K_X+B)). 
$$ 
So, the problem is reduced to the case when $(X, B)$ is lc but not klt. 
It will be treated in Step \ref{st3}. 
\end{step}
\begin{step}\label{st3}
We assume that 
$(X, B)$ is not klt. 
Let $p$ be a prime integer. 
We will prove that 
$\Bs|p^mL|=\emptyset$ for some positive integer $m$. 

By Theorem \ref{thm111}, 
$|p^{m_1}L|\ne\emptyset$ for some positive integer $m_1$. 
If $\Bs |p^{m_1}L|=\emptyset$, then 
there are nothing to prove. 
So, we can assume that $\Bs|p^{m_1}L|\ne\emptyset$. 
We take general members 
$D_1, \cdots, D_{n+1}\in |p^{m_1}L|$, 
where $n=\dim X$. 
Since $D_1, \cdots, D_{n+1}$ are general, 
$(X, B+D_1+\cdots+D_{n+1})$ is lc outside 
$\Bs|p^{m_1}L|\cup \Nlc (X, B)$. 
It is easy to see that $(X, B+D)$, 
where $D=D_1+\cdots +D_{n+1}$, is 
not lc at the generic point of every irreducible 
component of $\Bs|p^{m_1}L|$ (see Lemma \ref{lem134} below). 
Let $c$ be the log canonical threshold 
of $(X, B)$ with respect to 
$D$ outside 
$\Nlc(X, B)$. 
Then $(X, B+cD)$ is lc but not klt outside 
$\Nlc(X, B)$, $0<c<1$, 
and $\mathcal J_{NLC}(X, B+cD)=\mathcal J_{NLC}(X, B)$ 
(see Proposition \ref{prop63}). 
We note that 
\begin{align*}
(c(n+1)p^{m_1}+a)L-(K_X+B+cD)\sim _{\mathbb R}aL-(K_X+B) 
\end{align*} 
is $f$-ample. 
By the construction, there exists an lc center 
of $(X, B+cD)$ contained in $\Bs|p^{m_1}L|$. 
By Theorem \ref{thm111}, we can find 
$m_2>m_1$ such that 
$\Bs|p^{m_2}L|\subsetneq \Bs|p^{m_1}L|$. By the noetherian 
induction, there exists $m$ such that 
$\Bs|p^{m}L|=\emptyset$. 
\end{step}
\begin{step}\label{st4} 
Let $p'$ be a prime integer such that $p'\ne p$. 
Then, by Step \ref{st3} again, we can find a positive integer $m'$ such that 
$\Bs |p'^{m'}L|=\emptyset$. 
So, there exists a positive integer 
$m_0$ such that 
$|kL|$ is free for every $k\geq m_0$ 
by $\Bs|p^m L|=\emptyset$ and $\Bs|p'^{m'}L|=\emptyset$. 
\end{step}
This completes the proof. 
\end{proof}

We close this section with the following lemma. 
We used it in the proof of Theorem \ref{thm121}. 

\begin{lem}\label{lem134}
Let $X$ be a normal variety and let $B$ be an effective 
$\mathbb R$-divisor 
on $X$ such that 
$K_X+B$ is $\mathbb R$-Cartier. 
Let $P$ be a closed point of $X$ and let $P\in D_i$ 
be a Cartier divisor for every $i$. 
If $(X, B+\sum _{i=1}^{k}D_i)$ is log canonical at $P$, 
then $k\leq \dim X$. 
\end{lem}
\begin{proof}
The proof is by the induction 
on $\dim X$. 
The assertion is clear if $\dim X=1$. 
We put $S=D_1$. 
Let $\nu:S^{\nu}\to S$ be the normalization 
and let $B_{S^{\nu}}$ be the {\em{different}} 
of $(X, S+B)$ on $S^{\nu}$ (see Section \ref{sec13} below). 
So, we have $K_{S^{\nu}}+B_{S^{\nu}}=\nu^*
(K_X+S+B)$. 
Since $(X, B+S+\sum _{i=2}^kD_i)$ is log canonical 
at $P$, 
$(S^{\nu}, B_{S^{\nu}}+\sum _{i=2}^{k}\nu^*D_i)$ is log 
canonical 
at $Q\in \nu^{-1}(P)$. Thus, 
$k-1\leq \dim S^{\nu}$ by induction. 
This means that $k\leq \dim X$. 
\end{proof}

\section{Shokurov's differents}\label{sec13} 
Let us recall the definition and basic 
properties of Shokurov's {\em{differents}} 
following \cite[\S 3]{shokurov2} and 
\cite[9.2.1]{ambro2}.  

\begin{say}
Let $X$ be a normal variety and 
let $S+B$ be an $\mathbb R$-divisor on $X$ such that 
$K_X+S+B$ is $\mathbb R$-Cartier. 
Assume that 
$S$ is reduced and that 
$S$ and $B$ have no common irreducible 
components. 
Let $f:Y\to X$ be a resolution 
such that 
$$
K_Y+S_Y+B_Y=f^*(K_X+S+B)
$$ 
and $\Supp (S_Y+B_Y)$ is simple normal crossing and 
$S_Y$ is smooth, where 
$S_Y$ is the strict transform of $S$ on $Y$. 
Let $\nu:S^{\nu}\to S$ be the normalization. 
Then $f:S_Y\to S$ can be decomposed as 
$$
f:S_Y\overset{\pi}\longrightarrow S^{\nu}\overset{\nu}\longrightarrow S. 
$$
We define $B_{S_Y}=B_Y|_{S_Y}$. Then 
we obtain 
$$
(K_Y+S_Y+B_Y)|_{S_Y}=K_{S_Y}+B_{S_Y} 
$$ 
by adjunction. We put $B_{S^{\nu}}=\pi_*B_{S_Y}$. 
Then we have  
$$
K_{S^{\nu}}+B_{S^{\nu}}=\nu^*(K_X+S+B). 
$$ 
The $\mathbb R$-divisor $B_{S^{\nu}}$ on $S^{\nu}$ is called 
the {\em{different}} of $(X, S+B)$ on $S^{\nu}$. We can easily 
check that $B_{S^\nu}$ is independent of the resolution 
$f:Y\to X$. So, 
$B_{S^\nu}$ is a well-defined $\mathbb R$-divisor on $S^{\nu}$. 
We can check the following properties. 
\begin{itemize}
\item[(i)] $K_{S^\nu}+B_{S^\nu}$ is 
$\mathbb R$-Cartier and $K_{S^{\nu}}+B_{S^\nu}=\nu^*(K_X+S+B)$. 
\item[(ii)] If $B$ is a $\mathbb Q$-divisor, then 
so is $B_{S^{\nu}}$. 
\item[(iii)] $B_{S^\nu}$ is effective if $B$ is 
effective in a neighborhood of $S$. 
\item[(iv)] $(S^\nu, B_{S^\nu})$ is log canonical 
if $(X, S+B)$ is log canonical 
in a neighborhood of $S$. 
\item[(v)] Let $D$ be an $\mathbb R$-Cartier $\mathbb R$-divisor 
on $X$ such that 
$S$ and $D$ have no common irreducible components. Then 
we have 
$$
(B+D)_{S^\nu}=B_{S^\nu}+\nu^*D. 
$$
We sometimes write $D|_{S^{\nu}}=\nu^*D$ for simplicity. 
\end{itemize}
The properties except (iii) are obvious by the definition. 
We give a proof of (iii) for the reader's convenience. 
\end{say}

\begin{proof}[Proof of {\em{(iii)}}] 
By shrinking $X$, we can assume that 
$X$ is quasi-projective and 
$B$ is effective. 
By taking hyperplane cuts, we can also assume 
that 
$X$ is a surface. 
Run the log minimal model program 
over $X$ with respect to 
$K_Y+S_Y$. Let $C$ be a curve 
on $Y$ such that 
$(K_Y+S_Y)\cdot C<0$ and $f(C)$ is a 
point. 
Then $K_Y\cdot C<0$ because $S_Y$ is the strict transform 
of $S$. Therefore, 
each step of the log minimal model 
program over 
$X$ with 
respect to $K_Y+S_Y$ is a contraction of a $(-1)$-curve $E$ 
with $(K_Y+S_Y)\cdot E<0$. 
So, by replacing $(Y, S_Y)$ with 
the output of 
the above log minimal model program, 
we can assume that $Y$ is smooth, 
$(Y, S_Y)$ is plt, 
and $K_Y+S_Y$ is $f$-nef. We 
note that 
$S_Y$ is a smooth 
curve since $(Y, S_Y)$ is plt (cf.~\cite[Proposition 5.51]{km}). 
By the negativity lemma (see Lemma \ref{lem211}) and 
the assumption 
that $B$ is effective, 
$B_Y$ is effective. 
We note the following equality 
$$
-B_Y=K_Y+S_Y-f^*(K_X+S+B). 
$$ 
By adjunction, we obtain 
$$
(K_Y+S_Y+B_Y)|_{S_Y}=K_{S_Y}+B_Y|_{S_Y}. 
$$ 
It is obvious that $B_Y|_{S_Y}$ is effective. 
This implies that $B_{S^{\nu}}=B_Y|_{S_Y}$ is effective. 
\end{proof}
When $X$ is singular, $B_{S^\nu}$ is not necessarily zero 
even if $B=0$. 

\section{Rationality theorem}\label{sec14}
In this section, we prove the following rationality theorem, 
though it is a special case of \cite[Theorem 5.9]{ambro} (see 
also \cite[Theorem 3.68]{book}). 
In the traditional X-method, the rationality theorem for klt pairs 
is proved by the Kawamata--Viehweg vanishing theorem, 
Hironaka's resolution theorem, and Shokurov's non-vanishing theorem 
(see, for example, \cite[\S 3.4]{km}). 
Our proof of the rationality theorem given below only uses the vanishing 
theorem:~Theorem \ref{thm71}. 
We do not need the non-vanishing theorem (cf.~Theorem \ref{thm111}) 
nor Hironaka's resolution theorem in this section. 

\begin{thm}[Rationality theorem]\label{thm131} 
Let $X$ be a normal variety and let $B$ be an effective $\mathbb Q$-divisor on 
$X$ such that $K_X+B$ is $\mathbb Q$-Cartier. 
Let $\pi:X\to S$ be a projective morphism and let $H$ be a 
$\pi$-ample Cartier divisor on $X$. Assume that 
$K_X+B$ is not 
$\pi$-nef and that $r$ is a positive number such that 
\begin{itemize}
\item[$(1)$] $H+r(K_X+B)$ is $\pi$-nef 
but not $\pi$-ample, and 
\item[$(2)$] $(H+r(K_X+B))|_{\Nlc(X, B)}$ 
is $\pi|_{\Nlc (X, B)}$-ample. 
\end{itemize}
Then $r$ is a rational number, and 
in reduced form, $r$ has denominator at most $a(\dim X+1)$, where $a(K_X+B)$ 
is a Cartier divisor on $X$. 
\end{thm} 

Before the proof of Theorem \ref{thm131}, we recall the following lemmas. 
 
\begin{lem}[{cf.~\cite[Lemma 3.19]{km}}]\label{lem132}  
Let $P(x, y)$ be a non-trivial 
polynomial of degree $\leq n$ and 
assume that 
$P$ vanishes for all sufficiently large integral solutions of 
$0<ay-rx<\varepsilon$ for some fixed positive integer $a$ and 
positive $\varepsilon$ for some $r\in \mathbb R$. 
Then $r$ is rational, and 
in reduced form, $r$ has denominator $\leq 
a(n+1)/\varepsilon$. 
\end{lem}
\begin{proof}
We assume that $r$ is irrational. 
Then an infinite number of integral points in the $(x, y)$-plane 
on each side of the line $ay-rx=0$ are closer than 
$\varepsilon /(n+2)$ to that line. 
So there is a large integral solution $(x', y')$ 
with 
$0<ay'-rx'<\varepsilon /(n+2)$.  
In this case, we see that $$(2x', 2y'), \cdots, ((n+1)x', (n+1)y')$$ 
are also solutions by hypothesis. 
So $(y'x-x'y)$ divides $P$, since 
$P$ and $(y'x-x'y)$ have $(n+1)$ common zeroes. 
We choose a smaller $\varepsilon$ and repeat 
the argument. 
We do this $n+1$ times 
to get a contradiction. 

Now we assume that $r=u/v$ in lowest terms. 
For given $j$, let 
$(x', y')$ be a solution 
of $ay-rx=aj/v$. 
Note that an integral 
solution exists for every $j$. 
Then we have 
$a(y'+ku)-r(x'+akv)=aj/v$ for all $k$. 
So, as above, if $aj/v<\varepsilon$, 
$(ay-rx)-(aj/v)$ must divide $P$. So 
we can have at most $n$ such values of $j$. Thus $a(n+1)/v\geq \varepsilon$.  
\end{proof}

\begin{lem}\label{lem133}  
Let $C$ be a projective variety and let $D_1$ and $D_2$ be Cartier 
divisors on $X$. 
Consider the Hilbert polynomial 
$$
P(u_1, u_2)=\chi (C, \mathcal O_C(u_1 D_1+u_2D_2)). 
$$
If $D_1$ is ample, then 
$P(u_1, u_2)$ is a non-trivial 
polynomial of total degree $\leq \dim C$. 
It is because $P(u_1, 0)=h^0(C, \mathcal O_C(u_1D_1))\not\equiv 0$ if $u_1$ 
is sufficiently large. 
\end{lem}

\begin{proof}[Proof of {\em{Theorem \ref{thm131}}}]
Let $m$ be a positive integer such that 
$H'=mH$ is $\pi$-very ample. 
If $H'+r'(K_X+B)$ is $\pi$-nef but not 
$\pi$-ample, and $(H'+r'(K_X+B))|_{\Nlc (X, B)}$ is 
$\pi|_{\Nlc (X, B)}$-ample, 
then we have 
$$
H+r(K_X+B)=\frac{1}{m}(H'+r'(K_X+B)). 
$$ 
This gives $r=\frac{1}{m}r'$. 
Thus, $r$ is rational if and only if 
$r'$ is rational. 
Assume furthermore that $r'$ has 
denominator $v$. 
Then $r$ has denominator dividing $mv$. 
Since $m$ can be arbitrary sufficiently 
large integer, this implies that $r$ has 
denominator dividing $v$. 
Therefore, by replacing $H$ with 
$mH$, we can assume that $H$ is 
very ample over $S$. 

For each $(p, q)\in \mathbb Z^2$, 
let $L(p, q)$ denote the relative base 
locus of the linear system $M(p, q)$ on $X$ (with 
the reduced scheme structure), 
that is, 
$$
L(p, q)=\Supp (\Coker (\pi^*\pi_*
\mathcal O_X(M(p, q))\to \mathcal O_X(M(p, q)))), 
$$
where $M(p, q)=pH+qa(K_X+B)$. 
By 
the definition, $L(p, q)=X$ if and only if $\pi_*\mathcal O_X(M(p, q))= 0$. 

\setcounter{cla}{0}
\begin{cla}[{cf.~\cite[Claim 3.20]{km}}]\label{rat-c1}  
Let $\varepsilon$ be a positive number. 
For $(p, q)$ sufficiently 
large and $0<aq-rp<\varepsilon$, $L(p, q)$ is the same 
subset of $X$. 
We call this subset $L_0$. 
Let $I\subset \mathbb Z^2$ be the set of $(p, q)$ for 
which $0<aq-rp<1$ and $L(p, q)=L_0$. 
We note that $I$ contains all sufficiently 
large $(p, q)$ with $0<aq-rp<1$. 
\end{cla}

\begin{proof}
We fix $(p_0, q_0)\in \mathbb Z^2$ such that 
$p_0>0$ and $0<aq_0-rp_0<1$. 
Since $H$ is $\pi$-very ample, 
there exists a positive integer 
$m_0$ such that 
$\mathcal O_X(mH+ja(K_X+B))$ is $\pi$-generated 
for every $m>m_0$ and 
every $0\leq j\leq q_0-1$. 
Let $M$ be the round-up of 
$$
\Big(m_0+\frac{1}{r}\Big)\Big\slash\Big(\frac{a}{r}-\frac{p_0}{q_0}\Big). 
$$
If $(p', q')\in \mathbb Z^2$ such that 
$0<aq'-rp'<1$ and $q'\geq M+q_0-1$, 
then we can write 
\begin{align*}
p'H+q'a(K_X+B)=k(p_0H+q_0a(K_X+B))+(lH+ja(K_X+B))
\end{align*} 
for some $k\geq 0$, $0\leq j\leq q_0-1$ 
with $l> m_0$. 
It is because we can uniquely write $q'=kq_0+j$ with $0\leq j\leq q_0-1$. Thus, 
we have $kq_0\geq M$. So, 
we obtain 
\begin{align*}
l=p'-kp_0>\frac{a}{r}q'-\frac{1}{r}-(kq_0)\frac{p_0}{q_0}
\geq \Big(\frac{a}{r}-\frac{p_0}{q_0}\Big)M-\frac{1}{r}\geq m_0. 
\end{align*}
Therefore, 
$L(p', q')\subset L(p_0, q_0)$. 
By the noetherian induction, 
we obtain the desired closed subset $L_0\subset X$ and 
$I\subset \mathbb Z^2$. 
\end{proof}

\begin{cla}\label{rat-c2} 
We have $L_0\cap \Nlc (X, B)=\emptyset$. 
\end{cla}

\begin{proof}[Proof of {\em{Claim \ref{rat-c2}}}] 
We take $(\alpha, \beta)\in \mathbb Q^2$ such that 
$\alpha>0$, $\beta>0$, 
and $\beta a/\alpha>r$ is sufficiently close to 
$r$. 
Then $(\alpha H+\beta a(K_X+B))|_{\Nlc(X, B)}$ is 
$\pi|_{\Nlc(X, B)}$-ample because 
$(H+r(K_X+B))|_{\Nlc(X, B)}$ is 
$\pi|_{\Nlc(X, B)}$-ample. 
If $0<aq-rp<1$ and $(p, q)\in \mathbb Z^2$ is 
sufficiently large, then 
$$M(p, q)=mM(\alpha, \beta)+(M(p, q)-mM(\alpha, \beta))$$ such 
that $M(p, q)-mM(\alpha, \beta)$ is $\pi$-very ample and 
that $$m(\alpha H+\beta a (K_X+B))|_{\Nlc(X, B)}$$ is 
also $\pi|_{\Nlc(X, B)}$-very ample. 
It can be checked by the same argument as in the proof of 
Claim \ref{rat-c1}. 
Therefore, $\mathcal O_{\Nlc(X, B)}(M(p, q))$ is $\pi$-very 
ample. 
Since $$\pi_*\mathcal O_X(M(p, q))\to 
\pi_*\mathcal O_{\Nlc(X, B)}(M(p, q))
$$ is surjective 
by the vanishing theorem:~Theorem \ref{thm71}, 
we obtain $L(p, q)\cap \Nlc(X, B)=\emptyset$. 
We note that 
$$M(p, q)-(K_X+B)=pH+(qa-1)(K_X+B)$$ is $\pi$-ample because 
$(p, q)$ is sufficiently large and $aq-rp<1$.  
By Claim \ref{rat-c1}, we have 
$L_0\cap \Nlc(X, B)=\emptyset$. 
\end{proof}

\begin{cla}\label{rat-c3} 
We assume that $r$ is not rational or that 
$r$ is rational and has denominator $>a(n+1)$ in reduced form, 
where $n=\dim X$. 
Then, for $(p, q)$ sufficiently large and $0<aq-rp<1$, 
$\mathcal O_X(M(p, q))$ is $\pi$-generated 
at the generic point of every lc center 
of $(X, B)$. 
\end{cla}

\begin{proof}[Proof of {\em{Claim \ref{rat-c3}}}]
We note that $$M(p, q)-(K_X+B)=pH+(qa-1)(K_X+B).$$  If 
$aq-rp<1$ and $(p, q)$ is sufficiently large, 
then $M(p, q)-(K_X+B)$ is 
$\pi$-ample. 
Let $C$ be an lc center of $(X, B)$. 
We note that we can assume $C\cap 
\Nlc (X, B)=\emptyset$ by Claim \ref{rat-c2}. 
Then $P_{C_\eta}(p, q)=\chi (C_{\eta}, \mathcal 
O_{C_\eta}(M(p, q)))$ is a non-zero polynomial of degree 
at most $\dim C_{\eta}\leq \dim X$ by Lemma \ref{lem133}. 
Note that $C_\eta$ is the generic fiber of $C\to \pi(C)$. 
By Lemma \ref{lem132}, there exists $(p, q)$ such that 
$P_{C_\eta}(p, q)\ne 0$, 
$(p, q)$ sufficiently large, and $0<aq-rp<1$. 
By the $\pi$-ampleness of $M(p, q)-(K_X+B)$, $$P_{C_\eta}
(p, q)=
\chi (C_\eta, \mathcal O_{C_\eta}(M(p, q)))
=h^0(C_\eta, \mathcal O_{C_\eta}(M(p, q)))$$ and 
$$
\pi_*\mathcal O_X(M(p, q))\to \pi_*\mathcal O_C(M(p, q))
$$ 
is surjective by Theorem \ref{thm71}. 
We note that $C\cap \Nlc (X, B)=\emptyset$. 
Therefore, $\mathcal O_X(M(p, q))$ 
is $\pi$-generated 
at the generic point of $C$. 
By combining this 
with Claim \ref{rat-c1}, $\mathcal O_X(M(p, q))$ is 
$\pi$-generated at the generic point of every 
lc center of $(X, B)$ if 
$(p, q)$ is sufficiently large with 
$0<aq-rp<1$. 
So, we obtain Claim \ref{rat-c3}. 
\end{proof}

Note 
that $\mathcal O_X(M(p, q))$ is not $\pi$-generated 
for $(p, q)\in I$ 
because $M(p, q)$ is not $\pi$-nef. 
Therefore, $L_0\ne \emptyset$. 
We shrink $S$ to an affine open subset intersecting $\pi(L_0)$. 
Let $D_1, \cdots, D_{n+1}$ 
be general members of $\pi_*\mathcal O_X(M(p_0, q_0))
=H^0(X, \mathcal O_X(M(p_0, q_0)))$ with 
$(p_0, q_0)\in I$. 
We can check that $K_X+B+\sum _{i=1}^{n+1}D_i$ 
is not lc at the generic point of every irreducible component of $L_0$ 
by 
Lemma \ref{lem134}. 
On the other hand, $K_X+B+\sum _{i=1}^{n+1}D_i$ 
is lc outside $L_0\cup \Nlc(X, B)$ since 
$D_i$ is a general member of $|M(p_0, q_0)|$ for every $i$. 
Let $0<c<1$ be the log canonical threshold of $(X, B)$ 
with respect to $D=\sum _{i=1}^{n+1}D_i$ outside 
$\Nlc (X, B)$. 
Note that $c>0$ by Claim \ref{rat-c3}. 
Thus, the pair $(X, B+cD)$ has some 
lc centers contained in $L_0$. Let $C$ be an lc center 
contained in $L_0$. 
We note that $\mathcal J_{NLC}(X, B+cD)=\mathcal 
J_{NLC}(X, B)$ by Proposition \ref{prop63} 
and that $C\cap \Nlc (X, B+cD)=C\cap \Nlc (X, B)=\emptyset$. 
We consider 
$$K_X+B+cD=c(n+1)p_0H+(1+c(n+1)q_0a)(K_X+B). $$ 
Thus we have 
\begin{align*}
&pH+qa(K_X+B)-(K_X+B+cD)\\ &
=(p-c(n+1)p_0)H+(qa-(1+c(n+1)q_0a))(K_X+B).
\end{align*} 
If $p$ and $q$ are large enough and $0<aq-rp\leq aq_0-rp_0$, 
then $$pH+qa(K_X+B)-(K_X+B+cD)$$ is $\pi$-ample. 
It is because 
\begin{align*}
&(p-c(n+1)p_0)H+(qa-(1+c(n+1)q_0a))(K_X+B)\\
&=(p-(1+c(n+1))p_0)H+(qa-(1+c(n+1))q_0a)(K_X+B)
\\ &\ +p_0H+(q_0a-1)(K_X+B).  
\end{align*}

Suppose that $r$ is not rational. 
There must be arbitrarily large $(p, q)$ such that 
$0<aq-rp<\varepsilon =aq_0-rp_0$ and 
$\chi(C_\eta, \mathcal O_{C_\eta}(M(p, q)))\ne 0$ 
by Lemma \ref{lem132} because 
$P_{C_\eta}(p, q)=\chi (C_\eta, 
\mathcal O_{C_\eta}(M(p, q)))$ is 
a non-trivial 
polynomial of degree at most $\dim C_\eta$ by Lemma \ref{lem133}. 
Since $M(p, q)-(K_X+B+cD)$ is $\pi$-ample by 
$0<aq-rp<aq_0-rp_0$, 
we have $h^0(C_\eta, \mathcal O_{C_\eta}(M(p, q)))=
\chi (C_\eta, \mathcal O_{C_\eta}(M(p, q)))\ne 0$ by the 
vanishing theorem:~Theorem \ref{thm71}. 
By the vanishing theorem:~Theorem \ref{thm71}, 
$$
\pi_*\mathcal O_X(M(p, q))\to \pi_*\mathcal O_C(M(p, q))
$$ 
is surjective because $M(p, q)-(K_X+B+cD)$ is $\pi$-ample.
We note that $C\cap \Nlc (X, B+cD)=\emptyset$. 
Thus $C$ is not contained in $L(p, q)$. Therefore, 
$L(p, q)$ is a proper subset of $L(p_0, q_0)=L_0$, 
giving the desired contradiction. 
So now we know that $r$ is rational. 

We next suppose that the assertion of the theorem concerning 
the denominator of $r$ is false. 
We choose $(p_0, q_0)\in I$ such that 
$aq_0-rp_0$ is the maximum, 
say it is equal to $d/v$. 
If $0<aq-rp\leq d/v$ and 
$(p, q)$ is sufficiently large, then $\chi (C_\eta, 
\mathcal O_{C_\eta}(M(p, q)))=h^0(C_\eta, 
\mathcal O_{C_\eta}(M(p, q)))$ since 
$M(p, q)-(K_X+B+cD)$ is $\pi$-ample. 
There exists sufficiently 
large $(p, q)$ in the strip 
$0<aq-rp<1$ with $\varepsilon =1$ for which 
$h^0(C_\eta, \mathcal O_{C_\eta}(M(p, q)))
=\chi (C_\eta, \mathcal O_{C_\eta}(M(p, q)))\ne 0$ by 
Lemma \ref{lem132} since 
$\chi (C_\eta, \mathcal O_{C_\eta}(M(p, q)))$ is a 
non-trivial 
polynomial of degree at most $\dim C_\eta$ by Lemma \ref{lem133}. 
Note that $aq-rp\leq d/v=aq_0-rp_0$ holds automatically for 
$(p, q)\in I$. 
Since 
$$\pi_*\mathcal O_X(M(p, q))
\to \pi_*\mathcal O_C(M(p,q))$$ 
is surjective by the $\pi$-ampleness of $M(p,q)-(K_X+B+cD)$, 
we 
obtain the desired contradiction by the same reason as above. 
So, we finish the proof of the rationality theorem. 
\end{proof}

We close this section with an important remark, which is 
indispensable for the proof of the cone theorem:~Theorem \ref{thm144}.

\begin{rem}\label{rem13-saigo} 
In Theorem \ref{thm131}, it is sufficient to 
assume that $B$ is an effective $\mathbb R$-divisor on $X$ such 
that $K_X+B$ is $\mathbb R$-linearly equivalent to 
a $\mathbb Q$-Cartier $\mathbb Q$-divisor $\omega$ on $X$ with the 
condition 
that $a\omega$ is Cartier. 
All we have to do is to replace $a(K_X+B)$ with 
$a\omega$ in the proof of Theorem \ref{thm131}. 
We need this generalization in the proof of the cone theorem:~Theorem \ref{thm144}. 
\end{rem}

\section{Cone theorem}\label{sec15}

The main theorem of this section 
is the cone theorem.  
Before we state the main theorem, let us 
fix the notation. 

\begin{defn}\label{def141}
Let $X$ be a normal variety and let $B$ be an effective 
$\mathbb R$-divisor on $X$ such 
that $K_X+B$ is $\mathbb R$-Cartier. 
Let $\pi:X\to S$ be a projective morphism. 
We put 
$$
\overline {NE}(X/S)_{\Nlc(X, B)}=
\xIm (\overline {NE}(\Nlc(X, B)/S)
\to \overline {NE}(X/S)). 
$$
For an $\mathbb R$-Cartier $\mathbb R$-divisor $D$, we define 
$$
D_{\geq 0}=\{z\in N_1(X/S) \ | \ D\cdot z\geq 0\}. 
$$ 
Similarly, we can define $D_{>0}$, 
$D_{\leq 0}$, and 
$D_{<0}$. 
We also define 
$$
D^{\perp}=\{ z\in N_1(X/S) \ | \ D\cdot z=0\}. 
$$ 
We use the following notation 
$$
\overline {NE}(X/S)_{D\geq 0} =\overline {NE}(X/S)\cap 
D_{\geq 0}, 
$$ 
and similarly for $>0$, $\leq 0$, 
and $<0$. 
\end{defn}

\begin{defn}\label{def142}
An {\em{extremal face}} of $\overline {NE}(X/S)$ is a non-zero 
subcone $F\subset \overline {NE}(X/S)$ such that 
$z, z'\in F$ and $z+z'\in F$ imply that 
$z, z'\in F$. Equivalently, 
$F=\overline {NE}(X/S)\cap H^{\perp}$ for some 
$\pi$-nef $\mathbb R$-divisor $H$, which 
is called a {\em{supporting function}} of $F$. 
An {\em{extremal ray}} 
is a one-dimensional 
extremal face. 
\begin{itemize}
\item[(1)] An extremal face $F$ is called 
{\em{$(K_X+B)$-negative}} 
if $$F\cap \overline {NE}(X/S)_{K_X+B\geq 0} 
=\{ 0\}. $$ 
\item[(2)] An extremal face $F$ is called {\em{rational}} 
if we can choose 
a $\pi$-nef $\mathbb Q$-divisor $H$ as a support 
function of $F$. 
\item[(3)] An extremal face $F$ is called {\em{relatively ample at 
$\Nlc (X, B)$}} if $$F\cap \overline {NE}(X/S)_{\Nlc (X, B)}=\{0\}. 
$$
Equivalently, $H|_{\Nlc (X, B)}$ is $\pi|_{\Nlc (X, B)}$-ample for 
every supporting function $H$ of $F$. 
\item[(4)] An extremal face $F$ is called {\em{contractible 
at $\Nlc (X, B)$}} if it has a rational supporting function 
$H$ such that $H|_{\Nlc (X, B)}$ is $\pi|_{\Nlc (X, B)}$-semi-ample. 
\end{itemize}
\end{defn}

\begin{rem} 
If $X$ is complete but non-projective, 
then it sometimes happens that $\overline {NE}(X)=N_1(X)$ 
even when $X$ is smooth (cf.~\cite{fp}). 
Therefore, the projectivity is crucial 
for the log minimal model program.
\end{rem}

The following theorem is a direct consequence 
of Theorem \ref{thm121}. 

\begin{thm}[Contraction theorem]\label{thm143}
Let $X$ be a normal variety, let $B$ be an effective $\mathbb R$-divisor 
on $X$ such that 
$K_X+B$ is $\mathbb R$-Cartier, and 
let $\pi:X\to S$ be a projective 
morphism. 
Let $H$ be a $\pi$-nef 
Cartier divisor such that 
$F=H^{\perp}\cap \overline {NE}(X/S)$ is $(K_X+B)$-negative 
and contractible at $\Nlc (X, B)$. 
Then there exists a projective morphism 
$\varphi_F:X\to Y$ over $S$ with the following properties. 
\begin{itemize}
\item[$(1)$] 
Let $C$ be an integral curve on $X$ such that 
$\pi(C)$ is a point. 
Then $\varphi_F(C)$ is a point if and 
only if 
$[C]\in F$. 
\item[$(2)$] $\mathcal O_Y\simeq (\varphi_F)_*\mathcal O_X$.  
\item[$(3)$] Let $L$ be a line bundle on $X$ such 
that $L\cdot C=0$ for 
every curve $C$ with $[C]\in F$. 
Assume that $L^{\otimes m}|_{\Nlc(X, B)}$ 
is $\varphi_{F}|_{\Nlc (X, B)}$-generated for $m \gg 0$. 
Then there is a line bundle $L_Y$ on $Y$ such 
that $L\simeq \varphi^*_FL_Y$. 
\end{itemize} 
\end{thm}

\begin{proof}
By the assumption, $qH-(K_X+B)$ is 
$\pi$-ample for some positive integer $q$ and 
we can assume that $H|_{\Nlc (X, B)}$ is $\pi|_{\Nlc (X, B)}$-semi-ample. 
By Theorem \ref{thm121}, 
$\mathcal O_X(mH)$ is $\pi$-generated 
for some $m >0$. 
We take the Stein factorization of the associated 
morphism. Then, we have the contraction morphism 
$\varphi_F: X\to Y$ with the properties (1) and (2). 

We consider $\varphi_F:X\to Y$ and $\overline {NE}(X/Y)$. 
Then $\overline {NE}(X/Y)=F$, $L$ is 
numerically trivial over $Y$, 
and $-(K_X+B)$ 
is $\varphi_F$-ample. Applying the base point free 
theorem (cf.~Theorem 
\ref{thm121}) over $Y$, 
both $L^{\otimes m}$ and $L^{\otimes (m+1)}$ 
are pull-backs of line bundles on $Y$. Their difference 
gives a line bundle $L_Y$ such that 
$L\simeq \varphi^*_FL_Y$. 
\end{proof}

\begin{ex}
Let $S$ be a cone over a smooth cubic curve and let $\pi:X\to S$ 
be the blow-up at the vertex of $S$. 
Then $K_X+E=\pi^*K_S$, where $E$ is the $\pi$-exceptional 
curve. We put $B=2E$ and consider the pair $(X, B)$. 
In this case, $\varphi_F=\pi:X\to Y=S$ with 
$F=0^{\perp}\cap \overline {NE}(X/S)=\overline {NE}(X/S)=\mathbb R_{\geq 0}[E]$ 
is an example of 
contraction morphisms in Theorem \ref{thm143}. 
\end{ex}

The time is ripe to state one of the main theorems in this paper. 

\begin{thm}[Cone theorem]\label{thm144} 
Let $X$ be a normal variety, let $B$ be an effective $\mathbb R$-divisor 
on $X$ such that $K_X+B$ is $\mathbb R$-Cartier, and 
let $\pi:X\to S$ be a projective morphism. 
Then we have the following 
properties. 
\begin{itemize}
\item[$(1)$] $\overline {NE}(X/S)=\overline {NE}(X/S)_{K_X+B\geq 0} 
+\overline {NE}(X/S)_{\Nlc (X, B)}+\sum R_j$, 
where $R_j$'s are the $(K_X+B)$-negative 
extremal rays of $\overline {NE}(X/S)$ that are 
rational and relatively ample at $\Nlc (X, B)$. 
In particular, each $R_j$ is spanned by 
an integral curve $C_j$ on $X$ such that 
$\pi(C_j)$ is a point.  
\item[$(2)$] Let $H$ be a $\pi$-ample $\mathbb R$-divisor 
on $X$. 
Then there are only finitely many $R_j$'s included in 
$(K_X+B+H)_{<0}$. In particular, 
the $R_j$'s are discrete in the half-space 
$(K_X+B)_{<0}$. 
\item[$(3)$] Let $F$ be a $(K_X+B)$-negative extremal 
face of $\overline {NE}(X/S)$ that is 
relatively ample at $\Nlc (X, B)$. 
Then $F$ is a rational face. 
In particular, $F$ is contractible at $\Nlc (X, B)$. 
\end{itemize}
\end{thm}

\begin{proof}
First, we assume that $K_X+B$ is $\mathbb R$-linearly 
equivalent to a $\mathbb Q$-Cartier $\mathbb Q$-divisor on $X$ 
(see Remark \ref{rem13-saigo}).  
We can assume that $\dim _{\mathbb R}
N_1(X/S)\geq 2$ and $K_X+B$ is not $\pi$-nef. 
Otherwise, the theorem is obvious. 

\setcounter{step}{0}
\begin{step}\label{cone-st1} 
We have 
$$
\overline {NE}(X/S)=\overline 
{\overline {NE}(X/S)_{K_X+B\geq 0} 
+\overline {NE}(X/S)_{\Nlc (X, B)}+\sum _F F},  
$$ 
where $F$'s vary among all rational 
proper $(K_X+B)$-negative faces that are relatively 
ample at $\Nlc (X, B)$ and 
$\raise0.5ex\hbox{\textbf{-----}}$
denotes the closure 
with respect to 
the real topology. 
\end{step} 
\begin{proof}We put 
$$
\mathfrak B=\overline {\overline {NE}(X/S)_{K_X+B\geq 0} 
+\overline {NE}(X/S)_{\Nlc (X, B)}+\sum _F F}. 
$$ 
It is clear that $\overline {NE}(X/S)\supset \mathfrak B$. 
We note that each $F$ is spanned by curves 
on $X$ mapped to points on $S$ by 
Theorem \ref{thm143} (1). Supposing
$\overline {NE}(X/S)\ne \mathfrak B$, we shall 
derive a contradiction. 
There is a separating function $M$ which 
is Cartier and is not a multiple of $K_X+B$ in $N^1(X/S)$ such 
that $M>0$ on $\mathfrak B\setminus \{0\}$ and 
$M\cdot z_0<0$ for some $z_0\in \overline {NE}(X/S)$. 
Let $C$ be the dual cone of $\overline {NE}(X/S)_{K_X+B\geq 0}$, 
that is, 
$$
C=\{D\in N^1(X/S)\ | \ D\cdot z\geq 0 \ \text{for}\ z\in 
\overline {NE}(X/S)_{K_X+B \geq 0}\}. 
$$
Then $C$ is generated by $\pi$-nef divisors and $K_X+B$. 
Since $M>0$ on $\overline {NE}(X/S)_{K_X+B\geq 0}\setminus \{0\}$, 
$M$ is in the interior of $C$, and hence there 
exists a $\pi$-ample $\mathbb Q$-Cartier 
$\mathbb Q$-divisor $A$ such that 
$M-A=L'+p(K_X+B)$ in $N^1(X/S)$, where 
$L'$ is a $\pi$-nef $\mathbb Q$-Cartier $\mathbb Q$-divisor 
on $X$ and $p$ is a non-negative rational number. 
Therefore, $M$ is expressed in the form $M=H+p(K_X+B)$ in 
$N^1(X/S)$, where $H=A+L'$ is a $\pi$-ample 
$\mathbb Q$-Cartier $\mathbb Q$-divisor. 
The rationality theorem (see Theorem \ref{thm131}) 
implies that there exists a positive 
rational number $r<p$ such that $L=H+r(K_X+B)$ is 
$\pi$-nef but not 
$\pi$-ample, and $L|_{\Nlc (X, B)}$ is 
$\pi|_{\Nlc (X, B)}$-ample. 
Note that $L\ne 0$ in $N^1(X/S)$, since 
$M$ is not a multiple of $K_X+B$. 
Thus the extremal face $F_L$ associated to 
the supporting function $L$ is contained 
in $\mathfrak B$, which implies $M>0$ on $F_L$. 
Therefore, $p<r$. It is a contradiction. 
This completes the proof of our first claim. 
\end{proof}

\begin{step}
In the equality of Step \ref{cone-st1}, 
we can assume that every extremal face $F$ is one-dimensional. 
\end{step}
\begin{proof}
Let $F$ be a rational proper $(K_X+B)$-negative extremal face 
that is relatively ample at $\Nlc (X, B)$, and assume that $\dim F\geq 
2$. Let $\varphi_F:X\to W$ be the associated 
contraction. 
Note that $-(K_X+B)$ is $\varphi_F$-ample. 
By Step \ref{cone-st1}, 
we obtain 
$$
F=\overline {NE}(X/W)=\overline {\sum _G G}, 
$$ 
where the $G$'s are the rational proper 
$(K_X+B)$-negative extremal faces of $\overline {NE}(X/W)$. 
We note that $\overline {NE}(X/W)_{\Nlc (X, B)}=0$ because 
$\varphi_F$ embeds $\Nlc (X, B)$ into $W$. 
The $G$'s are also 
$(K_X+B)$-negative 
extremal faces of $\overline {NE}(X/S)$ that are 
ample at $\Nlc (X, B)$, 
and $\dim G<\dim F$. By induction, we obtain 
\begin{align*}\label{siki1} 
\overline {NE}(X/S)=\overline {\overline {NE}(X/S)_{K_X+B\geq 0} 
+\overline {NE}(X/S)_{\Nlc (X, B)} +\sum R_j},\tag{$\clubsuit$} 
\end{align*}
where the $R_j$'s are $(K_X+B)$-negative rational 
extremal rays. 
Note that each $R_j$ does not 
intersect 
$\overline {NE}(X/S)_{\Nlc (X, B)}$.  
\end{proof}

\begin{step}\label{cone-st3} 
The contraction theorem (cf.~Theorem \ref{thm143}) 
guarantees that 
for each extremal ray $R_j$ there exists a reduced irreducible curve $C_j$ on $X$ such that 
$[C_j]\in R_j$. 
Let $\psi_j:X\to W_j$ be the contraction 
morphism of $R_j$, and 
let $A$ be a $\pi$-ample 
Cartier divisor. 
We set 
$$
r_j=-\frac{A\cdot  C_j}{(K_X+B)\cdot C_j}. 
$$
Then $A+r_j(K_X+B)$ is $\psi_j$-nef but not $\psi_j$-ample, 
and 
$(A+r_j (K_X+B))|_{\Nlc (X, B)}$ 
is $\psi_j|_{\Nlc (X, B)}$-ample. 
By the rationality theorem 
$($see Theorem \ref{thm131}$)$, expressing 
$r_j =u_j /v_j$ with 
$u_j, v_j\in \mathbb Z_{>0}$ and 
$(u_j, v_j)=1$, we have 
the inequality $v_j\leq a(\dim X+1)$. 
\end{step}

\begin{step}
Now take $\pi$-ample Cartier divisors 
$H_1, H_2, \cdots, H_{\rho-1}$ such that 
$K_X+B$ and the $H_i$'s form a basis of $N^1(X/S)$, where 
$\rho =\dim _{\mathbb R}N^1(X/S)$. By 
Step \ref{cone-st3}, the intersection of  the extremal 
rays $R_j$ with the hyperplane 
$$\{z\in N_1(X/S)\ | \ a(K_X+B)\cdot z=-1\}$$ in 
$N_1(X/S)$ lie on the lattice 
$$
\Lambda=
\{z\in N_1(X/S)\ | \ a(K_X+B)\cdot z=-1, H_i \cdot z\in (a(a(\dim 
X+1))!)^{-1}\mathbb Z\}. 
$$ 
This implies that the extremal rays are discrete in the 
half space 
$$
\{z \in N_1(X/S) \ | \ (K_X+B)\cdot z<0\}. 
$$ 
Thus we can omit the closure sign 
$\raise0.5ex\hbox{\textbf{-----}}$ 
from 
the formula (\ref{siki1}) 
and this completes the proof of 
(1) when $K_X+B$ is $\mathbb R$-linearly equivalent to 
a $\mathbb Q$-Cartier $\mathbb Q$-divisor. 
\end{step} 
\begin{step}
Let $H$ be a $\pi$-ample $\mathbb R$-divisor on $X$. 
We choose $0<\varepsilon_i\ll 1$ for $1\leq i\leq \rho-1$ such 
that $H-\sum_{i=1}^{\rho-1}\varepsilon _i H_i$ is $\pi$-ample. 
Then the $R_j$'s included in $(K_X+B+H)_{<0}$ correspond to 
some elements of the above lattice $\Lambda$ for which 
$\sum_{i=1}^{\rho-1}\varepsilon _i H_i\cdot z<1/a$.
 Therefore, we obtain (2). 
\end{step}
\begin{step}\label{cone-st6} 
Let $F$ be a $(K_X+B)$-negative extremal face as in (3). 
The vector space $V=F^{\perp}\subset N^1(X/S)$ is 
defined over $\mathbb Q$ because 
$F$ is generated by some of the $R_j$'s. 
There exists a $\pi$-ample $\mathbb R$-divisor 
$H$ such that $F$ is contained in $(K_X+B+H)_{<0}$. 
Let $\langle F\rangle$ be the vector space 
spanned by $F$. 
We put 
$$
W_F=\overline {NE}(X/S)_{K_X+B+H\geq 0}+
\overline {NE}(X/S)_{\Nlc (X, B)}+\sum _{R_j\not\subset F}R_j. 
$$ 
Then $W_F$ is a closed cone, 
$\overline {NE}(X/S)=W_F+F$, 
and $W_F\cap \langle F\rangle=\{0\}$. The supporting 
functions of $F$ are the elements of $V$ that are 
positive on $W_F\setminus \{0\}$. 
This is a non-empty open set and thus it 
contains a rational element that, after scaling, 
gives a 
$\pi$-nef Cartier divisor $L$ such that 
$F=L^{\perp}\cap \overline {NE}(X/S)$. Therefore, 
$F$ is rational. 
So, we have (3). 
\end{step}
From now on, $K_X+B$ is $\mathbb R$-Cartier. 
\begin{step}
Let $H$ be a $\pi$-ample $\mathbb R$-divisor on $X$. 
We shall prove (2). 
We assume that there are infinitely many 
$R_j$'s in $(K_X+B+H)_{<0}$ and 
get a contradiction. 
There exists an affine open subset $U$ of 
$S$ such that 
$\overline {NE}(\pi^{-1}(U)/U)$ has infinitely many 
$(K_X+B+H)$-negative extremal rays. 
So, we shrink $S$ and can assume that 
$S$ is affine. 
We can write $H=E+H'$ such that $H'$ is $\pi$-ample, 
$\mathcal J_{NLC}(X, B+E)=\mathcal J_{NLC}(X, B)$, and 
$K_X+B+E$ is $\mathbb R$-linearly equivalent 
to a $\mathbb Q$-Cartier $\mathbb Q$-divisor. 
Since $K_X+B+H=K_X+B+E+H'$, we 
have 
$$
\overline {NE}(X/S)=\overline {NE}(X/S)_{K_X+B+H\geq 0}
+\overline {NE}(X/S)_{\Nlc (X, B)}+\sum _{\text{finite}}R_j. 
$$
It is a contradiction. 
Thus, we obtain (2). 
The statement (1) is a direct consequence of 
(2). Of course, 
(3) holds by Step \ref{cone-st6} 
once we obtain (1). 
\end{step}
So, we complete the proof of the cone theorem. 
\end{proof}

We close this section with the following elementary example. 

\begin{ex}
We consider $Y=\mathbb P^1\times \mathbb P^1$. 
Let $\pi_i:Y\to \mathbb P^1$ be the $i$-th projection for 
$i=1, 2$. 
Let $F_i$ be a fiber of $\pi_i$ for $i=1, 2$. 
We put $P=F_1\cap F_2$ and consider 
the blow-up $f:X\to Y$ at $P$. 
Let $E$ be the exceptional 
curve of $f$ and $C_i=f^{-1}_*F_i$ for 
$i=1, 2$. 
In this situation, we can check that 
$-K_X$ is ample, $\rho (X)=3$, and 
\begin{align*}
\overline {NE}(X)=\mathbb R_{\geq 0}[C_1] 
+\mathbb R_{\geq 0}[C_2]+\mathbb R_{\geq 0}[E]. 
\end{align*}
We put 
$$
B=\frac{3}{2}E+\frac{1}{2}C_1+C_2. 
$$ 
Then we have 
\begin{align*}
\overline {NE}(X)=\overline {NE}(X)_{K_X+B\geq 0}
+\overline {NE}(X)_{\Nlc (X, B)}+\mathbb R_{\geq 0}[C_2], 
\end{align*}
where $$\overline {NE}(X)_{\Nlc (X, B)}=\mathbb R_{\geq 0}[E], \quad
\overline {NE}(X)_{K_X+B\geq 0}=\mathbb R_{\geq 0}[C_1], $$
and $$C_2\cdot (K_X+B)<0. $$
\end{ex}

\section{Base point free theorem revisited}\label{sec-n17}

This section is a supplement to the base point free theorem:~Theorem \ref{thm121}. 
In the recent log minimal model program (cf.~\cite{sho-model}, 
\cite{bchm}, and so on), 
we frequently use $\mathbb R$-divisors. Therefore, 
the following theorem is useful.

\begin{thm}[Base point free theorem for $\mathbb R$-divisors]\label{thm171}
Let $(X, B)$ be a log canonical pair and 
let $\pi:X\to S$ be a projective 
morphism onto a variety $S$. 
Let $D$ be a nef $\mathbb R$-Cartier $\mathbb R$-divisor 
on $X$ such that $aD-(K_X+B)$ is ample for 
some real number $a>0$. Then 
$D$ is $\pi$-semi-ample. 
\end{thm}

\begin{proof}
We can assume that 
$a=1$ by replacing $D$ with $aD$. 
We put 
$$
F=\{ z\in \overline {NE}(X/S)\, | \, D\cdot z=0\}. 
$$ 
Then $F$ is a face of $\overline {NE}(X/S)$ and $(K_X+B)\cdot z<0$ for $z\in F$.
We claim that $F$ contains only finitely many $(K_X+B)$-negative 
extremal rays $R_1, \cdots, R_k$ of $\overline {NE}(X/S)$. 
If $F$ contains infinitely many $(K_X+B)$-negative extremal rays of $\overline {NE}(X/S)$, 
then it also holds after shrinking $S$ suitably. 
Therefore, we can assume that $S$ is affine. 
In this situation, $X$ is quasi-projective. 
We take a general small ample $\mathbb Q$-divisor $A$ on $X$ such that 
$D-(K_X+B+A)$ is ample and 
that $(X, B+A)$ is log canonical. 
Let $R$ be a $(K_X+B)$-negative extremal ray such that 
$R\subset F$. Then $R$ is a $(K_X+B+A)$-negative 
extremal ray since $D\cdot R=0$ and $D-(K_X+B+A)$ is ample. 
On the other hand, there are only finitely 
many $(K_X+B+A)$-negative extremal 
rays in $\overline {NE}(X/S)$ by Theorem \ref{thm144} (2). 
It is a contradiction. 
Therefore, $F$ is spanned by 
the extremal rays $R_1, \cdots, R_k$. 
We consider the finite dimensional real vector space $V=\underset {j}
\bigoplus \mathbb RD_j$, 
where $\sum _j D_j =\Supp D$ is the irreducible 
decomposition. Then 
$$
\mathcal R=\{ E\in V\, | \, E \ {\text{is $\mathbb R$-Cartier 
and $E\cdot z=0$ for every $z\in F$}} \} 
$$ 
is a rational affine subspace of $V$ and $D\in \mathcal R$. 
Thus, we can find positive real numbers $r_1, r_2, \cdots, r_m$ and 
nef $\mathbb Q$-Cartier $\mathbb Q$-divisors 
$E_1, E_2, \cdots, E_m$ such that 
$D=\sum _{i=1}^{m}r_iE_i$ and that 
$E_i-(K_X+B)$ is ample for every $i$ (cf.~Step \ref{cone-st6} in the 
proof of Theorem \ref{thm144}). 
By Theorem \ref{thm121}, 
$E_i$ is a semi-ample $\mathbb Q$-Cartier $\mathbb Q$-divisor 
for every $i$. 
Therefore, $D$ is semi-ample. 
\end{proof}

\section{Lengths of extremal rays}\label{sec16.5}
In this section, we discuss estimates of lengths of 
extremal rays. It is indispensable for 
the log minimal model program with scaling 
(see, for example, \cite{bchm}). 
Some results in this section have already been obtained in \cite{kollar2}, 
\cite{kollar3}, \cite{kawamata}, \cite{sho-model}, 
\cite{sho-7}, and \cite{birkar} with some extra assumptions. 
We note that the formulation of the main theorem of this section 
(cf.~Theorem \ref{thm-la}) is new.

Let us recall the following easy lemma. 

\begin{lem}[{cf.~\cite[Lemma 1]{sho-7}}]\label{lem145}
Let $(X, B)$ be a log canonical pair, 
where $B$ is 
an $\mathbb R$-divisor. Then there 
are positive real numbers 
$r_i$ and effective $\mathbb Q$-divisors $B_i$ for $1\leq i\leq l$ 
and a positive integer $m$ such that 
$\sum^{l}_{i=1}r_i=1$, $K_X+B=\sum ^{l}_{i=1}
r_i(K_X+B_i)$, $(X, B_i)$ is lc for every $i$, and 
$m(K_X+B_i)$ is Cartier for every $i$. 
\end{lem}
\begin{proof}
Let $\sum _kD_k$ be the irreducible decomposition 
of $\Supp B$. 
We consider the finite dimensional 
real vector space $V=\underset{k}{\bigoplus}\mathbb RD_k$. 
We put 
\begin{align*}
\mathcal Q=\left\{ D\in V \ |\  K_X+D \  {\text{is $\mathbb R$-Cartier}} 
\right\}. 
\end{align*}
Then it is easy to see that $\mathcal Q$ is an affine subspace of $V$ defined 
over $\mathbb Q$. 
We put 
\begin{align*}
\mathcal P=\left\{ D\in \mathcal Q \ |\  K_X+D \  {\text{is log canonical}} 
\right\}. 
\end{align*}
Thus by the definition of 
log canonicity, it is also easy to 
check that $\mathcal P$ is a closed convex rational 
polytope in $V$. 
We note that $\mathcal P$ is compact in the 
classical topology of $V$.  
By the assumption, $B\in \mathcal P$. 
Therefore, we can find the desired $\mathbb Q$-divisors 
$B_i\in \mathcal P$ and positive real numbers $r_i$. 
\end{proof}

The next result is essentially due to \cite{kawamata} and 
\cite[Proposition 1]{sho-7}. We will prove a more general result in Theorem \ref{thm-la} 
whose proof depends on Theorem \ref{prop146}. 

\begin{thm}\label{prop146}
Let $(X, B)$ be an lc pair and let 
$\pi:X\to S$ be a projective morphism onto a variety $S$. 
Let $R$ be a $(K_X+B)$-negative extremal ray 
of $\overline {NE}(X/S)$. 
Then we can find a rational curve $C$ on $X$ such that 
$[C]\in R$ and $$0<-(K_X+B)\cdot C\leq 2\dim X. $$
\end{thm}
\begin{proof}
By shrinking $S$, we can assume that 
$S$ is quasi-projective. 
By replacing $\pi:X\to S$ with 
the extremal contraction $\varphi_R:X\to Y$ over $S$, we can 
assume that the relative Picard number $\rho (X/S)=1$. 
In particular, $-(K_X+B)$ is $\pi$-ample. 
Let $K_X+B=\sum _{i=1}^{l}r_i(K_X+B_i)$ be as in 
Lemma \ref{lem145}. 
Without loss of generality, we can assume that 
$-(K_X+B_1)$ is $\pi$-ample and $-(K_X+B_i)=
-s_i(K_X+B_1)$ in $N^1(X/S)$ with 
$s_i\leq 1$ for every $i\geq 2$. 
Thus, it is sufficient to find a rational curve $C$ such that 
$\pi(C)$ is a point and that $-(K_X+B_1)\cdot C\leq 
2\dim X$. So, we can assume that $K_X+B$ is 
$\mathbb Q$-Cartier and lc. 
By Theorem \ref{thm91}, 
there is a birational morphism 
$f:(V, B_V)\to (X, B)$ such that $K_V+B_V=f^*(K_X+B)$, 
$V$ is $\mathbb Q$-factorial, and 
$(V, B_V)$ is dlt. 
By \cite[Theorem 1]{kawamata} and 
\cite[Theorem 10-2-1]{matsuki}, we can find 
a rational curve $C'$ on $V$ such that 
$-(K_V+B_V)\cdot C'\leq 2\dim V=2\dim X$ and 
that $C'$ spans 
a $(K_V+B_V)$-negative extremal ray. 
By the projection formula, the $f$-image of $C'$ is a desired 
rational curve. So, we finish the 
proof.    
\end{proof}

\begin{rem} 
It is conjectured that the 
estimate $\leq 2\dim X$ in Theorem \ref{prop146} 
should be replaced by $\leq \dim X+1$. 
When $X$ is smooth projective, it is true by Mori's 
famous result (cf.~\cite{mori-th}). See, for example, 
\cite[Theorem 1.13]{km}. 
When $X$ is a toric variety, it is also true by \cite{fuji-toric} 
and \cite{fuji-toric2}.  
\end{rem}

\begin{rem}
In the proof of Theorem \ref{prop146}, 
we need Kawamata's estimate on the length of an extremal rational 
curve (cf.~\cite[Theorem 1]{kawamata} and \cite[Theorem 10-2-1]{matsuki}). 
It depends on Mori's bend and break technique to create rational curves. 
So, we need the mod $p$ reduction technique there.  
\end{rem}

\begin{rem}
We give a remark on \cite{bchm}. We use the 
same notation as in \cite[3.8]{bchm}. 
In the proof of \cite[Corollary 3.8.2]{bchm}, we can assume that $K_X+\Delta$ 
is klt by \cite[Lemma 3.7.4]{bchm}. 
By perturbing the coefficients of $B$ slightly, we 
can further assume that 
$B$ is a $\mathbb Q$-divisor. By applying the usual 
cone theorem to the klt pair $(X, B)$, we obtain that 
there are only finitely many $(K_X+\Delta)$-negative 
extremal rays of $\overline {NE}(X/U)$. 
We note that \cite[Theorem 3.8.1]{bchm} is only used 
in the proof of \cite[Corollary 3.8.2]{bchm}. 
Therefore, we do not need the estimate of lengths of extremal rays in \cite{bchm}. 
In particular, we do not need mod $p$ reduction 
arguments for the proof of the main results in \cite{bchm}. 
\end{rem}

By the proof of Theorem \ref{prop146}, 
we have the following corollary. 

\begin{cor}\label{cor173}
Let $(X, B)$ be a log canonical pair 
and let $K_X+B=\sum _{i=1}^{l}r_i (K_X+B_i)$ and $m$ 
be as in {\em{Lemma \ref{lem145}}}. 
Let $\varphi:X\to Y$ be a projective surjective morphism 
with connected fibers such that 
the relative Picard number $\rho (X/Y)=1$. 
Then we can find a curve 
$C$ on $X$ such that 
$C$  spans $N_1(X/Y)$ and 
$$-(K_X+B_i)\cdot C=\frac{n_i}{m}$$ with 
$n_i\leq 2m\dim X$ for every $i$. 
Of course, we have 
$$
-(K_X+B)\cdot C=\sum _i \frac{r_in_i}{m}\leq 2\dim X. 
$$ 
If $-(K_X+B_i)$ is $\varphi$-ample for some $i$, then 
we can find a rational curve $C$ in the above statement. 
We note that $\varphi$ is not necessarily 
assumed to be a $(K_X+B)$-negative extremal contraction. 
\end{cor}

The following important lemma is a 
very special case of \cite[6.2. First Main Theorem]{sho-model}. 

\begin{lem}\label{thm-sho}
Let $(X, B)$ be a log canonical pair and 
let $\pi:X\to S$ be a projective morphism 
onto a variety $S$. 
We take $\sum _k D_k$ such that 
$\Supp B\subset \sum _k D_k$, where 
$D_i$ is an irreducible Weil divisor for every $i$ and 
$D_i\ne D_j$ for every $i\ne j$. 
We put 
\begin{align*}
\mathcal P=\left\{ \sum _k d_k D_k \ ;\  0\leq d_k\leq 1 \ {\text{for 
all}} \ k\  {\text{and}}\  K_X+\sum _k d_k D_k\  {\text{is lc}} 
\right\}. 
\end{align*}
Then $\mathcal P$ is a closed convex rational 
polytope. 

Let $\{R_j\}$ be any set of $(K_X+B)$-negative extremal 
rays of the lc pair $(X, B)$ over $S$. 
We put 
\begin{align*}
\mathcal N=\bigcap_j \left\{
\sum _k d_k D_k\in \mathcal P; (K_X+\sum _k d_k D_k)\cdot R_j\geq 0 
\right\}. 
\end{align*}
Then $\mathcal N$ is a closed convex subset of $\mathcal P$. 

We take $B'\in \mathcal P$. 
Let $\mathcal F$ be the minimal face of $\mathcal P$ containing $B'$. 
Assume that $(K_X+B')\cdot R_j >0$ for 
every $j$. 
Then there is an open subset $U$ of $\mathcal F$ in the 
classical topology such that 
$B'\in U\subset \mathcal N\cap \mathcal F$. 
In particular, 
we can write 
$$
K_X+B'=\sum _{i=1}^{d+1}r'_i (K_X+B'_i)
$$ 
with the following properties. 
\begin{itemize}
\item[(a)] $d=\dim \mathcal F$. 
\item[(b)] $B'_i \in \mathcal F$ for every $i$. 
\item[(c)] $m'(K_X+B'_i)$ is Cartier for some positive integer 
$m'$ for every $i$. 
\item[(d)] $\sum _{i=1}^{d+1}r'_i=1$ and 
$0\leq r'_i\leq 1$ for 
every $i$. 
\item[(e)] $(K_X+B'_i)\cdot R_j>0$ for every $i$ and $j$. 
\end{itemize}
\end{lem}

\begin{proof}
It is obvious that $\mathcal P$ is a closed convex rational 
polytope (see the proof of 
Lemma \ref{lem145}). 
By the definition, 
$\mathcal N$ is a closed convex subset of $\mathcal P$. 
Since $\mathcal F$ is a face of $\mathcal P$ and 
contains $B'$, we can take  
a $d$-dimensional rational simplex spanned by $\Delta_i$ for 
$1\leq i\leq d+1$ in 
$\mathcal F$ containing $B'$ inside it. 
Thus, we can write 
$$
K_X+B'=\sum _{i=1}^{d+1}r_i (K_X+\Delta_i) 
$$ 
such that $\sum _{i=1}^{d+1}r_i =1$ and $0<r_i <1$ for every $i$, 
and $m(K_X+\Delta_i)$ is Cartier for every $i$, where 
$m$ is a positive integer. 

We take an extremal ray $R_j$. 
By Corollary \ref{cor173}, 
we can find a curve $C_j$ on $X$ such that 
$C_j$ spans $R_j$ and 
that $m(K_X+\Delta_i)\cdot C_j =n_{ij}$ with 
$n_{ij}\geq -2m\dim X$ for every $i$. 
By the assumption, we have 
$$
(K_X+B')\cdot C_j=\sum _i \frac{r_in_{ij}}{m}>0. 
$$ 
We define 
$$
\alpha=\inf \left\{\sum _{i}\frac{r_i n_i}{m}
>0\  ;  n_i \geq -2m\dim X\ \text{and}\  
n_i \in \mathbb Z \ \text{for every} \ i \right\}.  
$$
Then we obtain $\alpha>0$. 
We put 
$$
c=\frac{\alpha}{2\dim X+\alpha +1}>0. 
$$ 
It is obvious that 
$$
B'+c(\Delta_i-B')\in \mathcal F 
$$ for every $i$ since $0<c<1$ and that 
$$
(K_X+B'+c(\Delta_i-B'))\cdot C_j>0 
$$ for 
every $i$ and $j$ by the definition of $c$. 
Thus, the $d$-dimensional simplex spanned 
by $B'+c(\Delta_i-B')$ for $1\leq i\leq d+1$ is contained in $\mathcal 
N\cap \mathcal F$ and contains 
$B'$ in its interior. 
So, the interior of the above simplex is a desired open 
set contained in $\mathcal N\cap \mathcal F$. 
Thus, we can write $$K_X+B'=\sum _{i=1}^{d+1}r'_i (K_X+B'_i)$$ with 
the required properties.   
\end{proof}

\begin{rem}
In \cite[6.2. First Main Theorem]{sho-model}, 
it is proved that $\mathcal N$ is a closed convex 
rational polytope. 
We recommend the reader to see \cite[Section 3]{birkar2} for details. 
The arguments in \cite[Section 3]{birkar2} work for lc pairs 
by Theorem \ref{prop146} (see, for example, \cite{book}). 
\end{rem}

By Corollary \ref{cor173} 
and Lemma \ref{thm-sho}, 
Lemma 2.6 in \cite{birkar} holds for lc pairs. 
It may be useful for the log minimal model 
program with scaling. 
We follow Birkar's proof in \cite{bp}. 

\begin{thm}[{cf.~\cite[Lemma 2.6]{birkar}}]\label{bir-prop} 
Let $(X,B)$ be an lc pair, let $B$ be an $\mathbb R$-divisor, 
and let $\pi:X\to S$ be a projective morphism between algebraic 
varieties. Let $H$ be an effective $\mathbb R$-Cartier 
$\mathbb R$-divisor on $X$ such that 
$K_X+B+H$ is $\pi$-nef and $(X, B+H)$ is lc. 
Then, either $K_X+B$ is also $\pi$-nef or there 
is a $(K_X+B)$-negative 
extremal ray $R$ such that $(K_X+B+\lambda H)\cdot R=0$, 
where 
$$\lambda:=\inf\{ t\geq 0 \, |\,  K_X+B+tH \ {\text{is $\pi$-nef}} \,\}. 
$$ 
Of course, $K_X+B+\lambda H$ is $\pi$-nef.  
\end{thm}
\begin{proof}
Assume that $K_X+B$ is not $\pi$-nef. 
Let $\{R_j\}$ be the set of $(K_X+B)$-negative extremal rays 
over $S$. 
Let $C_j$ be the rational 
curve spanning $R_j$ with the estimate as in 
Corollary \ref{cor173} for every $j$. 
We put 
$\mu=\underset{j}{\sup} \{\mu_j\}$, 
where 
$$
\mu_j =\frac{-(K_X+B)\cdot C_j}{H\cdot C_j}. 
$$ 
Obviously, $\lambda=\mu$ and 
$0<\mu\leq 1$. 
So, it is sufficient to 
prove that $\mu=\mu_l$ for some $l$. 
By Corollary \ref{cor173}, 
there are positive real numbers $r_1, \cdots, r_l$ and a positive integer 
$m$, which are independent of $j$, such that 
$$
-(K_X+B)\cdot C_j =\sum _{i=1}^{l}\frac{r_in_{ij}}{m}>0, 
$$ 
where $n_{ij}$ is an integer with $n_{ij}\leq 2m\dim X$ for 
every $i$ and $j$. 
If $(K_X+B+H)\cdot R_l=0$ for some 
$l$, then there are nothing to 
prove since $\lambda=1$ and 
$(K_X+B+H)\cdot R=0$ with $R=R_l$. 
Thus, we assume that $(K_X+B+H)\cdot R_j>0$ for every $j$. 
Therefore, we can 
apply Lemma \ref{thm-sho} and 
obtain 
$$
K_X+B+H=\sum _{p=1}^{q}r'_p(K_X+\Delta_p), 
$$ 
where $r'_1, \cdots, r'_q$ are positive real numbers, $(X, \Delta_p)$ is lc for every $p$, 
$m'(K_X+\Delta_p)$ is Cartier for some positive 
integer $m'$ and every $p$, 
and $(K_X+\Delta_p)\cdot C_j> 0$ for every $p$ and $j$. 
So, 
we obtain 
$$
(K_X+B+H)\cdot C_j=\sum _{p=1}^{q}\frac{r'_p n'_{pj}}{m'} 
$$ 
with $0< n'_{pj}=m'(K_X+\Delta_p)\cdot C_j\in \mathbb Z$.  
Note that $m'$ and $r'_p$ are independent of $j$ for 
every $p$. 
We also note that 
\begin{align*}
\frac{1}{\mu_j}=\frac{H\cdot C_j}{-(K_X+B)\cdot C_j}&=
\frac{(K_X+B+H)\cdot C_j}{-(K_X+B)\cdot C_j}+1
\\ &
=\frac{m\sum _{p=1}^q r'_p n'_{pj}}{m'\sum _{i=1}^lr_j n_{ij}}+1. 
\end{align*} 
Since 
$$\sum _{i=1}^{l}\frac{r_i n_{ij}}{m}>0 
$$ for every $j$ and $n_{ij}\leq 2m\dim X$ with $n_{ij}\in \mathbb Z$ for 
every $i$ and $j$, 
the number of the set $\{n_{ij}\}_{i, j}$ is finite. 
Thus, 
$$\inf _j \left\{\frac{1}{\mu_j}\right\}=\frac{1}{\mu_l}$$ 
for some $l$. Therefore, we obtain $\mu=\mu_l$. 
We finish the proof. 
\end{proof}

The following picture helps the reader 
to understand Theorem \ref{bir-prop}. 

\vspace{5mm} 
 
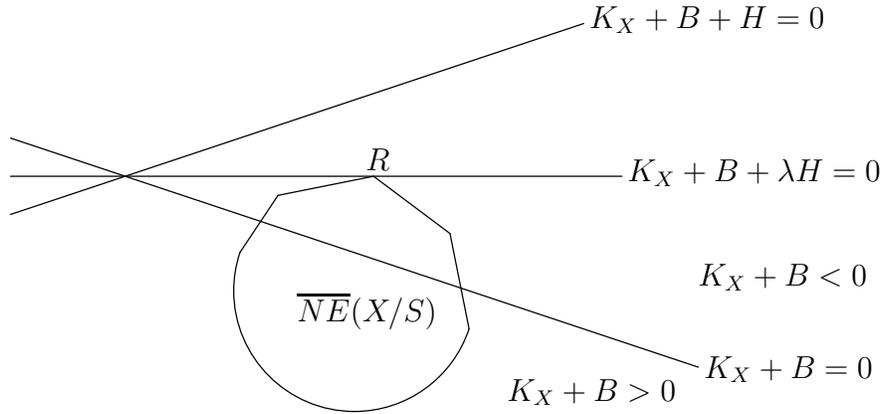
\begin{center}
\unitlength 0.1in
\begin{picture}( 46.0000, 21.5200)(  2.0000,-24.3200)
%
\special{pn 8}%
\special{pa 200 1200}%
\special{pa 3400 1200}%
\special{fp}%
\special{pa 200 1400}%
\special{pa 3200 400}%
\special{fp}%
\special{pa 200 1000}%
\special{pa 3800 2200}%
\special{fp}%
%
\special{pn 8}%
\special{ar 2000 1800 632 632  0.3217506 3.4633432}%
%
\special{pn 8}%
\special{pa 2600 2000}%
\special{pa 2500 1500}%
\special{fp}%
\special{pa 2500 1500}%
\special{pa 2100 1200}%
\special{fp}%
\special{pa 2100 1200}%
\special{pa 1600 1300}%
\special{fp}%
\special{pa 1600 1300}%
\special{pa 1400 1600}%
\special{fp}%
\put(17.0000,-20.0000){\makebox(0,0)[lb]{$\overline{NE}(X/S)$}}%
\put(20.6000,-11.7000){\makebox(0,0)[lb]{$R$}}%
\put(32.3000,-4.5000){\makebox(0,0)[lb]{$K_{X}+B+H=0$}}%
\put(34.3000,-12.5000){\makebox(0,0)[lb]{$K_{X}+B+\lambda H=0$}}%
\put(38.0000,-18.0000){\makebox(0,0)[lb]{$K_{X}+B<0$}}%
\put(38.3000,-22.8000){\makebox(0,0)[lb]{$K_{X}+B=0$}}%
%
\put(48.0000,-10.0000){\makebox(0,0)[lb]{}}%
\put(28.0000,-24.0000){\makebox(0,0)[lb]{$K_{X}+B>0$}}%
\end{picture}%
\end{center}
\vspace{5mm} 

The main result of this section is an estimate 
of lengths of extremal rays which are relatively ample at non-lc loci 
(cf.~\cite{kollar2}, \cite{kollar3}). 

\begin{thm}\label{thm-la} 
Let $X$ be a normal variety, 
let $B$ be an effective $\mathbb R$-divisor on $X$ such that 
$K_X+B$ is $\mathbb R$-Cartier, and let 
$\pi:X\to S$ be a projective morphism onto a variety $S$. 
Let $R$ be a $(K_X+B)$-negative extremal ray of 
$\overline {NE}(X/S)$ which is relatively ample at $\Nlc(X, B)$. 
Then we can find a rational curve $C$ on $X$ such that 
$[C]\in R$ and $$0<-(K_X+B)\cdot C\leq 2\dim X.$$ 
\end{thm}
\begin{proof}
By shrinking $S$, we can assume that $S$ is quasi-projective. 
By replacing $\pi:X\to S$ with the extremal contraction 
$\varphi_R:X\to Y$ over $S$ (cf.~Theorem \ref{thm144} (3)), 
we can assume that 
the relative Picard number $\rho (X/S)=1$ and that 
$\pi$ is an isomorphism in a neighborhood of $\Nlc (X, B)$. In particular, 
$-(K_X+B)$ is $\pi$-ample. 
By Theorem \ref{thm91}, there is a projective birational 
morphism $f:Y\to X$ such that 
\begin{itemize}
\item[(i)] $K_Y+B_Y=f^*(K_X+B)+\underset{a(E, X, B)<-1}{\sum} (a(E, X, B)+1)E$ 
where $B_Y=f^{-1}_*B+\underset{E:{\text{$f$-exceptional}}}\sum E$,  
\item[(ii)] $(Y, B_Y)$ is a $\mathbb Q$-factorial dlt pair, and 
\item[(iii)] $D=B_Y+F$ with $F=-\underset{a(E, X, B)<-1}{\sum}{(a(E, X, B)+1)E\geq 0}$. 
\end{itemize}
We note that $K_Y+D=f^*(K_X+B)$. 
Therefore, we have 
$$f_*(\overline{NE}(Y/S)_{K_Y+D\geq 0})\subseteq 
\overline {NE}(X/S)_{K_X+B\geq 0}=\{0\}.$$ 
We also note that 
$$
f_*(\overline {NE}(Y/S)_{\Nlc (Y, D)})=\{0\}. 
$$ 
Thus, there is a $(K_Y+D)$-negative extremal ray $R'$ of $\overline {NE}(Y/S)$ which 
is relatively ample at $\Nlc(Y, D)$. 
By Theorem \ref{thm144} (1), $R'$ is spanned by a curve $C^{\dag}$. 
Since $-(K_Y+D)\cdot C^{\dag}>0$, we see that 
$f(C^{\dag})$ is a curve. 
If $C^{\dag}\subset \Supp F$, 
then $f(C^{\dag})\subset \Nlc (X, B)$. 
It is a contradiction because $\pi\circ f(C^{\dag})$ is a point. 
Thus, $C^{\dag}\not\subset \Supp F$. 
Since $-(K_Y+B_Y)=-(K_Y+D)+F$, 
we can see that $R'$ is a $(K_Y+B_Y)$-negative extremal ray of $\overline {NE}(Y/S)$. 
Therefore, we can find a rational curve $C'$ on $Y$ such that 
$C'$ spans $R'$ and that 
$$
0<-(K_Y+B_Y)\cdot C'\leq 2\dim X
$$ 
by Theorem \ref{prop146}. 
By the above argument, we can easily see that 
$C'\not\subset \Supp F$. 
Therefore, we obtain 
\begin{align*}
0<-(K_Y+D)\cdot C'&=-(K_Y+B_Y)\cdot C'-F\cdot C'\\
&\leq -(K_Y+B_Y)\cdot C'\leq 2\dim X. 
\end{align*}
Since $K_Y+D=f^*(K_X+B)$, $C=f(C')$ is a rational 
curve on $X$ such that 
$\pi(C)$ is a point and $0<-(K_X+B)\cdot C\leq 2\dim X$. 
\end{proof}

\begin{rem}
In Theorem \ref{thm-la}, we can easily prove 
$0<-(K_X+B)\cdot C\leq \dim X+1$ when $\dim X\leq 2$. For 
details, see \cite[Proposition 3.7]{fujino16}. 
\end{rem}

\section{Ambro's theory of quasi-log varieties}\label{sec16}

In this section, we make some comments on Ambro's theory of 
quasi-log varieties. We strongly recommend the reader to see \cite{fuji-lec} 
for an introduction to the theory of quasi-log varieties. 

In the acknowledgements in \cite{ambro}, 
Ambro wrote \lq\lq The motivation behind 
this work is his (Professor Shokurov's) idea that log 
varieties and their LCS loci should be treated on an equal 
footing.\rq\rq  \ 
So, in the theory of quasi-log varieties, we have to treat highly 
reducible non-equidimensional 
varieties (see Example \ref{ex-18} below). 
Therefore, our approach explained 
in this paper is completely different from 
the theory of quasi-log varieties. 
We recommend the reader to compare our proof of the 
base point free theorem for projective lc surfaces in Section \ref{sec2}
with Ambro's proof 
(see, for example, \cite[Section 4]{fuji-lec}). 

Let us explain some results of the theory of quasi-log varieties 
which can not be covered by our approach. 

\begin{say}
Let $(X, B)$ be a projective log canonical pair and 
let $\{C_i\}$ be any set of lc centers of the pair 
$(X, B)$. 
We put $W=\bigcup C_i$ with the reduced scheme structure. 
Then $[W, \omega]$ is a {\em{qlc pair}}, where 
$\omega=(K_X+B)|_W$. For the definition of 
{\em{qlc pairs}}, see 
\cite[Definition 3.29]{book} or \cite[Definition 3.1]{fuji-lec}. 

\begin{ex}\label{ex-18}
Let $V$ be a projective toric variety and let 
$D$ be the complement of the big torus. 
Then $(V, D)$ is log canonical 
and $K_V+D\sim 0$. 
In this case, every torus invariant closed 
subvariety $W$ of $V$ with 
$\omega=0$ is a qlc pair. 
In particular, $W$ is not necessarily pure-dimensional 
(cf.~\cite[\S 5]{fuji-poly}). 
\end{ex}

We can prove the cone theorem for $[W, \omega]$. 

\begin{thm}[Cone theorem]
We have 
$$
\overline {NE}(W)=\overline {NE}(W)_{\omega\geq 0}+\sum _j R_j. 
$$ 
\end{thm}
For the details, see \cite[3.3.3 Cone Theorem]{book}. 
We can also prove the base point free theorem. 

\begin{thm}[Base point free theorem] 
Let $L$ be a nef Cartier divisor on $W$ such that 
$aL-\omega$ is ample for some 
$a>0$. 
Then $|mL|$ is base point free for $m \gg 0$. 
\end{thm} 
See, for example, \cite[3.3.1 Base point free theorem]{book}. 
By these theorems, we have the following statement. 

\begin{thm}[Contraction theorem] 
Let $F$ be an $\omega$-negative extremal face of $\overline {NE}(W)$. 
Then there is a contraction morphism $\varphi_F:W\to V$ with the 
following properties. 
\begin{itemize}
\item[(i)] Let $C$ be an integral curve on $W$. Then 
$\varphi_F(C)$ is a point if and only if 
$[C]\in F$. 
\item[(ii)] $\mathcal O_V\simeq (\varphi_F)_*\mathcal O_W$. 
\item[(iii)] Let $L$ be a line bundle 
on $W$ such that $L\cdot C=0$ for every 
curve $C$ with 
$[C]\in F$. 
Then there is a line bundle 
$L_V$ on $V$ such that $L\simeq \varphi^*_FL_V$. 
\end{itemize}
\end{thm}
\end{say}
For the details of the theory of quasi-log varieties, 
see \cite{book}. The book \cite{book} treats 
some various other topics 
which can not be covered by this paper. 

\section{Related topics}\label{sec19} 

In this final section, we briefly explain some related topics obtained 
by the author for the reader's convenience.

In this paper, we did not describe 
the notion of singularities of pairs. However, 
it is very important when we read some papers on the log minimal model program. 
We think that \cite{what} helps the reader to 
understand the subtlety of the notion of dlt pairs. 

The reader can find that all the injectivity, vanishing, and 
torsion-free theorems in this paper 
are discussed in full generality in \cite[Sections 2 and 3]{book}. 
They heavily depend on the theory of mixed Hodge structures 
on compact support cohomology groups of reducible varieties. 

We omitted the explanation of the log minimal model 
program for log canonical pairs. 
It is because the framework is the same as 
for klt pairs. The reader can find it in \cite[Section 3]{book}. 
We note that the existence problem of log canonical flips 
is still open in dimension $\geq 5$ and 
the termination of log canonical flips follows from 
the termination of klt flips. For the details, see \cite[Section 3]{book}. 

In \cite{fuji-I}, we prove an effective version of the base point free 
theorem for log canonical pairs. 
It is a log canonical version of Koll\'ar's effective freeness. 
In \cite{fuji-II}, the Angehrn--Siu type effective 
base point free theorems 
are proved for log canonical pairs. 
The reader can find that the proof of our non-vanishing theorem 
(cf.~Theorem \ref{thm111} and \cite[Theorem 1.1]{non-va}) 
grew out from the arguments in \cite{fuji-I} and \cite{fuji-II}. 

In \cite{non-lc}, we systematically treat the basic properties of 
non-lc ideal sheaves, especially, 
the restriction theorem of 
non-lc ideal sheaves for normal divisors. 
It is a generalization of 
Kawakita's inversion of adjunction on log canonicity. 
See also \cite{ft} for further discussions on various 
analogues of non-lc ideal sheaves.  
 
In \cite{fuji-finite}, 
we prove the finite generation of the log canonical 
ring for log canonical pairs in dimension four and discuss related topics. 
It induces the existence of fourfold 
log canonical flips. 

In \cite{fujino16}, we discuss the minimal model theory for log surfaces. 
The results in \cite{fujino16} 
are obtained under much weaker assumptions than everybody expected. 
The paper \cite{fujino16} is an ultimate application of our new approach to the log minimal 
model program. 
\ifx\undefined\bysame
\newcommand{\bysame|{leavemode\hbox to3em{\hrulefill}\,}
\fi

\end{document}